\documentclass[10pt]{article}
\usepackage{amsfonts}
\usepackage{amssymb,amsmath,amsthm, amsfonts}

\def\sqr#1#2{{\vcenter{\vbox{\hrule height.#2pt
              \hbox{\vrule width.#2pt height#1pt \kern#1pt \vrule width.#2pt}
          \hrule height.#2pt}}}}
\def\sqr#1#2{{\vcenter{\vbox{\hrule height.#2pt
              \hbox{\vrule width.#2pt height#1pt \kern#1pt \vrule width.#2pt}
              \hrule height.#2pt}}}}
\def\3n{\negthinspace \negthinspace \negthinspace }
\def\2n{\negthinspace \negthinspace }
\def\1n{\negthinspace }

\def\={\buildrel \triangle \over =}

\def\min{\mathop{\rm min}}

\def\sup{\mathop{\rm sup}}

\def\inf{\hbox{\rm inf$\,$}}

\def\|{\Big |}
\def\({\Big (}
\def\){\Big )}
\def\[{\Big[}
\def\]{\Big]}
\def\be{\begin{equation}}
\def\bel{\begin{equation}\label}
\def\ee{\end{equation}}
\def\bt{\begin{theorem}}
\def\bcd{\begin{condition}}
\def\ecd{\end{condition}}
\def\et{\end{theorem}}
\def\bc{\begin{corollary}}
\def\ec{\end{corollary}}
\def\bde{\begin{definition}}
\def\ede{\end{definition}}
\def\bl{\begin{lemma}}
\def\el{\end{lemma}}
\def\bp{\begin{proposition}}
\def\ep{\end{proposition}}
\def\bex{\begin{example}}
\def\eex{\end{example}}
\def\br{\begin{remark}}
\def\er{\end{remark}}
\def\ba{\begin{array}}
\def\ea{\end{array}}
\def\ed{\end{document}}

\def\square#1{\vbox{\hrule\hbox{\vrule height#1%
     \kern#1\vrule}\hrule}}
\def\rectangle#1#2{\vbox{\hrule\hbox{\vrule height#1%
     \kern#2\vrule}\hrule}}

\font\tenbb=msbm10 \font\sevenbb=msbm7 \font\fivebb=msbm5

\newfam\bbfam
\scriptscriptfont\bbfam=\fivebb \textfont\bbfam=\tenbb
\scriptfont\bbfam=\sevenbb

\newtheorem{lemma}{Lemma}[section]
\newtheorem{remark}{Remark}[section]
\newtheorem{example}{Example}[section]
\newtheorem{theorem}{Theorem}[section]
\newtheorem{corollary}{Corollary}[section]

\newtheorem{definition}{Definition}[section]
\newtheorem{proposition}{Proposition}[section]
\newtheorem{condition}{Condition}[section]

\makeatletter
   
   \@addtoreset{equation}{section}
\makeatother

\begin{document}
\title{Nonlinear stochastic differential games involving a major player and
a large number of collectively acting minor agents\footnote{The work has
been supported by the NSF of P.R.China (Nos. 11071144, 11171187, 11222110), Shandong Province (Nos. BS2011SF010, JQ201202), Program for New Century Excellent Talents in University (NCET, 2012), 111 Project (No. B12023). Juan Li is the corresponding author.}}
\author{Rainer Buckdahn$^{1, 3}$,\, \,  Juan Li$^{2}$,\, \,  Shige Peng$^{3}$\\
{\small $^1$Laboratoire de Math\'{e}matiques LMBA, CNRS-UMR 6205, Universit\'{e} de
Bretagne Occidentale,}\\
 {\small 6, avenue Victor-le-Gorgeu, CS 93837, 29238 Brest cedex 3, France.}\\
{\small $^2$School of Mathematics and Statistics, Shandong University (Weihai), Weihai 264209, P. R. China.;}\\
{\small $^3$School of Mathematics, Shandong University, Jinan 250100, P. R. China.}\\
{\small{\it E-mails: rainer.buckdahn@univ-brest.fr, juanli@sdu.edu.cn, peng@sdu.edu.cn.}}
\date{August 15, 2013} }\maketitle \noindent{\bf
Abstract}\hskip4mm The purpose of this paper is to study 2-person
zero-sum stochastic differential games, in which one player is a
major one and the other player is a group of $N$ minor agents which
are collectively playing, statistically identical and have the same
cost-functional. The game is studied in a weak formulation; this
means in particular, we can study it as a game of the type
``feedback control against feedback control''.  The payoff/cost
functional is defined through a controlled backward stochastic
differential equation, for which driving coefficient is assumed to satisfy strict concavity-convexity with respect to the control
parameters. This ensures the existence of saddle point feedback
controls for the game with $N$ minor agents. We study the limit
behavior of these saddle point controls and of the associated
Hamiltonian, and we characterize the limit of the saddle point
controls as the unique saddle point control of the limit mean-field stochastic differential game.


\section{\large{Introduction}}

\ \ \ \ \ \ \ In this paper we study a particular type of 2-person
zero-sum stochastic differential games, where one
player is a major one, who plays against a group of $N$ collectively
acting minor agents which each of them participate with the same
percentage and are statistically identical. We study the stochastic
differential game in its weak formulation and with a
pay-off/cost-functional given through a backward stochastic
differential equation (BSDE), which allows to consider the game of
the type ``feedback control against feedback control''. Under
suitable assumptions on the driving coefficient of the BSDE we show
for the game with $N$ minor agents the existence of  saddle point
feedback controls, which can be characterized as Stackelberg
feedback strategy, where the major player is the leader and the
collectively acting minor agents are the follower. We investigate
the limit of these saddle point feedback controls and of the
associated Hamiltonian of stochastic differential game and
characterize the limit saddle point controls as unique saddle point
controls of the limit stochastic differential game which turns out
to be of mean-field type.

Mean-field stochastic differential games obtained as limit of
stochastic differential games between $N$ statistically identical
players were studied by Lasry and Lions in a series of pioneering
papers (\cite{LL1}, \cite{LL2}, \cite{LL3}, \cite{LL4}). In their papers they
investigated so-called approximate Nash equilibria, obtained as
distributed closed-loop strategies given by solving the limit
problem. This limit problem consists of a coupled system, formed by
a Hamilton-Jacobi-Bellman equation and an equation of Kolmogorov
type, the first one with a terminal and the second one with an
initial condition. In subsequent works several authors studied
different applications coming from such different domains as, for
instance, Statistical Mechanics and Physics, finance and management
of exhaustible resources. On the other hand, motivated by problems
occurring in large communication networks but also in stochastic
differential games involving a large number of players, Huang,
Malham\'e and Caines \cite{HMC} introduced a similar concept, that
of Nash Certainty Equivalence.

The problem of investigating stochastic differential games with a
large number of players or with one major player at one side and $N$
minor players on the other side, which participate each of them in
the same proportion and in a symmetric way at the game, leads in the
limit, as $N$ tends to infinity, to an averaging over the minor
players and, thus, to a mean-field limit game. Such a problem,
without the control part, relates with the McKean-Vlasov theory of
chaos propagation (we refer the reader to the scholarly paper by
Sznitman \cite{S}). In \cite{CD} and \cite{CDL} Carmona and Delarue,
and Carmona, Delarue and Lachapelle, respectively, study approximate
Nash-equilibria for $N$-person non-zero stochastic games with
mean-field interaction between the players and they discuss the
limit behavior as the number of players converges to infinity. In
\cite{CD} the authors embed the Mean-Field Game strategy developed
by Lasry and Lions, in an analytical approach into a purely
probabilistic framework, which transforms the strongly coupled
system consisting of a Hamilton-Jacobi-Bellman equation with
terminal condition and a Kolmogorov-type equation with initial
condition, considered by Lasry and Lions, into a strongly coupled
forward-backward stochastic differential equation (FBSDE) of
McKean-Vlasov type. The authors of \cite{CD} show that the solutions
of this FBSDE together with the associated FBSDE value function
allow to obtain a set of distributed strategies which turn out to be
an $\varepsilon_N$-approximate Nash equilibrium for the $N$-person
non-zero sum stochastic differential game, where $\varepsilon$
converges to zero as $N$ tends to infinity. Such a kind of argument,
but for simpler models can be also found in \cite{BSYY} an in
\cite{C}. In \cite{CDL} Carmona, Delarue and Lachapelle make
comparing studies between stochastic differential games with
mean-field interactions on one side and the characterization of
optimal strategies for the associated mean-field linear-quadratic
McKean-Vlasov stochastic control problem on the other side.

In their recent work \cite{NC} Nourian and Caines study with a different approach that those chosen by Lasry and Lions and Carmona
and Delarue $\varepsilon$-Mean-Field games. They consider
$N+1$-person non-zero sum stochastic differential games with one
major player and $N$ symmetric minor players. The stochastic
dynamics are non-linear and with mean-field interaction, the forward
equation for each player is driven by its own Brownian motion and
each player control only his own dynamics and his own running cost;
the running cost are non-linear and with mean-field interaction. The
study of $\varepsilon_N$-Nash equilibriums for the $N+1$-person game
for large $N$ leads the authors to a strongly coupled stochastic
mean-field system composed of a stochastic Hamilton-Jacobi-Bellman
equation with terminal condition and two McKean-Vlasov equations
with stochastic coefficients, describing the state of the major
player as well as the measure determining the mean-field behavior
of the minor agents. Let us emphasize that the paper \cite{NC}
represent an extension to the framework of non-linear mean-field
stochastic differential equations, which was preceded by papers by
Huang \cite{H} but also by Huang together with Caines and Malham\'e
\cite{HMC} and with Nourian \cite{NCMH}, in order to mention only
these important works of a longer list of papers.

In the present work we study a somehow different framework which,
although is related with the works discussed above, in particular
with \cite{NC}. But unlike \cite{NC} we consider the $N$ minor
agents as collectively acting, with a common cost functional. This
allows to consider the game as 2-person zero-sum stochastic
differential game. Moreover, it will be studied in a weak form.
Stochastic differential games in the weak form have been studied by
Hamad\`{e}ne and Lepeltier \cite{HaL} but also by Hamad\`{e}ne in
different works, see, e.g., \cite{Ha1} and \cite{Ha2}.

In order to be more precise, for independent Brownian motions
$\tilde{W}^0, \tilde{W}^1,\dots, \tilde{W}^N$ and given initial
positions $x^{(N)}=(x_0,x_1,\dots,x_N)$ we consider the weak solution
$X^{(N)}=(X^{0,N},X^{1,N},\dots,X^{N,N})$ of the following system of
dynamics with given feedback controls $u=u(X^{(N)})$ - for the major
player and $v^{(N)}(=(v^1,\dots,v^N))=v(X^{(N)})$ for the $N$ collectively
acting minor agents:
\be\label{SDE_iaa}\begin{array}{llll} dX_s^{0,N}&=&\displaystyle
\sigma_0(X^{0,N})d\tilde{W}^0_s+\frac{1}{N}
\sum_{\ell=1}^Nb_0(X^{0,N}_s,X^{\ell,N}_s,Z^{0,N}_s)u_sds,\, &
X_t^{0,N}=x_0,\\
dX_s^{j,N}&=&\displaystyle
\sigma_1(X^{j,N})d\tilde{W}^j_s+\frac{\varepsilon_N}
{N}\sum_{\ell=1}^Nb_1(X^{0,N}_s,X^{\ell,N}_s, Z^{j,N}_s)v_s^\ell
ds,\, & X_t^{j,N}=x_j,
\end{array}\ee

\noindent $1\le j\le N,$ and we associate the nonlinear payoff/cost functional
$J(t,x;u,v):=Y^N_t$ defined through the following BSDE:
\be\label{BSDE_iaa}\begin{array}{rll} dY^N_s&=&
-\displaystyle\frac{1}{N}\sum_{\ell=1}^Nf(X_s^{0,N},
X_s^{\ell,N},Y^N_s, Z_s^{0,N},Z_s^{\ell,N},u_s,v^\ell_s)ds +
\sum_{\ell=0}^NZ_s^{\ell,N}d\tilde{W}_s^{\ell},\, s\in[t,T],\\
Y^N_T&=& \displaystyle\frac{1}{N}\sum_{\ell=1}^N\Phi(X_T^{0,N},
X_T^{\ell,N}).
\end{array}\ee
We see that here the major player can control his own dynamics
$X^{0,N}$, and each of the minor agents controls his dynamics but
also those of all the other minor agents and, together with the
major player; moreover, all players can control the
pay-off/cost-functional $J(t,x;u,v)$. While the objective of the
major player is to maximize $J(t,x;u,v)$, the collectively acting
minor agents want to minimize their common cost functional
$J(t,x;u,v)$.

The choice of the weak formulation of the problem of stochastic
differential games allows to shift with the help of a Girsanov
transformation the doubly controlled drift terms of the dynamics of
the game with $N$ collectively acting minor players into the BSDE
defining the pay-off/cost functional. This reduces the study of the
limit behavior of the game and of the saddle point feedback
controls for this game to the investigation of the limit behavior of
the corresponding BSDEs with non-feedback controls and the limiting
Mean-Field BSDE. In order to guarantee the existence of saddle point
controls we impose on the driving coefficient $f$ of the BSDE a
strict concavity-convexity assumption with respect to the control
parameters $(u,v)$. The specificity of our approach using the
Girsanov transformation necessitates the factor $\varepsilon_N$
which is supposed to be of order $O(N^{-3/4})$, as $N\rightarrow
+\infty.$

The problem of a limit approach for mean-field BSDEs as well as
mean-field BSDEs themselves were studied by the authors in
\cite{BDLP} (together with Djehiche) and in \cite{BLP}, but unlike
here without controls.

The paper is organized as follows: In Section 2 a short recall for
the convergence of the above system (\ref{SDE_iaa}) and
(\ref{BSDE_iaa}) in the case without control will be given. In Section
3 the stochastic differential game with one major player and $N$
collectively acting minor agents is introduced, the assumptions on
the coefficients are given and the existence of saddle points
controls which are of feedback form is discussed. They are
characterized as a Stackelberg feedback strategy. This
characterization admits estimates for the saddle point feedback
controls which will be used in what follows. Section 5 is devoted to
the study of the limit mean-field game. For this the convergence of
the saddle point controls for the game with $N$ minor agents is
proved and the limit controls are shown to be the unique saddle
point controls of the limit game. In order to improve the
readability of the work, the proofs of several lemmas have been
postponed to the Appendix.

\section{Preliminaries. The N+1 players system without control}

\ \ \ \ \ \ \ \ \ In this short section we consider first briefly the case of a
``stochastic differential game with $N$ minor agents'' without
control and recall its limit behavior. We restrict for this to the
(particular) case we will need for our discussion in the Sections 3
and 4. The more interested reader is referred to \cite{S}, and for
the BSDE part, for instance, to \cite{BDLP}.

Let $(\Omega, {\cal F},P)$ be a complete probability space, endowed
with a sequence of independent $d$-dimensional Brownian motions
$W^j=(W^j_s)_{s\in[0,T]} , j\ge 0$, where $T>0$ is an arbitrarily
fixed time horizon. We denote by $\mathbb{F}=({\cal
F}_t)_{t\in[0,T]}$ the filtration generated by $W^j,\ j\ge 0$, and
augmented by the $P$-null sets. Given bounded Lipschitz coefficients
$b_{0}: {\mathbb R}^d\times {\mathbb R}^d\rightarrow {\mathbb R}^d,\
\sigma_0: {\mathbb R}^d\times {\mathbb R}^d\rightarrow {\mathbb
R}^{d\times d}$, $b_1: {\mathbb R}^d\times {\mathbb R}^d\times
{\mathbb R}^d\rightarrow {\mathbb R}^{d\times d}$ and $\sigma_1: {\mathbb R}^d\times {\mathbb R}^d\times {\mathbb R}^d\rightarrow
{\mathbb R}^{d\times d}$, we consider for an arbitrarily chosen
initial time $t\in [0,T]$ and initial positions $x_0,\dots, x_N\in
{\mathbb R}^d$ the system of $N\ge 2$\ $d$-dimensional coupled
stochastic differential equations (SDEs):
\be\label{3.1}\begin{array}{lll}
&dX_s^{0,N}=\frac{1}{N}\sum_{\ell=1}^Nb_0(X_s^{0,N},X_s^{\ell,N})
ds+\frac{1}{N}
\sum_{\ell=1}^N\sigma_0(X_s^{0,N},X_s^{\ell,N})dW_s^0,\, s\in
[t,T],\\
&X_t^{0,N}=x_0,\end{array}\ee
\be\label{3.2}\begin{array}{lll}
&dX_s^{j,N}=\frac{1}{N}\sum_{\ell=1}^Nb_1(X_s^{0,N},X_s^{j,N},
X_s^{\ell,N})ds+\frac{1}{N}
\sum_{\ell=1}^N\sigma_1(X_s^{0,N},X_s^{j,N},X_s^{\ell,N})dW_s^j,\
s\in [t,T],\\
&X_t^{j,N}=x_j,\ 1\le j\le N,\end{array}\ee \noindent associated
with the backward stochastic differential equation (BSDE):
\be\label{3.3}\begin{array}{lll}
&dY_s^{N}=-\frac{1}{N}\sum_{\ell=1}^Nf(X_s^{0,
N},X_s^{\ell,N},Y_s^{N},Z_s^{0,N},Z_s^{\ell,N})ds
+Z_s^{0,N}dW_s^{0}+\sum_{j=1}^NZ_s^{j,N}dW_s^{j},\ s\in[t,T],\\
&Y_T^{N}=\frac{1}{N} \sum_{\ell=1}^N\Phi(X_T^{0,
N},X_T^{\ell,N}),\end{array}\ee \noindent where the functions $f: {\mathbb R}^d\times {\mathbb R}^d\times
{\mathbb R}\times {\mathbb R}^d\times {\mathbb R}^d\rightarrow
{\mathbb R}$ and $\Phi: {\mathbb R}^d\times {\mathbb R}^d\rightarrow
{\mathbb R}$ are assumed to be bounded and Lipschitz in all its
variables.

Let us now discuss the forward equation (\ref{3.1})-(\ref{3.2}) and the
backward equation (\ref{3.3}) separately.

\subsection{Limit behavior of the forward stochastic system}

\ \ \ \ \ \ \ \ \ The objective of this section is to discuss briefly the limit
behavior of the system (\ref{3.1}) and (\ref{3.2}):
$$\label{3.1x}\begin{array}{lll}
&dX_s^{0,N}=\frac{1}{N}\sum_{\ell=1}^Nb_0(X_s^{0,N},X_s^{\ell,N})ds+
\frac{1}{N} \sum_{\ell=1}^N\sigma_0(X_s^{0,N},X_s^{\ell,N})dW_s^0,\,
s\in [t,T],\\
&X_t^{0,N}=x_0;\end{array}$$
$$\label{3.2x}\begin{array}{lll}
&dX_s^{j,N}=\frac{1}{N}\sum_{\ell=1}^Nb_1(X_s^{0,N},X_s^{j,N},
X_s^{\ell,N})ds+\frac{1}{N} \sum_{\ell=1}^N\sigma_1(X_s^{0,N},
X_s^{j,N},X_s^{\ell,N})dW_s^j,\ s\in [t,T],\\
& X_t^{j,N}=x_j, 1\le j\le N,\end{array}$$

\noindent as $N$ tends to $+\infty.$ The limit behavior of such
systems as well as the associated limit McKean-Vlasov SDEs have been
already largely discussed in the literature. For completeness we
first state the following classical existence and uniqueness result.

\bp\label{p1} Under our standard assumptions, that is, the coefficients are bounded and Lipschitz in their variables, we have for any initial datum
$(t, x^{(N)})\in[0,T]\times ({\mathbb R}^{d})^{N+1}$,\ $x^{(N)}=(x_0, x_1,\dots, x_N)$, the existence and the uniqueness of
the solution $X^{(N)}=(X^{0,N}, X^{1,N},\dots, X^{N,N})$\ in the space
${\cal S}_{\mathbb{F}}^2(t,T; {\mathbb R}^d)^{N+1}$, where ${\cal
S}_{\mathbb{F}}^2(t,T; {\mathbb R}^d)$ denotes the space of all
continuous $\mathbb{F}$-adapted, ${\mathbb R}^d$-valued processes
which supremum of the Euclidean norm over the interval $[t,T]$ is square integrable. \ep

Let us now suppose that $(x_j)_{j\ge 0}\subset {\mathbb R}^d$ is such that, for some
$\overline{x}\in {\mathbb R}^d$,
\be\label{3.4}\frac{1}{N}\sum_{j=1}^N|x_j-\overline{x}|^2\rightarrow
0,\quad\mbox{as }N\rightarrow +\infty.\ee
In order to justify the choice of this condition, we let ${\cal
P}_2({\mathbb R}^d)$ be the space of the Borel probability measures on ${\mathbb R}^d$ with
finite second moments. Wishing that the convergence of the
above system (\ref{3.1})-(\ref{3.2}) can be measured in terms of the
(Monge-Kantorovich-)Wasserstein distance $d_{2}$ of second order,
$$d_2(\mu,\nu):=\inf\{E[|\xi-\eta|^2]: \xi,\ \eta\in L^0({\cal F}; {\mathbb R}^d)
\ \mbox{with}\ P_\xi=\mu,\ P_\eta=\nu\},\quad \mu,\ \nu\in{\cal
P}_2({\mathbb R}^d).$$
We define $\nu^N=\frac{1}{N}\sum_{j=1}^N\delta_{x_j}$, with $\delta_{x_j}$
denoting the Dirac measure with mass in $x_j$. Given another
probability $\nu\in {\cal P}_2({\mathbb R}^d)$, the convergence
$d_2(\nu^N,\nu)\rightarrow 0,$ as $N\rightarrow +\infty$, is
equivalent to the weak convergence of $\nu^N$ to $\nu$ as well as
that of their second moments. Preferring, for simplicity, the choice
$\nu=\delta_{\overline{x}}$, for some $\overline{x}\in {\mathbb R}^d,$ the
convergence $d_2(\nu^N,\nu)\rightarrow 0$ is equivalent to
$(\ref{3.4})$.

The $L^2$-limit system of the above systems of SDEs
(\ref{3.1})-(\ref{3.2}) (see Proposition 2.2) is given by
\begin{equation}\label{aaa}\begin{array}{llll}
d\overline{X}_s^0 & = & b_0(\overline{X}_s^0,{\mu}_s)ds+\sigma_0
(\overline{X}_s^0,{\mu}_s)dW_s^0, & s\in[t,T],\ \overline{X}_t^0=x_0,\\
d\overline{X}_s^j & = & b_1(\overline{X}_s^0,\overline{X}_s^j,{\mu}_s)ds
+\sigma_1(\overline{X}_s^0,\overline{X}_s^j,{\mu}_s)dW_s^j, & s\in[t,T],
\ \overline{X}_t^j=\overline{x},\ j\ge 1,\\
\end{array}\end{equation}
where
\begin{equation}\label{aaaz}\mu_s(dy)=P\{\overline{X}_s^1\in dy\ |
\ {\cal F}_T^{W^0}\}
\end{equation}

\noindent is the conditional distribution law of $\overline{X}_s^1$
knowing the $\sigma$-field ${\cal F}_T^{W^0}$ which is generated by $W^0$
over the time interval $[0,T]$ and augmented by all $P$-null sets.

Given a bounded measurable function $h$ over ${\mathbb R}^d$ we use the
notation $h(\mu_s)=\int_{{\mathbb R}^d}h(x)\mu_s(dx).$ By $L$ we denote the
second order operator \be\label{3.63}
L[\mu, x]\varphi(x'):=\frac{1}{2}tr\left(\sigma_1\sigma_1^*(x, x', \mu)
D^2\varphi(x')\right)+b_1(x, x', \mu) D\varphi(x'),\ x,\ x'\in {\mathbb
R}^d,\ee defined for probability measures $\mu$ on ${\mathbb R}^d$
and functions $\varphi\in C^2_K({\mathbb R}^d)$, and by $L(x,\mu)^*$
we denote its dual operator applying to the probability measures on
${\mathbb R}^d$. It is well known that $\mu=(\mu_s)_{s\in[t,T])}$ can be
characterized as the weak solution of the PDE with stochastic
coefficients \be\label{3.65}\displaystyle\frac{d}{ds}\mu_s=L[\mu_s, \overline{X}_s^0]^*\mu_s,\ s\in(t,T],\ \mu_t(dy)=\delta_{\overline{x}}(dy),\ee i.e.,
\be\label{3.64}\displaystyle\frac{d}{ds}\langle\mu_s, \varphi\rangle=
\langle\mu_s, L[\mu_s, \overline{X}_s^0]\varphi \rangle=\langle
L[\mu_s,\overline{X}_s^0]^*\mu_s, \varphi \rangle,\ \mbox{ for all}\ \varphi\in C^2_K({\mathbb R}^d)\ee ($C^2_K({\mathbb R}^d)$ denotes
the space of $C^2$-functions with compact support in ${\mathbb R}^d$).

Let us first remark the following proposition.

\bp\label{p2} Under our standard assumptions the above system $(\ref{aaa})$ has a unique
solution $\overline{X}^j=(\overline{X}^j_s)_{s\in[t,T]}\in {\cal
S}_{\mathbb{F}}^2(t,T; {\mathbb R}^d),\ j\ge 0.$\ Moreover, $(\ref{aaa})$ can be equivalently rewritten in the following form:
\be\label{aab}\begin{array}{llll}
d\overline{X}_s^0 & = & E[b_0(\overline{X}_s^0,\overline{X}_s^{\ell_0})|
{\cal F}^{W^0}_T]ds+ E[\sigma_0(\overline{X}_s^0,\overline{X}_s^{\ell_0})|
{\cal F}^{W^0}_T]dW_s^0, & \overline{X}_t^0=x_0,\\
d\overline{X}_s^j & = & E[b_1(\overline{X}_s^0,\overline{X}_s^j,
\overline{X}_s^{\ell_j})|{\cal F}^{W^0,W^j}_T]ds+ E[\sigma_1(\overline{X}_s^0,
\overline{X}_s^j, \overline{X}_s^{\ell_j})|{\cal F}^{W^0,W^j}_T]dW_s^j, &
\overline{X}_t^j=\overline{x},\ j\ge 1,\\
\end{array}\ee
where $s$ runs the time interval $[t,T)$ and $\ell_j\ge 1$ is arbitrary but different from $j$, $j\ge 0.$
\ep

The equivalence between the first equation of $(\ref{aaa})$ and that
of $(\ref{aab})$ is evident, since we can replace $\overline{X}^1$
in (\ref{aaaz}) by $\overline{X}^{\ell_1}$ without changing $\mu_s$.
The equivalence between the second equation of $(\ref{aaa})$ and
that of $(\ref{aab})$ follows from the fact that
$\overline{X}^{\ell_j}$ defined by $(\ref{aab})$ is
$\mathbb{F}^{W^0,W^{\ell_j}}$-adapted and, knowing ${\cal
F}^{W^0}_T$ its law doesn't depend on $\ell_j\ge 1.$ Hence, since,
for $\ell_j\ge 1$ different from $j$, $W^0,W^j$ and $W^{\ell_j}$ are
independent,

\centerline{$P\{\overline{X}_s^{\ell_j}\in dy\|\ {\cal F}^{W^0,W^{j}}_T\}
=P\{\overline{X}_s^{\ell_j}\in dy\|\ {\cal F}^{W^0}_T\}=\mu_s(dy)$.}

\begin{proof} The proof is a direct consequence of the above observation.
Indeed, we have
\be\label{aabx}\begin{array}{llll}
d\overline{X}_s^0 & = & E[b_0(\overline{X}_s^0,\overline{X}_s^{1})|
{\cal F}^{W^0}_T]ds+ E[\sigma_0(\overline{X}_s^0,\overline{X}_s^{1})|
{\cal F}^{W^0}_T]dW_s^0, & \overline{X}_t^0=x_0,\\
d\overline{X}_s^1 & = & E[b_1(\overline{X}_s^0,\overline{X}_s^1,
\overline{X}_s^{2})|{\cal F}^{W^0,W^1}_T]ds+ E[\sigma_1
(\overline{X}_s^0,\overline{X}_s^1, \overline{X}_s^{2})|
{\cal F}^{W^0,W^1}_T]dW_s^1, & \overline{X}_t^1=\overline{x},\\
d\overline{X}_s^2 & = & E[b_1(\overline{X}_s^0,\overline{X}_s^2,
\overline{X}_s^{1})|{\cal F}^{W^0,W^2}_T]ds+ E[\sigma_1(\overline{X}_s^0,
\overline{X}_s^2, \overline{X}_s^{1})|{\cal F}^{W^0,W^2}_T]dW_s^2, &
\overline{X}_t^2=\overline{x};\\
\end{array}\ee
it's a finite-dimensional SDE with Lipschitz coefficients, and
standard estimates to show the existence and the uniqueness. Once
having the processes $\overline{X}^0$, $\overline{X}^1$\ and $\overline{X}^2$, we can
obtain the unique solution processes $\overline{X}^j, j\ge 3,$ in
(\ref{aab}) by choosing, for instance, $\ell_j=1.$
\end{proof}

The limit system $(\ref{aaa})$, or equivalently $(\ref{aab})$, is
related with $(\ref{3.1})$\ and $(\ref{3.2})$ through the following convergence
property.

\bp\label{p3} Under our standard assumptions we have that, for all $m\ge 1$, there is
some constant $C_m\in {\mathbb R}_+$ such that, for all $N\geq 0$,\ $0\leq
\ell\leq N,$ \be\label{3.25}\displaystyle
E[\sup_{r\in[t,T]}\big(|X_r^{\ell,N}-\overline{X}_r^{\ell}
|^{2}+\frac{1}{N}\sum_{\ell=1}^N|X_r^{\ell,N}-
\overline{X}_r^\ell|^2\big)^m]\le
\frac{C_m}{N^m}+C_m(\frac{1}{N}\sum_{\ell=1}^N
|x_\ell-\overline{x}|^2)^m.\ee \ep

For the reader's convenience we give the proof in Appendix 1. By adapting the argument of the proof of the above proposition we
also see the following lemma.

\bl\label{lemma3.47} Under our standard assumptions we have, for all bounded Lipschitz
functions $h: {\mathbb R}^d\times {\mathbb R}^d\rightarrow {\mathbb
R}$ and $g:{\mathbb R}^d\rightarrow {\mathbb
R}$, and for all $m\ge 1$, the existence of a real constant $C_m$ such
that, for all $N\ge 1,$ \be\label{3.47} \
\sup_{s\in[t,T]}E[\|\frac{1}{N}\sum_{\ell=1}^Nh(X_s^{0,N},
X_s^{\ell,N})- E[h(\overline{X}_s^{0},\overline{X}_s^{1})|{\cal F}_T^{W^{0}}]\|^{2m}]\le
C_m(\frac{1}{N}+\frac{1}{N} \sum_{\ell=1}^N
|x_\ell-\overline{x}|^2)^m, \ee
and
\be\label{3.47x} \
E[\|\frac{1}{N}\sum_{\ell=1}^N g(X_s^{\ell,N})-
E[g(\overline{X}_s^{1})|{\cal F}_T^{W^{0}}]\|^{2m}
|{\cal F}_T^{W^{0}}]\le C_m(\frac{1}{N}+\frac{1}{N}
\sum_{\ell=1}^N|x_\ell-\overline{x}|^2)^m, \quad s\in[t,T].
\ee
\el

After having reviewed the limit behavior of the forward equation
(\ref{3.1})-(\ref{3.2}) let us come now to the backward one.

\subsection{Limit behavior of the backward stochastic differential
equations}

\ \ \ \ \ \ \ \ In this subsection we discuss briefly the limit behavior of the
solution $(Y^{N},Z^{N}=(Z^{j,N})_{0\le j\le N})\in{\cal
S}^2_{\mathbb{F}}(t,T)\times L^2_{\mathbb{F}}(t,T;{\mathbb R}^d)^{N+1}$ of
BSDE
\be\label{3.50}\begin{array}{lll}& dY_s^{N}=-\frac{1}{N}
\sum_{\ell=1}^Nf(X_s^{0,N},X_s^{\ell,N},Y_s^{N}, Z_s^{0, N},Z_s^{\ell,
N})ds+\sum_{j=0}^NZ_s^{j, N}dW_s^{j},\ \ s\in[t,T],\\
& Y_T^{N}=\frac{1}{N} \sum_{\ell=1}^N\Phi(X_T^{0,
N},X_T^{\ell,N}),\\
\end{array}\ee
\noindent as $N$ tends to infinity, and we show $(Y^N, Z^{0, N})$\ converges to
the unique solution $(\overline{Y}, \overline{Z}^0)\in {\cal
S}^2_{\mathbb{F}^{W^0}}(t, T)\times L^2_{\mathbb{F}^{W^0}}(t, T;{\mathbb R}^d)$ of the BSDE:
\be\label{3.50a}\begin{array}{lll}& &  d\overline{Y}_s=
-f(\overline{X}^0_s, \mu_s, \overline{Y}_s, \overline{Z}_s^0, 0)ds
+\overline{Z}_s^0dW^0_s,\ \ s\in[t,T],\\
& & \overline{Y}_T=\Phi(\overline{X}_T^0,\mu_T);
\end{array}\ee
and $Z^{\ell, N}$\ converges to 0, for every $\ell\geq 1$\ (recall that $\mu_s(dy)=P\{\overline{X}_s^1\in dy\ | \ {\cal
F}_T^{W^0}\}$). We refer to the fact that it is by now standard that
under our assumptions on the coefficients the BSDEs
(\ref{3.50}) and (\ref{3.50a}) have a unique solution. In analogy to
the forward system we also see that the limit BSDE (\ref{3.50a}) can
be equivalently written in the form
\be\label{3.51}\overline{Y}_s=E[\Phi(\overline{X}_T^0,\overline{X}_T^\ell)
|{\cal{F}}_T^{W^0}]+\int_s^TE[f(\overline{X}_r^0,
\overline{X}_r^\ell,\overline{Y}_r,
\overline{Z}_r^{0},0)|{\cal{F}}_T^{W^0}]dr-\int_s^T\overline{Z}^0_r d W_r^0,\
s\in [t, T],\ee  for any $\ell \ge 1$. Moreover, we can have the
following statement on convergence:

\bp\label{BSDE_a} Under our standard assumptions on the coefficients we have
for all $m\ge 1$ the existence of a constant $C_m$ such that
\be\label{3.52}\begin{array}{lll}
 & &E\left[\sup_{s\in[t,T]}|Y_s^N-\overline{Y}_s|^{2m}+\left(\int_t^T|Z_s^{0,N}
  -\overline{Z}^0_s |^2ds+\sum_{\ell=1}^N\int_t^T|Z_s^{\ell,N}|^2dr\right)^m\right]\\
& &\ \le C_m\left(\frac{1}{N}+\frac{1}{N}
\sum_{\ell=1}^N| x_\ell -\overline{x}|^2\right)^m,\ \ \mbox{for all}\
N\ge 1.\\ \end{array}\ee \ep \noindent For the reader's convenience the proof is given in the Appendix 2.

\section{The stochastic differential game with N+1 participants}

\ \ \ \ \ \ Let $U={\mathbb R}^k$ be the control state space for the major player and
$V={\mathbb R}^m$ for the minor agents. To simplify the notation, we will
suppose from now on that the dimension $d$ used in the preceding
section is equal to $1$.

Our objective is to study the limit behavior of the stochastic
differential game between a major player and $N$ collectively
behaving minor agents, when $N\rightarrow +\infty.$

We denote by ${\mathbb{F}}^N$ the filtration generated by the Brownian
motion $W^{(N)}=(W^0,W^1,\dots, W^N)$ and augmented by all $P$-null
sets, and by $L^\infty_N(t,T;{\mathbb R}^\ell)$ $(\ell\ge 1)$ we denote the
space of all bounded, measurable functionals $\gamma: [t,T]\times
C([t,T])\rightarrow {\mathbb{R}}^\ell$ which are non-anticipating, i.e., for
all $s\in[t,T]$ and all $\psi,\ \psi'\in C([t,T])$ with
$\psi(r)=\psi'(r),\ r\in[t,s]$, it holds $\gamma_r(\psi)=\gamma_r(\psi'),\
r\in[t,s]$. Given initial positions $x^{(N)}=(x_0,x_1,\dots,x_N)$, a
feedback control $u\in L^\infty_N(t,T;U)$ of the major player, and
feedback controls $v^j\in L^\infty_N(t,T;V),\ 1\le j\le N,$ for the
$N$ minor agents, the dynamics of the major player $X^{0,N}$ and
those of the minor agents $X^{\ell,N}$ are defined by the system
\be\label{SDE_aa}\begin{array}{llll} dX_s^{0,N}&=&
\sigma_0(X^{0,N}_s)d\tilde{W}^0_s+\frac{1}{N}
\sum_{\ell=1}^Nb_0(X^{0,N}_s,X^{\ell,N}_s,Z^{0,N}_s)u_s(X^{(N)}_{\cdot})ds,\,
& X_t^{0,N}=x_0,\\
dX_s^{j,N}&=& \sigma_1(X^{j,N}_s)d\tilde{W}^j_s+
\frac{\varepsilon_N} {N}\sum_{\ell=1}^Nb_1(X^{0,N}_s,X^{\ell,N}_s,
Z^{j,N}_s)v_s^\ell(X^{(N)}_{\cdot}) ds,\, & X_t^{j,N}=x_j,
\end{array}\ee
$1\le j\le N,$\ driven by an $N+1$-dimensional Brownian motion $(\tilde{W}^0, \tilde{W}^1,\cdots, \tilde{W}^N)$. And the nonlinear payoff/cost functional is defined through the BSDE \be\label{BSDE_aa}\begin{array}{rll} dY^N_s&=&
- \frac{1}{N}\sum_{\ell=1}^Nf(X_s^{0,N},
X_s^{\ell,N}, Y^N_s, Z_s^{0,N}, Z_s^{\ell,N}, (u_s,v^\ell_s)(X^{(N)}_{\cdot}))ds +
\sum_{\ell=0}^NZ_s^{\ell,N}d\tilde{W}_s^{\ell},\ s\in[t,T],\\
Y^N_T&=& \frac{1}{N}\sum_{\ell=1}^N\Phi(X_T^{0,N},
X_T^{\ell,N}).
\end{array}\ee
Here $X^{(N)}=(X^{0,N},\dots,X^{N,N})$, and $\varepsilon_N>0$ is a
factor which converges to zero, as $N$ tends to $\infty$ (its role
will be discussed later). We want to study the above system in a
weak sense, i.e., we allow the driving $N$-dimensional Brownian
motion $(\tilde{W}^0,\dots, \tilde{W}^N)$ to depend on the control
processes. Assuming for simplicity that $\sigma_0\equiv 1,\
\sigma_1\equiv 1$, we use the Girsanov transformation

\smallskip

\centerline{$\begin{array}{lll}
d\tilde{W}^0_s& =& dW_s^0-I_{[t,T]}(s)\frac{1}{N}\sum_{\ell=1}^Nb_0
(X^{0,N}_s,X^{\ell,N}_s,Z^{0,N}_s)u_s(X^{(N)}_{\cdot})ds,\\
d\tilde{W}^j_s& =& dW_s^j-I_{[t,T]}(s)\frac{\varepsilon_N}{N}
\sum_{\ell=1}^Nb_1(X^{0,N}_s,X^{\ell,N}_s, Z^{j,N}_s)
v_s^\ell(X^{(N)}_{\cdot}) ds,
\end{array}$}

\smallskip

\noindent in order to reduce the study of the above system to the
resolution of the following simplified system:
\be\label{FBSDE_aa}\begin{array}{rll}
X_s^{j,N}&=& X_s^{j}=x_j+W_s^{j}-W_t^{j},\, \,  0\le j\le N,\\
dY^N_s&=&
- \{\frac{1}{N}\sum_{\ell=1}^Nf(X_s^{0,N},
X_s^{\ell,N},Y^N_s,Z_s^{0,N},Z_s^{\ell,N},u_s,v^\ell_s)+ (\frac{1}{N}
\sum_{\ell=1}^Nb_0(X^{0,N}_s,X^{\ell,N}_s,Z^{0,N}_s)u_s)Z^{0,N}_s\\
&
&+ \sum_{j=1}^N (\frac{\varepsilon_N}{N}\sum_{\ell=1}^Nb_1
(X^{0,N}_s,X^{\ell,N}_s,Z^{j,N}_s)v_s^\ell\big)Z^{j,N}_s\}ds
+\sum_{\ell=0}^NZ_s^{\ell,N}d{W}_s^{\ell},\ \ s\in[t,T],\\
Y^N_T& =& \frac{1}{N}\sum_{\ell=1}^N\Phi(X_T^{0,N},X_T^{\ell,N}).
\end{array}\ee
Here $u_s:=u_s(x^{(N)}+W_{\cdot}^{(N)}-W_t^{(N)}),\ v_s^\ell:=v_s^\ell(x^{(N)}+W_{\cdot}^{(N)}-W_t^{(N)}),\ s\in[t,T],\
1\le\ell\le N$, are now open-loop controls, and we can work for the
system (\ref{FBSDE_aa}) with open-loop controls $u\in{\cal
U}_N:=L_{\mathbb{F}^N}^\infty(t, T; U)$ and $v^\ell\in{\cal
V}_N:=L_{\mathbb{F}^N}^\infty(t, T; V)$, $1\le\ell\le N$, as long as our
saddle point controls are feedback controls.

Our objective is to study this latter stochastic differential game
and its limit behavior as $N\rightarrow +\infty.$ In order to
abbreviate the notation, given $\xi^{(N)}:=(x^{(N)}=(x_0,\dots,x_N)
,y,z^{(N)}=(z_0,\dots,z_N))\in {\mathbb R}^{N+1}\times {\mathbb R} \times {\mathbb R}^{N+1},\ u\in
U$ and $v^{(N)}=(v_1,...,v_N)\in V^N,$ we denote the driving
coefficient of the BSDE in (\ref{FBSDE_aa}) by

\smallskip

\centerline{$\begin{array}{ll} &  H_N(\xi^{(N)},u,v^{(N)})\\
& \displaystyle =\frac{1}{N}\sum^N_{\ell=1}f(x_0,x_\ell,y,z_0,z_\ell,u,v_\ell)
+\left(\frac{1}{N} \sum^N_{\ell=1}b_0(x_0,x_\ell,z_0)u\right)z_0+\varepsilon_N
\sum^N_{i=1}\left(\frac{1}{N}\sum^N_{\ell=1}b_1(x_0,x_\ell,z_i)v_\ell\right)z_i.
\end{array}$}

\smallskip

\noindent Let us make throughout our paper the following assumptions
on the coefficients involved in the definition of the Hamiltonian $H_N$:

\medskip

\noindent\textbf{Assumptions.} 1) $f=f(x_0,x_1,y,z_0,z_1,u,v):{\mathbb R}^5\times U \times V
\rightarrow {\mathbb R}$ is a continuous function with $f(\xi,.,.)\in C^1,\ \xi:=(\eta, z_1)\in {\mathbb R}^5,\ \eta:=(x_0,x_1,y,z_0)\in {\mathbb R}^4$, such that\\
$\bullet$ Ai) For some $\lambda >0:$
\be\begin{array}{lll} \((D_uf)(\xi,u,v)-(D_uf)(\xi,u',v), u-u'\)&
\leq &-\lambda |u-u'|^2,\\
\((D_vf)(\xi,u,v)-(D_vf)(\xi,u,v'), v-v'\)&\geq&\lambda |v-v'|^2,
\end{array}\ee
for all $\xi\in {\mathbb R}^5,\ (u,v),\ (u',v')\in U\times V;$\\
$\bullet$ Aii) For some constants $C>0$ and $\mu$ ($0<\mu<\lambda),$
\be\begin{array}{lll}
& &|f(\xi,0,0)|\leq C;\\
& &|f(\xi,u,v)-f(\xi',u,v)|\leq C|\xi-\xi'|;\\
& &|(D_uf)(\xi,u,v)-(D_uf)(\xi',u',v')|\leq C(|\xi-\xi'|+|u-u'|)+\mu|v-v'|,\ \
|(D_uf)(\xi,0,0)|\leq C;\\
& &|(D_vf)(\xi,u,v)-(D_vf)(\xi',u',v')|\leq C(|\eta-\eta'|+|v-v'|)+\mu|u-u'|,\ \
|(D_vf)(\xi,0,0)|\leq C;\\
\end{array}\ee
2) The coefficients $b_0:\mathbb{R}^3\rightarrow U$ and
$b_1:{\mathbb R}^3\rightarrow V$ are bounded and continuous, such that
$b_i(x_0,x_1,z)z$ is bounded and
Lipschitz in $(x_0,x_1,z)$, $i=0,1.$\\
3) The function $\Phi:\mathbb{R}\rightarrow \mathbb{R}$ is bounded
and Lipschitz.

\br 1. From {\rm Ai)} we see $f(\xi,\cdot,\cdot): U\times V \rightarrow
{\mathbb R}$ is a strictly concave-convex function. This strict
concavity-convexity will be crucial for the study of saddle point
controls for the associated stochastic differential game and the
study of their behavior, as $N$ tends to $+\infty.$

An easy example for such a function $f$ satisfying our assumptions
{\rm Ai)} and {\rm Aii)} is the following one for $k=m=1:$

\centerline{$f(\xi,u,v)=f_0(\xi,u,v)-\frac{1}{2}\alpha
|u|^2+\frac{1}{2}\gamma|v|^2,$}

\noindent where $\alpha,\ \gamma>0$ are strictly positive and $f_0\in
C^{1,2,2}_b({\mathbb R}^5\times U\times V)$ is such that
$D^2_{uu}f_0(\xi,u,v)\le (\alpha-\lambda)I_U(u)$ and
$D^2_{vv}f_0(\xi,u,v)\ge (\lambda-\gamma)I_V(v)$, $(\xi,u,v)\in
{\mathbb R}^5\times U\times V$, for some $\lambda>0,$ and
$|D^2_{uv}f_0(\xi,u,v)|\le \mu$,  $(\xi,u,v)\in {\mathbb R}^5\times U\times
V$, for some $\mu<\lambda$.

 2. Let us also give an example of a function $b_i(x_0,x_1,z)$
satisfying Assumptions 2): For this let $\tilde{b}_i: {\mathbb R}^3\rightarrow {\mathbb R}$ be a
bounded Lipschitz function and consider

\centerline{$\displaystyle b_i(x_0,x_1,z)=\tilde{b}_i(x_0,x_1,z)
\cdot\frac{1}{1+|z|}.$} \er

We also observe that from our assumptions on the function $f$ we get
\be\label{4.66}\begin{array}{rcl}
& &{\rm i)}\ |f(\xi,u,v)|\leq C(1+|u|^2+|v|^2);\\
& &{\rm ii)}\ |f(\xi,u,v)-f(\tilde{\xi},\tilde{u},\tilde{v})|\leq
C|\xi-\tilde{\xi}|+C(|u|+|\tilde{u}|+ |v|+|\tilde{v}|)(|u-\tilde{u}|+
|v-\tilde{v}|),
\end{array}\ee
for all $\xi=(x_{0},x_{1},y,z_{0},z_{1}),\ \tilde{\xi}=(\tilde{x}_{0},
\tilde{x}_{1},\tilde{y},\tilde{z}_{0},\tilde{z}_{1})\in\mathbb{R}^5$\
and $u,\ \tilde{u}\in U,\ v,\ \tilde{v}\in V.$\ Hence, taking into
account the above assumption on the coefficients $b_0$ and $b_1$, we
see that, for all
$\xi^{(N)}=(x^{(N)}=(x_0,\dots,x_N),y,z^{(N)}=(z_0,\dots,z_N)),\
\tilde{\xi}^{(N)}=(\tilde{x}^{(N)}=(\tilde{x}_0,
\dots,\tilde{x}_N),\tilde{y},\tilde{z}^{(N)}=(\tilde{z}_0,\dots,
\tilde{z}_N)),\ \xi^{(N,0)}=(x^{(N)}=(x_0,\dots,x_N),y,0)
\in\mathbb{R}^{N+1}\times\mathbb{R}\times\mathbb{R}^{N+1}$, and
$u\in U,\ v^{(N)}\in V^N$, we have
$$\begin{array}{rcl}
&\bullet &|H_N(\xi^{(N,0)},u,v^{(N)})|\leq C(1+|u|^2+\frac{1}
{N}|v^{(N)}|^2);\\
&\bullet &|H_N(\xi^{(N)},u,v^{(N)})-
H_N(\tilde{\xi}^{(N)}, u,v^{(N)})|\\
& &\leq C\left(1+|u|+\varepsilon_N\sum_{\ell=1}^N|v_\ell|\right)\left(
|x_0-\tilde{x}_0|+|y-\tilde{y}|+|z_0-\tilde{z}_0|+
\frac{1}{N}\sum_{\ell=1}^N(|x_\ell-\tilde{x}_\ell|+|z_\ell-\tilde{z}_
\ell|)\right)\\
& &\qquad +C\varepsilon_N\sum_{\ell=1}^N|v_\ell||x_\ell-\tilde{x}_\ell|.
\end{array}$$

The above properties of $H_N$ allow to apply the classical result on the
existence and uniqueness of the solution for BSDE:

\bp\label{BSDE_N-existence} For any admissible controls $u\in{\cal
U}_N=L^\infty_{\mathbb{F}^N}(t,T;U)$ and
$v^{(N)}=(v^{1,N},\dots,v^{N,N})\in{\cal
V}_N^N=L^\infty_{\mathbb{F}^N}(t,T;V)^N$, the BSDE in
(\ref{FBSDE_aa}) admits a unique solution
$(Y^N, Z^N=(Z^{0,N},\dots,Z^{N,N}))\in{\cal
S}^2_{\mathbb{F}^N}(t,T)\times L^2_{\mathbb{F}^N}(t,T)^{N+1}.$ In
order to indicate that this solution is associated with the controls
$(u,\ v^{(N)})$, we also write

\centerline{$\displaystyle (Y^{u,v^{(N)}},\ Z^{u,v^{(N)}}=
(Z^{0,u,v^{(N)}},\dots, Z^{N,u,v^{(N)}})):=(Y^N,\ Z^N=(Z^{0,N},
\dots, Z^{N,N}))$.} \ep

As already explained, we have two objectives: First we want to study
for each fixed $N\ge 1$ saddle point controls for our stochastic
differential game with $N$ collectively acting minor agents playing
against one major player, i.e, we are looking for a couple of
controls
$(\overline{u}^N,\ \overline{v}^{(N)}=(\overline{v}^{1,N},\dots,
\overline{v}^{N,N}))\in{\cal U}_N\times{\cal V}_N^N$ such that
\be\label{Saddle_point_N} Y_t^{u,\overline{v}^{(N)}}\le
Y_t^{\overline{u}^N,\overline{v}^{(N)}}\le
Y_t^{\overline{u}^N,v^{(N)}},\ \mbox{ for all }(u, v^{(N)})\in{\cal
U}_N\times{\cal V}_N^N.\ee In a second step we are interested in the
limit behavior of the saddle point controls
$(\overline{u}^N,\overline{v}^{(N)})$ when $N$ tends to $+\infty$,
and we want to characterize the limit controls as saddle point
controls for the limit stochastic differential game.

Let us begin with the study of the existence of saddle point
controls and its characterization for the game with $N+1$
participants. For this end, we have first to point out some useful
properties of the Hamiltonian $H_N$. So we observe that, as a direct
consequence of the above assumptions on our coefficients, we have
\be\label{4.2}\begin{array}{lll}
&\((D_uH_N)(\xi^{(N)},u,v^{(N)})-(D_uH_N)(\xi^{(N)},\tilde{u},v^{(N)}),
u-\tilde{u}\)\leq -\lambda |u-\tilde{u}|^2;\\
&\((D_{v^{(N)}}H_N)(\xi^{(N)},u,v^{(N)})-(D_{v^{(N)}}H_N)(\xi^{(N)},u,
\tilde{v}^{(N)}),v^{(N)}-\tilde{v}^{(N)}\)\\
&\ \ \  \geq \frac{\lambda}{N}\sum^N_{l=1}|v_l-\tilde{v}_l|^2 = \frac{\lambda}{N}|v^{(N)}-\tilde{v}^{(N)}|^2.\end{array}\ee

From ($\ref{4.2}$) it follows immediately that, for all
$\xi^{(N)}\in \mathbb{R}^{N+1}\times \mathbb{R}\times
\mathbb{R}^{N+1},$ the function $H_N(\xi^{(N)}, \cdot, \cdot):U\times
V^N\rightarrow \mathbb{R}$ has a unique saddle point. More
precisely, there exists a couple of Borel measurable feedback
controls $(\bar{u}^N, \bar{v}^{(N)}): {\mathbb R}^{N+1}\times {\mathbb
R}\times {\mathbb R}^{N+1}\rightarrow U\times V^N$, such that for
all $(u, v^{(N)})\in U\times V^N$,
\be\label{saddlepoint00}\displaystyle
H_N(\xi^{(N)}, u, \bar{v}^{(N)}(\xi^{(N)}))\leq
H_N(\xi^{(N)},\bar{u}^N(\xi^{(N)}),\bar{v}^{(N)}(\xi^{(N)}))\leq
H_N(\xi^{(N)},\bar{u}^N(\xi^{(N)}),v^{(N)}).\ee

\noindent However, the special form of our Hamiltonian $H_N$ allows
to compute the form of this saddle point in a more precise manner.
For this end, we let
$\tilde{v}^{(N)}=(\tilde{v}^N_1,\cdots,\tilde{v}^N_N): {\mathbb
R}^{N+1}\times {\mathbb R}\times {\mathbb R}^{N+1}\times
U\rightarrow V^N $  be a Borel measurable function such that
\be\label{4.4}\begin{array}{ll} &H_N(\xi^{(N)},u,\tilde{v}^{(N)}
(\xi^{(N)},u))=\displaystyle\inf_{v^{(N)}\in V^N}H_N(\xi^{(N)},u,v^{(N)}) \\
&=\displaystyle \frac{1}{N}\sum^N_{l=1}\inf_{v_l\in
V} \{f(x_0,x_l,y,z_0,z_l,u,v_l)+\varepsilon_N\sum^N_{i=1}
b_1(x_0,x_l,z_i)z_iv_l\}+\frac{1}{N}\(\sum^N_{l=1}
b_0(x_0,x_l,z_0)z_0\)u,\\
\end{array}\ee
and $\tilde{v}_N: {\mathbb R}\times{\mathbb R}\times{\mathbb R}
\times  {\mathbb R}^{N+1} \times U \rightarrow V $\ a measurable
function with \be\label{4.7a}\begin{array}{lll}& &
\displaystyle f(\xi_\ell,u,\tilde{v}_N(\xi_\ell^{(N)},u))+\varepsilon_N
\sum^N_{i=1}b_1(x_0,x_\ell,z_i)z_i\tilde{v}_N(\xi_\ell^{(N)},u)\\
& &=\displaystyle\inf_{v_l\in V}(f(\xi_\ell,u,v_\ell)+\varepsilon_N
\sum^N_{i=1}b_1(x_0,x_l,z_i)z_iv_\ell),
\end{array}\ee

\noindent for all $\xi_\ell^{(N)}=(x_0,x_\ell,y,z^{(N)})$ and all $u\in
U.$ Obviously,
\be\label{vtilde1}
\tilde{v}^N_\ell(\xi^{(N)},u)=\tilde{v}_N(\xi_\ell^{(N)},u),\ 1\le \ell\le N.
\ee
Let us also consider a measurable function $ \tilde{u}^N:
{\mathbb R}^{N+1}\times {\mathbb R}\times {\mathbb R}^{N+1}
\longrightarrow U$ such that
\be\label{4.6}\begin{array}{lll}& &H_N(\xi^{(N)},\tilde{u}^N(\xi^{(N)}),
\tilde{v}^{(N)}(\xi^{(N)}, \tilde{u}^N(\xi^{(N)})))\\
& &=\sup_{u\in
U}H_N(\xi^{(N)},u,\tilde{v}^{(N)}(\xi^{(N)},u)) \(=\sup_{u\in
U}\displaystyle\inf_{v^{(N)}\in V^N}
H_N(\xi^{(N)},u,v^{(N)})\),
\end{array}\ee
for all $\xi^{(N)}\in\mathbb{R}^{N+1}\times\mathbb{R}\times\mathbb{R}^{N+1}.$

\smallskip

For the above introduced functions we have the following result:

\bl\label{lemma4.7} Under our standard assumptions the unique saddle
point $(\bar{u}^N(\xi^{(N)}),\bar{v}^{(N)}(\xi^{(N)}))$ of
$H_N(\xi^{(N)},\cdot,\cdot)$ is of the form
\be(\bar{u}^N(\xi^{(N)}),\ \bar{v}^{(N)}(\xi^{(N)}))=(\tilde{u}^N(\xi^{(N)}),\ \tilde{v}^{(N)}(\xi^{(N)},\tilde{u}^N(\xi^{(N)}))),\
\xi^{(N)}\in\mathbb{R}^{N+1}\times\mathbb{R}\times\mathbb{R}^{N+1},\ee
i.e., it corresponds to a Stackelberg feedback strategy for a
2-player zero-sum game with the major player as leader and the
collectively acting minor agents as followers; it is the optimal
feedback strategy for the major player, if the collectively acting
minor agents react in an optimal way. \el

For the convenience of the reader we give the proof in Appendix 3.

\smallskip

With the help of our standard assumptions and Lemma \ref{lemma4.7}
we can now derive the following estimates for our saddle point
controls: \bl\label{lemma4.8} There exists some constant
$C\in\mathbb{R}$ such that, for all $N\ge
1,\ {\xi}^{(N)}=(x^{(N)}, y, z^{(N)}),\ \tilde{\xi}^{(N)}=(\tilde{x}^{(N)}, \tilde{y}, \tilde{z}^{(N)})\in\mathbb{R}^{N+1}
\times\mathbb{R}\times\mathbb{R}^{N+1},$
\be\label{4.23}\begin{array}{llll}
\hskip -0.7cm &{\rm{(i)}}\ |\overline{u}^N(\xi^{(N)})|&\leq& \min
\{C(1+\varepsilon_N\sum^N_{i=1}|z_i|), C(1+\varepsilon_N N)\};\\
\hskip -0.7cm &{\rm{(ii)}}\ |\bar{v}^N_{\ell}(\xi^{(N)})|(=|\tilde{v}_N
(\xi^{(N)}_{\ell},\tilde{u}^N(\xi^{(N)}))|)
&\leq& \min\{C(1+\varepsilon_N\sum^N_{i=1}|z_i|), C(1+\varepsilon_N N)\},
\ \ 1\leq {\ell} \leq N.\end{array}\ee
And
\be\label{4.33}\begin{array}{llll}
\hskip -0.7cm &{\rm{(iii)}}|\bar{u}^N(\xi^{(N)})-\bar{u}^N
(\tilde{\xi}^{(N)})|&\leq & C(|x_0-\tilde{x}_0|+|y-\tilde{y}|+|z_0-
\tilde{z}_0|)\\
\hskip -0.7cm& & &+C(1+\varepsilon_N N)(|x_0-\tilde{x}_0|+\frac{1}{N}
\sum^N_{{\ell}=1}|x_{\ell}-\tilde{x}_{\ell}|+\frac{1}{N}\sum^N_{{\ell}
=1}|z_{\ell}-\tilde{z}_{\ell}|);\\
\hskip -0.7cm &{\rm{(iv)}}|\bar{v}^N_{\ell}(\xi^{(N)})-\bar{v}^N_{\ell}
(\tilde{\xi}^{(N)})|&\leq& C(|x_0-\tilde{x}_0|+|y-\tilde{y}|+|z_0-
\tilde{z}_0|)\\
\hskip -0.7cm & & &+C(1+\varepsilon_N
N)(|x_0-\tilde{x}_0|+|x_{\ell}-
\tilde{x}_{\ell}|+\frac{1}{N}\sum^N_{i=1}(|x_i-\tilde{x}_i|+|z_i-
\tilde{z}_i|)).\end{array}\ee \el
For the proof of this lemma the reader is referred to Appendix 3.

With the help of the couple of feedback saddle point controls
$(\overline{u}^N,\ \overline{v}^{(N)}=(\overline{v}^N_1,\dots,\overline{v}^N_N))$
we now introduce the function
$$\overline{H}_{N}(\xi^{(N)}):=H_{N}(\xi^{(N)},{\bar{u}}^{N}
(\xi^{(N)}),\bar{v}^{(N)}(\xi^{(N)})),\ \xi^{(N)}\in
\mathbb{R}^{N+1}\times\mathbb{R}\times\mathbb{R}^{N+1},$$

\noindent which has the following properties:

\bl\label{lemma4.68} 1) Using the notation $\xi^{(N,0)}=(x^{(N)},y,
(z_{0},0,\cdots,0))\in\mathbb{R}^{N+1}\times\mathbb{R}\times\mathbb{R}^{N+1},$
we have, for some constant $C\in\mathbb{R}$ independent of $N\ge 1$,

\centerline{$\displaystyle |\overline{H}_{N}(\xi^{(N,0)})|\leq C.$}

\noindent 2) There is a constant $C_N$ (depending on $N\ge 1$) such
that, for all $\xi^{(N)},\ \tilde{\xi}^{(N)}\in \mathbb{R}^{N+1}\times\mathbb{R}\times\mathbb{R}^{N+1},$

\smallskip

\centerline{$\displaystyle
|\overline{H}_{N}(\xi^{(N)})-\overline{H}_{N}(\tilde{\xi}^{(N)})| \le
C_N|\xi^{(N)}-\tilde{\xi}^{(N)}|.$} \el The proof of this lemma
follows immediately from the assumptions on $b_0$ and $b_1$, and our
estimates (\ref{4.66}) for $f$,  combined with the estimates given
in Lemma \ref{lemma4.8}. Consequently, we can have the following
result:

\bl\label{lemma4.73} Under our assumptions, there exists a unique
solution $(\overline{Y}^{N}, \overline{Z}^{(N)}=
(\overline{Z}^{0,N},\dots,\overline{Z}^{N,N}))\in{\cal
S}^2_{\mathbb{F}^N}(t,T)\times L^2_{\mathbb{F}^N}(t,T)^{N+1}$ of the
BSDE
\begin{equation}\label{BSDE100}\begin{array}{rcl}
{\rm d}\overline{Y}_{s}^{N}&=&-\overline{H}_{N}(X_{s}^{(N)},
\overline{Y}_{s}^N,\overline{Z}_{s}^{(N)}){\rm
d}s+\sum_{j=0}^N\overline{Z}_{s}^{j,N}{\rm d}W_{s}^{j},\
\ s\in[t,T],\\
\overline{Y}_{T}^{N}&=&\frac{1}{N}\sum\limits_{{\ell}=1}^N
\Phi(X_{T}^{0},X_{T}^{{\ell,N}}).
\end{array}\end{equation}
\el
By putting
$$\overline{u}^N_s:=\overline{u}^N(X^{(N)}_s,\overline{Y}^N_s,
\overline{Z}^{(N)}_s),\ \ \overline{v}^{(N)}_s:=\overline{v}^{(N)}
(X^{(N)}_s,\overline{Y}^N_s,\overline{Z}^{(N)}_s),\, s\in[t,T],$$

\noindent and $\overline{u}^N_s:=u,\,
\overline{v}^{(N)}_s:=v^{(N)},\, s\in[0,t),$ for arbitrarily chosen
$u\in U$ and $v^{(N)}\in V^N$, we define admissible control
processes $\overline{u}^N\in{\cal U}_N,\
\overline{v}^{(N)}=(\overline{v}^{(1,N)},\dots,
\overline{v}^{(N,N)})\in{\cal V}_N^N$ (abusing notation we use
nearly the same notations as for the corresponding feedback
controls). The following statement is now a standard statement
following from the uniqueness of the solution of BSDE in
(\ref{FBSDE_aa}), from (\ref{saddlepoint00}) and the comparison
theorem for BSDEs.

\bp The above introduced couple of admissible controls
$(\overline{u}^N,\overline{v}^{(N)})\in{\cal U}_N\times{\cal V}_N^N$
forms a saddle point for the pay-off/cost functional ${\cal
U}_N\times{\cal V}_N^N\ni (u,v^{(N)})\rightarrow Y_t^{u,v^{(N)}}$.
More precisely, we even have

\smallskip

\centerline{$Y_s^{u,\overline{v}^{(N)}}\le Y_s^{\overline{u}^N,
\overline{v}^{(N)}}\le Y_s^{\overline{u}^N,v^{(N)}},\, s\in[t,T],$}

\smallskip

\noindent $P$-a.s., for all $(u,v^{(N)})\in{\cal U}_N\times
{\cal V}_N^N$. Moreover,

\smallskip

\centerline{$(Y^{\overline{u}^N,\overline{v}^{(N)}},Z^{\overline{u}^N,
\overline{v}^{(N)}})=(\overline{Y}^N,\overline{Z}^{(N)})$.}
\ep

\section{The limit game}
\subsection{The limit backward stochastic differential equation}

\ \ \ \ \ \ \ \ After having shown in the preceding subsection the existence of
saddle point controls $(\overline{u}^N,\overline{v}^{(N)})\in{\cal
U}_N\times{\cal V}_N^N$ for the game with $N+1$ agents, the
objective of this subsection is to introduce a BSDE, for which we
will prove later that it is the limit of the BSDEs for $N+1$
participants, if the saddle point controls are played, and we will
study the associated saddle point control processes for this limit
BSDE.

Let us begin with the introduction of the Hamiltonian for our limit
BSDE. For this end, we observe that Assumption Ai) allows to select
a measurable function $\overline{v}:\mathbb{R}^4\times U\rightarrow
V$ such that \be\label{a11} f(\xi,u,\overline{v}(\xi,u))=\inf_{v\in
V}f(\xi,u,v),\ (\xi=(x_0,x_1,y,z_0,0),u)\in {\mathbb R}^5\times U.\ee To
shorten the notation, we identify in what follows
$\xi=(x_0,x_1,y,z_0,0)$ with $(x_0,x_1,y,z_0)$ and we put
$f_{\overline{v}}(\xi,u):=f(\xi,u,\overline{v}(\xi,u)).$ From a
straight-forward computation using the assumptions Ai) and Aii) on
$f$, we obtain

\bl\label{lemma4.3} Under our standard assumptions on $f$, we have:\\
1. There is some $C\in\mathbb{R}$ such that, for all $(\xi,u),\ (\xi',
u')\in\mathbb{R}^4\times U,$
\be\label{4.39}\begin{array}{rcl}
& & |\bar{v}(\xi,u)|\leq C+\frac{\mu}{\lambda}|u|;\\
& & |\bar{v}(\xi,u)-\bar{v}(\xi',u')|\leq C|\xi-\xi'|+\frac{\mu}{\lambda}|u-u'|.
\end{array}\ee
2. For all $\xi\in\mathbb{R}^4,$ $f_{\overline{v}}(\xi,.)\in C^1$ and
\be\label{4.44}
(D_{u}f_{\bar{v}})(\xi,u)=(D_{u}f)(\xi,u,\bar{v}(\xi,u)),\
(\xi,u)\in {\mathbb{R}}^{4}\times U.\ee
\el
With the help of the function $f_{\bar{v}}$ we introduce the Hamiltonian
\be\label{4.49}\begin{array}{rcl}
\overline{H}(s,x_{0},y,z_{0},u):&=&E[f_{\bar{v}}(x_{0},
\overline{X}_{s}^{1},y,z_{0},0, u)+
b_{0}(x_{0},\overline{X}_{s}^{1},z_{0})z_{0}\cdot u|{\cal F}_T^{W^0}]\\
&=&E[f_{\bar{v}}(x_{0},\overline{X}_{s}^{1},y,z_{0},0,u)|{\cal
F}_T^{W^0}]+ E[b_{0}(x_{0},\overline{X}_{s}^{1},z_{0})z_{0}|{\cal
F}_T^{W^0}]\cdot u.
\end{array}\ee
We observe that $\overline{H}(s,x_{0},y,z_{0},u)$ is a continuous
random field which, for every fixed $(x_{0},y,z_{0},u)$, is
$\mathbb{F}^{W^0}$-progressively measurable. Moreover, we have the
following further properties of $\overline{H}$:

\bl\label{lemma4.4} 1) There exists some constant $C\in\mathbb{R}$
such that, $P$-a.s., for all $\xi,\ \xi'\in\mathbb{R}^3,\ u,\ u'\in U$,
\be\label{4.82}\begin{array}{rcl}
& &{\rm i)}\ |\overline{H}(s,\xi,u)|\leq C(1+|u|^2);\\
& &{\rm ii)}\ |\overline{H}(s,\xi,u)-\overline{H}(s,\xi',u')|
\leq C(1+|u|+|u'|)(|\xi-\xi'|+|u-u'|).\\
\end{array}\ee
\noindent 2) $P$-almost all trajectories of $\overline{H}(s,\xi,u)$ are
continuously differentiable in $u$, and\\
\be\label{4.57-1}(D_{u}\overline{H})(s,\xi,u)=E[(D_{u}f)(\xi,
\overline{X}_{s}^{1},u,\bar{v}(\xi,\overline{X}_{s}^{1},u))+
b_{0}(x_{0},\overline{X}_{s}^{1},z_{0})z_{0}|{\cal F}_T^{W^0}],\ee
\mbox{for all} $\xi:=(x_{0},y,z_{0})$,  where, with our convention
for $\xi$, $(\xi,\overline{X}_{s}^{1}):=(x_{0},\overline{X}_{s}^{1},
y,z_{0},0).$\\
3) We have, $P$-a.s., for all $(s, \xi)$ and all $u,\ u'\in
U$,\be\label{4.57}\begin{array}{rcl}
\((D_{u}\overline{H})(s,\xi,u)-(D_{u}\overline{H})(s,\xi,u'), u-u'\)\leq
-(\frac{\lambda^2-\mu^2}{\lambda})|u-u'|^2.
\end{array}\ee
\el
\noindent For the proof the reader is referred to Appendix 4.

A consequence of the preceding lemma is the strict concavity of the
function $\overline{H}(s,\xi,.):U\rightarrow {\mathbb R},$ uniform in
$(s,\ \xi=(x_0,y,z_0))\in[0,T]\times\mathbb{R}^3$ and uniform over
$\Omega.$ This implies that the random field
$\overline{u}:\Omega\times [0,T]\times {\mathbb R}^{3}\rightarrow
U$, $\mathbb{F}^{W^0}$-progressively measurable for frozen
$\xi=(x_0,y,z_0),$\ is uniquely defined by the relation
\be\label{4.58} \overline{H}(s,\xi,\bar{u}(s,\xi))=\sup_{u\in
U}\overline{H}(s,\xi,u),\ (s,\xi)\in [0,T]\times {\mathbb R}^{3}.\ee

\bl\label{lemma4.61} Under our standard assumptions\ $\bar{u}: \Omega\times[0,T]\times
{\mathbb R}^{3}\rightarrow U$ has the following properties:\\
\be\label{4.61}\begin{array}{lll}
|\bar{u}(s,\xi)|&\leq &C\\
|\bar{u}(s,\xi)-\bar{u}(s,\xi')|&\leq & C|\xi-\xi'|,
\end{array}\ee
for some constant $C\in\mathbb{R}$ not depending on $\ s\in [0,T],\
\xi,\ \xi'\in {\mathbb R}^{3}$. \el
\noindent For the proof the reader is referred to Appendix 4.

Let us now introduce the Hamiltonian of our limit BSDE: We
put \be\label{4.81}
\overline{H}(s,\xi)=\overline{H}(s,\xi,\bar{u}(s, \xi)),\
\xi=(x_{0},y,z_{0})\in\mathbb{R}^3. \ee \noindent Then we get
immediately from the Lemmas \ref{lemma4.4} and \ref{lemma4.61} that,
for some constant $C$ independent of $s\in[0,T],$ $\
\xi,\ \xi'\in\mathbb{R}^3,$ \be\label{4.83}\begin{array}{lll}
|\overline{H}(s,\xi)|&\leq & C;\\
|\overline{H}(s,\xi)-\overline{H}(s,\xi')|&\leq & C|\xi-\xi'|.
\end{array}\ee
Consequently, as $\overline{H}(.,\xi)$ is
$\mathbb{F}^{W^0}$-progressively measurable, we have the following

\bl\label{lemma4.85} There exists a unique solution
$(\overline{Y}, \overline{Z}^{0})\in{\cal S}^2_{\mathbb{F}^{W^0}}(t,T)\times
L^2_{\mathbb{F}^{W^0}}(t,T)$ of the BSDE \be\label{4.85}
d\overline{Y}_{s}=-\overline{H}(s,X_{s}^{0},\overline{Y}_{s},\overline{Z}_{s}^{0})
ds+\overline{Z}_{s}^{0}dW^0_{s},\
s\in [t,T),\ \
\overline{Y}_{T}=E[\Phi(X_{T}^{0},\overline{X}_{T}^1)|{\mathcal{F}}_{T}^{W^0}]
.\ee
\el

\subsection{The convergence to the limit game}

\ \  \ \ \ \ \ The main result in this subsection concerns the convergence of the
pay-off/cost functional of the game with $N+1$ participants under
saddle point controls to the solution of BSDE (\ref{4.85}). More
precisely, we have

\bt\label{Theorem 4.1} Let $(\overline{Y}, \overline{Z}^{0})\in {\cal
S}^2_{\mathbb{F}^{W^0}}(t,T)\times L^2_{\mathbb{F}^{W^0}}(t,T)$ be the
unique solution of the BSDE (\ref{4.85}) and
$(\overline{Y}^{N}, \overline{Z}^{(N)}=
(\overline{Z}^{0,N},\dots,\overline{Z}^{N,N}))\in{\cal
S}^2_{\mathbb{F}^N}(t,T)\times L^2_{\mathbb{F}^N}(t,T)^{N+1}$ that
of the BSDE (\ref{BSDE100}) introduced in Lemma \ref{lemma4.73}. Then, under the
assumption that $\varepsilon_N=O(N^{-3/4}),$ as $N\rightarrow
+\infty$, there exists a constant $C\in\mathbb{R},$ $N_0\ge 1$ and,
for all $m\ge 1,$  a constant $C_m\in\mathbb{R}$, such that for all
$N\ge N_0,$ \be\label{4.153}\begin{array}{lll}
&& {\rm i)}\displaystyle \ |\overline{Y}_{s}^{N}- \overline{Y_{s}}|
\leq C(\frac{1}{N}+\frac{1}{N}\sum_{l=1}^{N}|x_{l}-\overline{x}|^{2}
)^{\frac{1}{2}},\ s\in [t,T];\\
&& {\rm ii)}\displaystyle \ E[(\int_{t}^{T}|\overline{Z}_{s}^{0,N}-
\overline{Z}_{s}^{0}|^{2}ds+\sum_{l=1}^{N}\int_{t}^{T} |
\overline{Z}_{s}^{l,N}|^{2}ds)^{m}|\mathcal {F}_{t}^{W^{0}}]\leq
C_{m}(\frac{1}{N}+\frac{1}{N}\sum_{l=1}^{N}|x_{l}-\bar{x}|^{2})^{m}.
\end{array}\ee

\et Taking into account that we have supposed that
$\frac{1}{N}\sum_{l=1}^{N}|x_{l}-\bar{x}|^{2}\rightarrow 0,$ as
$N\rightarrow +\infty$, an obvious consequence of the theorem is the
following

\bc\label{corollary4.1} Under our standard assumptions on the
coefficients as well as the condition that
$\varepsilon_N=O(N^{-3/4}),$ as $N\rightarrow +\infty$, we have for
all $m\ge 1$, in $L^{m},$
$$|\overline{Y}_{s}^{N}- \overline{Y_{s}}|^2+\int_{t}^{T}|
\overline{Z}_{s}^{0,N} -\overline{Z}_{s}^{0}|^{2}ds+\sum_{l=1}^{N}
\int_{t}^{T} |\overline{Z}_{s}^{l,N}|^{2}ds\ \longrightarrow 0.$$
\ec

The proof of the theorem is prepared by two auxiliary lemmas. The
first lemma analyzes the limit behavior of the couple of saddle
point controls
\begin{equation}\label{saddlepoint}\overline{u}_s^N=\overline{u}^N(X_s^{(N)},
\overline{Y}_s^N,\overline{Z}_s^{(N)}),\ \ \overline{v}^{(N)}_s=
(\overline{v}^{(\ell,N)}_s=\overline{v}^N_\ell(X^{(N)}_s,
\overline{Y}^N_s,\overline{Z}^{(N)}_s))_{1\le\ell\le N},\end{equation}

\noindent as the number $N$ of collectively playing minor agents
tends to infinity. Putting
\begin{equation}\label{saddlepoint1}\overline{u}_s:=\overline{u}(s,X_s^0,\overline{Y}_s,
\overline{Z}^0_s),\ \ \overline{v}^j_s=\overline{v}(X_s^0,\overline{X}_s^j,
\overline{Y}_s,\overline{Z}_s^0, \overline{u}_s),\ s\in[t,T],\end{equation}

\noindent (recall the definition of $\overline{u}(.)$ by
(\ref{4.58}) and that of $\overline{v}(.)$ by (\ref{a11})), we can
establish the following result:

\bl\label{auxlemma11} Under our standard assumption of Lipschitz
continuity  on the coefficients we have
\be\label{4.119}\begin{array}{rcl}
&&|\overline{u}_s^N-\overline{u}_s|+|\overline{v}^{(j,N)}_s-
\overline{v}^j_s|\\
&\leq&C(|\overline{Y}_s^N-\overline{Y}_s|+|\overline{Z}_s^{0,N}
-\overline{Z}_s^0| +(1+\varepsilon_N\cdot
N)\frac{1}{N}\sum_{\ell=1}^N|\overline{Z}_s^{\ell,N}|)+R^N_s,
\end{array}\ee
$s\in[t,T],\ 1\le j\le N,\  N\geq 1,$ where, for all $m\ge 1,$
the remainder $R^N_s$ satisfies the following estimate for some
constant $C_m$ only depending on $m$:
$$E[|R^N_s|^{2m}|{\cal F}_T^{W^0}]
\le C_m(\frac{1}{N}+\frac{1}{N}\sum_{\ell=1}^N|x_\ell-
\overline{x}|^2)^m,\, s\in[t,T],\ N\ge 1.$$
\el

\noindent The proof of this lemma is presented in Appendix 5.

The above lemma allows to characterize the limit behavior of
the Hamiltonian
\begin{equation}\label{H}\overline{H}_N(\xi^{(N)}):
=H_N(\xi^{(N)},\overline{u}^N(\xi^{(N)}),\overline{v}^N(\xi^{(N)}))\end{equation}

\noindent along the saddle points $(\overline{u}^N(\xi^{(N)}),\
\overline{v}^N(\xi^{(N)}))$, as $N$ tends $+\infty.$

\bl\label{lemma4.138} Under the standard assumptions on the coefficients
the following estimate with some suitable constant $C$ holds true:
\be\label{4.138}\begin{array}{lll}
&&|\overline{H}_N(X_s^{(N)},\overline{Y}^N_s,\overline{Z}^{(N)}_s)-
\overline{H}(s,X_s^0,\overline{Y}_s,\overline{Z}_s^0)|\\
&\leq&C(|\overline{Y}_s^N-\overline{Y}_s|+|\overline{Z}_s^{0,N}-
\overline{Z}_s^0|)+C(\frac{1}{N}+\varepsilon_N+\varepsilon_N^2N)
\sum_{\ell=1}^N|\overline{Z}_s^{\ell,N}|+R^N_s.
\end{array}\ee
As in Lemma \ref{auxlemma11}, the remainder $R^N_s$ satisfies the following
estimate:$$E[|R^N_s|^{2m}|{\cal F}_T^{W^0}]
\le C_m(\frac{1}{N}+\frac{1}{N}\sum_{\ell=1}^N|x_\ell-\overline{x}|^2)^m,
\, s\in[t,T],\ N\ge 1,\ m\ge 1.$$
\el
The proof of this lemma is given in Appendix 5.

With the both preceding lemmas we are now able to prove Theorem \ref{Theorem 4.1}.

\begin{proof} (of Theorem \ref{Theorem 4.1}). Recall that $(\overline{Y}^N,\
\overline{Z}^{(N)})$ is introduced in Lemma \ref{lemma4.73} as the unique
solution of the BSDE (\ref{BSDE100}) with Hamiltonian $\overline{H}_N$, while the couple
of processes $(\overline{Y},\overline{Z}^0)$ is the unique solution of BSDE
(\ref{4.85}) with Hamiltonian $\overline{H}$. By applying It\^{o} formula to
$|\overline{Y}_s^N-\overline{Y}_s|^2$, we get from
BSDE standard estimates

\be\label{4.141a}\begin{array}{lll}
&&|\overline{Y}_s^N-\overline{Y}_s|^2+E[\int_s^T|\overline{Z}_r^{0,N}-
\overline{Z}_r^0|^2r+\sum_{\ell=1}^N\int_s^T|\overline{Z}_r^{\ell,N}|^2dr|
\mathcal{F}_s^{W^0}]\\
&\leq&CE[|\frac{1}{N}\sum_{\ell=1}^N\Phi(X_T^0,X_T^\ell)-E[\Phi(X_T^0,
\overline{X}_T^\ell)|\mathcal{F}_T^{W^0}]|^2|\mathcal{F}_s^{W^0}]\\
& &+CE[\int_s^T|\overline{Y}_r^N-\overline{Y}_r||\overline{H}_N(X_r^{(N)},
\overline{Y}^N_r,\overline{Z}_r^{(N)})-\overline{H}(r,X_r^0,\overline{Y}_r,
\overline{Z}_r^0)|dr|\mathcal{F}_s^{W^0}].
\end{array}\ee
From Lemma \ref{lemma3.47} we know that

\smallskip

\centerline{$E[|\frac{1}{N}\sum_{\ell=1}^N\Phi(X_T^0,X_T^\ell)-E[\Phi(X_T^0,
\overline{X}_T^\ell)|\mathcal{F}_T^{W^0}]|^2| \mathcal{F}_s^{W^0}]\le
C(\frac{1}{N}+\frac{1}{N}\sum_{\ell=1}^N|x_\ell-\overline{x}|^2).$}

\smallskip

\noindent Consequently, Lemma \ref{lemma4.138} yields that, for $\rho>0$ small
(we will specify $\rho$ later) and a constant $C_\rho$ depending on $\rho$,
\be\label{4.141}\begin{array}{lll}
&&|\overline{Y}_s^N-\overline{Y}_s|^2+E[\int_s^T|\overline{Z}_r^{0,N}-
\overline{Z}_r^0|^2dr+\sum_{\ell=1}^N\int_s^T|\overline{Z}_r^{\ell,N}|^2dr|
\mathcal{F}_s^{W^0}]\\
&\leq&C_\rho(\frac{1}{N}+\frac{1}{N}\sum_{\ell=1}^N|x_\ell-\overline{x}|^2)
+C_\rho\int_s^TE[|\overline{Y}_r^N-\overline{Y}_r|^2|\mathcal{F}_s^{W^0}]dr+
\frac{1}{2}E[\int_s^T|\overline{Z}_r^{0,N}-\overline{Z}_r^0|^2dr|
\mathcal{F}_s^{W^0}]\\
& &+C_\rho\int_s^TE[|R^N_r|^2|\mathcal{F}_s^{W^0}]dr+\rho(\frac{1}{N}+
\varepsilon_N+\varepsilon_N^2N)^2E[\int_s^T(\sum_{\ell=1}^N|
\overline{Z}_r^{\ell,N}|)^2dr|\mathcal{F}_s^{W^0}],
\end{array}\ee
where the estimate of the remainder $R^N_r,\ r\in[t,T],$ is given by Lemma
\ref{lemma4.138}. Obviously,
$$\big(\frac{1}{N}+\varepsilon_N+\varepsilon_N^2N\big)^2
\big(\sum_{\ell=1}^N|\overline{Z}_r^{\ell,N}|\big)^2\leq C
\big(\frac{1}{N}+N\varepsilon_N^2+N^3\varepsilon_N^4\big)\sum_{\ell=1}^N|
\overline{Z}_r^{\ell,N}|^2,$$

\noindent and it is here, where we have to use our assumption that
$\varepsilon_N=O(N^{-3/4}),$ as $N\rightarrow \infty$. Indeed, this
assumption allows to get
$$\rho\big(\frac{1}{N}+\varepsilon_N+\varepsilon_N^2N\big)^2\big(
\sum_{\ell=1}^N|\overline{Z}_s^{\ell,N}|\big)^2\leq\frac12\sum_{\ell=1}^N|
\overline{Z}_s^{\ell,N}|^2,\ \ N\geq1,\ s\in[t,T],$$

\noindent for a sufficiently small chosen $\rho>0.$ Hence, with such a
choice of $\rho$ we obtain
$$\begin{array}{lll}
&&|\overline{Y}_s^N-\overline{Y}_s|^2+E[\int_s^T|\overline{Z}_r^{0,N}-
\overline{Z}_r^0|^2dr+\sum_{\ell=1}^N\int_s^T|\overline{Z}_r^{\ell,N}|^2dr|
\mathcal{F}_s^{W^0}]\\
&\leq&C(\frac{1}{N}+\frac{1}{N}\sum_{\ell=1}^N|x_\ell-\overline{x}|^2)+
CE[\int_s^T|\overline{Y}_r^N-\overline{Y}_r|^2dr|\mathcal{F}_s^{W^0}],
\quad s\in[t,T],\ N\geq 1.
\end{array}$$

\noindent Consequently, from Gronwall's inequality we have for all
$N\ge 1,$ $P$-a.s.,
\be\label{4.147}|\overline{Y}_s^N-\overline{Y}_s|^2+E\big[\int_s^T|
\overline{Z}_r^{0,N}-\overline{Z}_r^0|^2dr+\sum_{\ell=1}^N\int_s^T|
\overline{Z}_r^{\ell,N}|^2dr|\mathcal{F}_s^{W^0}\big]\leq C\big(\frac{1}{N}
+\frac{1}{N}\sum_{\ell=1}^N|x_\ell-\overline{x}|^2\big),\, s\in[t,T],\ee

\noindent i.e., estimate (\ref{4.153})-i) is proved. It still remains to show ii). For this end, we consider the difference
between the BSDEs solved by $(\overline{Y}^N,\overline{Z}^{(N)})$ and
by $(\overline{Y},\overline{Z}^{0})$, and we apply  to this difference
the Burkholder-Davis-Gundy inequality. Thus, using Lemma \ref{lemma3.47}
and Lemma \ref{lemma4.138}, we get, for every $m\ge 1$ and some
constant $C_m$ depending on $m$,
\be\label{4.149}\begin{array}{lll}
&&E[(\int_s^{s+\delta}|\overline{Z}_r^{0,N}-\overline{Z}_r^0|^2dr+
\sum_{\ell=1}^N\int_s^{s+\delta}|\overline{Z}_r^{\ell,N}|^2dr)^m|
\mathcal{F}_t^{W^0}]\\
&\leq&C_m E[\sup_{r\in[s,s+\delta]}\big|\int_s^r(\overline{Z}^{0,N}_v-
\overline{Z}^0_v)dW^0_v+\sum_{\ell=1}^N\int_s^{r}
\overline{Z}_v^{\ell,N}dW^\ell_v|^{2m}|{\cal F}_t^{W^0}]\\
&\leq& C_m E\big[\sup_{r\in[s,{s+\delta}]}|\overline{Y}_r^N-
\overline{Y}_r|^{2m}+\big(\int_s^{s+\delta}|\overline{H}_N(X_r^{(N)},
\overline{Y}_r^N,\overline{Z}_r^{(N)})-\overline{H}(r,X_r^0,\overline{Y}_r,
\overline{Z}^0_r)|dr\big)^{2m}|{\cal F}_t^{W^0}\big]\\
&\leq&C_m(\frac{1}{N}+\frac{1}{N}\sum_{\ell=1}^N|x_\ell-\overline{x}|^2)^m
+C_mE\big[\big(\int_s^{s+\delta}(|\overline{Z}^{0,N}_r\hskip -1mm
-\overline{Z}^0_r|\hskip -0.5mm+\hskip -0.5mm(\frac{1}{N}\hskip -0.5mm+
\hskip -0.5mm\varepsilon_N
\hskip -0.5mm+\hskip -0.5mm\varepsilon_N^2N)\sum_{\ell=1}^N\hskip -0.5mm|
\overline{Z}_r^{\ell,N}|)dr\big)^{2m}|{\cal F}_t^{W^0}\big]\\
&\leq&C_m(\frac{1}{N}+\frac{1}{N}\sum_{\ell=1}^N|x_\ell-\overline{x}|^2)^m
+C_m \delta^m E[(\int_s^{s+\delta}|\overline{Z}_r^{0,N}\hskip -0.5mm-\hskip
-0.5mm\overline{Z}_r^0|^2dr+\sum_{\ell=1}^N\int_s^{s+\delta}|
\overline{Z}_r^{\ell,N}|^2dr)^{m}|\mathcal{F}_t^{W^0}],\\
\end{array}\ee
for all $t\leq s<s+\delta\leq T$ and all $N\geq N_0$, for some $N_0\geq 1$
large enough (recall that $\varepsilon_N=O(N^{-3/4})$). Now inequality (\ref{4.153})-ii) follows easily.
\end{proof}

After having proved Theorem \ref{Theorem 4.1}, we can combine it with Lemma
\ref{auxlemma11}, in order to improve its statement concerning the convergence
of the saddle point controls of the game with $N$ collectively acting minor
agents, when $N$ tends to $+\infty$. Then we obtain easily the following
result:

\bt\label{theorem111} Under our standard assumptions as well as the condition
that $\varepsilon_N=O(N^{-3/4}),$ as $N\rightarrow +\infty$, we have that for
all $m\ge 1$  there is a constant $C_m\in\mathbb{R}$, such that for all
$N\ge N_0,$
\be\label{4.161}\begin{array}{lll}
& {\rm i)}&\ E[(\int_{t}^{T}(|\overline{u}_{s}^{N}-\overline{u}_{s}|+|\overline{
v}_s^{(j,N)}-\overline{v}_s^j|)^{2}ds)^{m}|\mathcal {F}_{t}^{W^{0}}]\leq C_{m}
(\frac{1}{N} +\frac{1}{N}\sum_{\ell=1}^{N}|x_{\ell}-\overline{x}|^{2})^{m},\,
1\le j\le N;\\
& {\rm ii)}&\ E[(\int_{t}^{T}| \frac{1}{N}\sum_{\ell=1}^N\psi(\overline{v}^{
\ell,N}_s)-E[\psi(\overline{v}_{s})|{\cal F}_T^{W^0}]|^{2}ds)^{m}|
\mathcal{F}_{t}^{W^{0}}]\\
& &\ \ \leq C_{m}(\frac{1}{N}+\frac{1}{N}\sum_{\ell=1}^{N}|x_{\ell}-
\overline{x}|^{2})^{m},\ 1\le j\le N,\ \mbox{for all}\ \mbox{bounded
Lipschitz functions} \ \psi,
\end{array}\ee
where $\overline{v}_s:=\overline{v}(X_{s}^{0},\overline{X}_{s}^{1},
\overline{Y}_{s},\overline{Z}_{s}^{0},\overline{u}_{s}),\ s\in [t, T]$
(recall that $\overline{v}^j_s=\overline{v}(X_s^0, \overline{X}_s^j, \overline{Y}_s,
\overline{Z}_s^0, \overline{u}_s)$).
As, due to our assumption, $\frac{1}{N}\sum_{\ell=1}^{N}|x_{\ell}-
\overline{x}|^{2}\rightarrow 0$, this means in particular, that the
left-hand sides of the above estimates converge to zero, as $N$ tends
to $+\infty.$
\et
\br Theorem \ref{theorem111} describes the limit behavior of the saddle
point controls of the game with $N$ collectively acting minor agents.
While statement {\rm i)} says that the saddle point control of the major player
$\overline{u}^N$ and, for all $j\ge 1$, that of the $j$-th minor agent
$\overline{v}^{j,N}$ converge, for all $m\ge 1$, in $L^m(\Omega,L^2([t,T]))$
to the processes $\overline{u}$ and $\overline{v}^j$, respectively, statement
{\rm ii)} says that, if we identify in the limit the collectively acting minor
agents with a limit player, whose dynamics is $\overline{X}^1_s=\overline{x}+
W^1_s-W^1_t,\, s\in[t,T]$, we get the associated control process
$\overline{v}_s=\overline{v}(X_{s}^{0},\overline{X}_{s}^{1},
\overline{Y}_{s},\overline{Z}_{s}^{0},\overline{u}_{s})$ as a weak limit. More
precisely, with the convention
$E[\delta_{\overline{v}_s}|{\cal F}_T^{W^0}](\psi)=E[\delta_{\overline{v}_s}
(\psi)|{\cal F}_T^{W^0}]=E[\psi({\overline{v}_s})|{\cal F}_T^{W^0}]$, for
$\psi\in C_b(\mathbb{R})$,
$$\frac{1}{N}\sum_{\ell=1}^N\delta_{\overline{v}_s^{\ell, N}}\rightarrow
E[\delta_{\overline{v}_s}|{\cal F}_T^{W^0}]\ \mbox{weakly, in}\ L^m(\Omega,L^2
([t, T])).$$
\er
\begin{proof} Statement i) is a direct consequence from Lemma \ref{auxlemma11}
combined with Theorem \ref{Theorem 4.1}, while statement ii) of the theorem
follows easily from i) and Lemma \ref{lemma3.47}.
\end{proof}

Recall the definition of $\overline{u}$ and also that
\be\label{a1a1}\overline{v}_s=\overline{v}(X_{s}^{0},\overline{X}_{s}^{1},
\overline{Y}_{s},\overline{Z}_{s}^{0},\overline{u}_{s}),\, s\in [t,T].\ee
\noindent From Lemma \ref{lemma4.3} and the boundedness of the process
$\overline{u}$ (Lemma \ref{lemma4.61}) it follows that of the process $\overline{v}$. Consequently,
from (\ref{4.161}) we have, for all $m\ge 1$ and all $N\ge N_0,$
\be
E[(\int_{t}^{T}| \frac{1}{N}\sum_{\ell=1}^N\overline{v}^{\ell,N}_s-
E[\overline{v}_{s}|{\cal F}_T^{W^0}]|^{2}ds)^{m}|\mathcal {F}_{t}^{W^{0}}]\leq C_{m}(\frac{1}{N}+\frac{1}{N}\sum_{\ell=1}^{N}|x_{\ell}-\overline{x}|^{2})^{m}.
\ee
Our objective is to characterize the couple $(\overline{u},\ \overline{v})$ obtained in
Theorem \ref{theorem111} as saddle point control for a limit 2-person zero-sum
stochastic differential game.

\bigskip

In order to define this 2-person zero-sum game, we introduce the function
\be\label{4.169} F(x_{0}, x_{1}, y, z_{0}, u, v):=f(x_{0}, x_{1}, y, z_{0},
0, u, v)+b_{0}(x_{0}, x_{1}, z_{0})z_{0}u,\ee
$(x_{0}, x_{1}, y, z_{0}, u, v)\in\mathbb{R}^4\times U\times V$. We consider as
space of admissible controls for Player 1 the set
$\overline{\cal U}=L^\infty_{\mathbb{F}^{W^0}}(t,T;U)$ and for Player 2 the set
$\overline{\cal V}=L^\infty_{\mathbb{F}^{1}}(t,T;V)(=L^\infty_{\mathbb{F}^{W^0,
W^1}}(t,T;V)$ (recall that $\mathbb{F}^{N}$ is the filtration generated by the
Brownian motions $W^0, W^1,\dots, W^N$ and augmented by all $P$-null sets).
Given a couple of admissible controls $(u,v)\in\overline{\cal U}\times
\overline{\cal V}$, we consider the BSDE
\be\label{4.170}\begin{array}{lll}
d\overline{Y}^{u,v}_s&=&-E[F(X_s^0,\overline{X}_s^1,\overline{Y}_s^{u,v},
\overline{Z}_s^{u,v},u_s,v_s)|{\cal F}_T^{W^0}]ds+ \overline{Z}^{u,v}_sdW_s^0,
\ s\in[t,T],\\
\overline{Y}^{u,v}_T&=&E[\Phi(X_T^0,\overline{X}_T^1)|{\cal F}_T^{W^0}],
\end{array}\ee
governed by

\centerline{$\begin{array}{llll}
X_s^0& = &x_0+W_s^0-W_t^0 &\mbox{--the dynamics of Player 1 (the major
player)}\\
\overline{X}_s^1&=&\overline{x}+W^1_s-W^1_t &\mbox{--the dynamics
of Player 2 (the collectively acting minor agents }
\end{array}$}
\ \ \ \hskip4.3cm who, in the averaging limit, amalgamate to a single player).

\smallskip

\noindent From our standard assumptions on the functions $f$ and $b_0$
we have

\smallskip

\centerline{$\begin{array}{rcl}
|F(x_{0}, x_{1}, y, z_{0}, u, v)|&\leq &C(1+|u|^{2}+|v|^{2}),\\
|F(x_{0}, x_{1}, y, z_{0}, u, v)-F(x_{0}, x_{1}, y', z_{0}', u, v)|
&\leq &C(1+|u|)(|y-y'|+|z_{0}-z_{0}'|),
\end{array}$}

\smallskip

\noindent for all $(x_0, x_1, u, v),\ y,\ y'$ and $z_{0},\ z_{0}'.$\ Consequently,
from standard BSDE arguments we have

\bl\label{lemma5.1} Under our standard assumptions, for any $(u,v)\in
\overline{\cal U}\times\overline{\cal V}$, BSDE (\ref{4.170}) has a unique
solution $(\overline{Y}^{u,v},\overline{Z}^{u,v})
\in{\cal S}^2_{\mathbb{F}^{W^0}}(t,T)\times L^2_{\mathbb{F}^{W^0}}(t,T).$
Moreover, for the controls $(\overline{u},\overline{v})\in
\overline{\cal U}\times\overline{\cal V}$ defined respectively in (\ref{saddlepoint1}) and
(\ref{a1a1}), the couple $(\overline{Y}^{\overline{u},\overline{v}},
\overline{Z}^{\overline{u},\overline{v}})$ coincides with the unique
solution $(\overline{Y},\overline{Z}^{0})$ of BSDE (\ref{4.85}).
\el

\begin{proof} The fact that BSDE (\ref{4.170}) has a unique solution
$(\overline{Y}^{u,v},\overline{Z}^{u,v})\in{\cal S}^2_{\mathbb{F}^{W^0}}
(t,T)\times L^2_{\mathbb{F}^{W^0}}(t,T)$ is an easy consequence of the
above properties of the function $F$.

Indeed, thanks to them, given $(u, v)\in\overline{\cal U}\times
\overline{\cal V}$, the coefficient $E[F(X_s^0,\overline{X}_s^1,
y,z_0,u_s,v_s)|{\cal F}_T^{W^0}]$ is, for all $(y, z_0)\in \mathbb{R}^2$,
$\mathbb{F}^{W^0}$-progressively measurable, Lipschitz in $(y, z_0)$ and
bounded, uniformly in $(s, \omega)\in[0,T]\times\Omega.$ The second part
of the statement, the assertion that $(\overline{Y}^{\overline{u},\overline{v}}, \overline{Z}^{\overline{u},\overline{v}})=(\overline{Y},
\overline{Z}^{0})$, follows the fact that the driving coefficients of the
both corresponding BSDEs coincide. Indeed, taking into account the
definition of the controls $\overline{u}$ and $\overline{v}$, we have
$$\begin{array}{lll}
& &E[F(X_{s}^{0},{\overline{X}}_{s}^{1},\overline{Y}_{s},\overline{Z}_{s}^{0},
\bar{u}_{s},\bar{v}_{s})|\mathcal {F}_{T}^{W_{0}}]=E[f_{\bar{v}}(X_{s}^{0}, {\overline{X}}_{s}^{1},\overline{Y}_{s},
\overline{Z}_{s}^{0},\bar{u}_{s})+b_{0}(X_{s}^{0}, {\overline{X}}_{s}^{1},
\overline{Z}_{s}^{0})\overline{Z}_{s}^{0}.\bar{u}_{s}|
\mathcal {F}_{T}^{W_{0}}]\\
&=&\overline{H}(s,X_{s}^{0},\overline{Y}_{s},\overline{Z}_{s}^{0},
\bar{u}_{s})=\overline{H}(s, {X}_{s}^{0},\overline{Y}_{s},
\overline{Z}_{s}^{0})\ \ (\mbox{from}\ (\ref{4.81})),\ \ s\in [t,T]
\end{array}$$
(recall that due to our convention $f_{\overline{v}}(x_0,x_1,y,z_0,u)=
f_{\overline{v}}(x_0,x_1,y,z_0,0,u)$).

\smallskip

\noindent The stated result follows now from the uniqueness of the solution
for BSDE (\ref{4.85}).
\end{proof}

With the help of BSDE (\ref{4.170}) let us now introduce the pay-off/cost
functional for our 2-person zero-sum limit game:
\begin{equation}J(u,v):=\overline{Y}^{u,v}_t,\ \ (u,v)\in\overline{\cal U}
\times\overline{\cal V}.\end{equation}

\smallskip

\noindent Player 1 (the major player) wants to maximize his payoff $J(u,v)$
through the controls $u\in\overline{\cal U}$, while Player 2---the amalgamated collectively acting minor agents---wants to minimize the
loss $J(u,v)$ by using the controls $v\in\overline{\cal V}$.

For the such defined game we have the following characterization:

\bt\label{theorem150bb} The limit $(\overline{u},\overline{v})\in
\overline{\cal U}\times\overline{\cal V}$\ in the sense of Theorem
\ref{theorem111} of the couples of saddle point controls $(\overline{u}^N,
\overline{v}^{(N)})$ of the game with $N$ minor agents is a saddle point
control for the limit stochastic differential game defined above:
\be\label{150bb} \overline{Y}^{u,\overline{v}}_s\le \overline{Y}^{\overline{u},
\overline{v}}_s\le \overline{Y}^{\overline{u},v}_s,\, s\in[t,T],\,
P\mbox{-a.s.},\ (u,v)\in\overline{\cal U}\times\overline{\cal V},
\ee
i.e., in particular, it holds

\centerline{$J(u,\overline{v})\le J(\overline{u},\overline{v})\le
J(\overline{u},v),\ (u,v)\in\overline{\cal U}\times\overline{\cal V}.$}
\et

\begin{proof} \underline{Step 1}: $\overline{Y}^{\overline{u},
\overline{v}}_s\le \overline{Y}^{\overline{u},v}_s,\, s\in[t,T],\,
P\mbox{-a.s.},\ v\in\overline{\cal V}.$

\smallskip

\noindent Indeed, given any $v\in\overline{\cal V}$, we have

\smallskip

\centerline{$\begin{array}{rcl}
F(X_{s}^{0},\overline{X}_{s}^{1},\overline{Y}_{s},\overline{Z}_{s}^{0},
\bar{u}_{s},\bar{v}_{s})&=&f_{\bar{v}}(X_{s}^{0},\overline{X}_{s}^{1},
\overline{Y}_{s},\overline{Z}_{s}^{0},0,\bar{u}_{s})+b_{0}(X_{s}^{0},
\overline{X}_{s}^{1},\overline{Z}_{s}^{0})\overline{Z}_{s}^{0}\bar{u}_{s}\\
&\leq &f(X_{s}^{0},\overline{X}_{s}^{1},\overline{Y}_{s},\overline{Z}_{s}^{0},0,
\bar{u}_{s},v_{s})+b_{0}(X_{s}^{0},\overline{X}_{s}^{1},\overline{Z}_{s}^{0})
\overline{Z}_{s}^{0}\bar{u}_{s}\\
&=&F(X_{s}^{0},\overline{X}_{s}^{1},\overline{Y}_{s},\overline{Z}_{s}^{0},
\bar{u}_{s},v_{s}),\ \ s\in [t,T],
\end{array}$}

\noindent and, thus,

\centerline{$E[F(X_{s}^{0},\overline{X}_{s}^{1},\overline{Y}_{s},
\overline{Z}_{s}^{0},\bar{u}_{s},\bar{v}_{s})|\mathcal{F}_{T}^{W^{0}}]
\leq E[F(X_{s}^{0},\overline{X}_{s}^{1},\overline{Y}_{s},\overline{Z}_{s}^{0},
\bar{u}_{s},v_{s})|\mathcal{F}_{T}^{W^{0}}],\ s\in[t,T].$}

\noindent  This allows to conclude with the help of the comparison
theorem for BSDEs.

\smallskip

\noindent \underline{Step 2}: $\overline{Y}^{u,\overline{v}}_s\le
\overline{Y}^{\overline{u},\overline{v}}_s,\ s\in[t,T],\ P\mbox{-a.s.},
\ u\in\overline{\cal U}.$

\noindent Let $u\in \overline{\mathcal {U}}.$ Then, using the definition
$(\overline{u}_s)_{0\le s\le T}$ and $(\overline{v}_s)_{0\le s\le T}$\ in (\ref{saddlepoint1}) and (\ref{a1a1}), respectively, we have
\be\label{4.183}\begin{array}{rcl}
F(X_{s}^{0},\overline{X}_{s}^{1},\overline{Y}_{s},\overline{Z}_{s}^{0},u_{s},
\bar{v}_{s})&=& F(X_{s}^{0},\overline{X}_{s}^{1},\overline{Y}_{s},\overline{Z}_{s}^{0},
u_{s},\bar{v}(X_{s}^{0}, \overline{X}_{s}^{1},\overline{Y}_{s},\overline{Z}_{s}^{0},
\bar{u}(s,X_{s}^{0},\overline{Y}_{s},\overline{Z}_{s}^{0}))) \\
&=&F_{\bar{v}}(s, X_{s}^{0}, \overline{X}_{s}^{1}, \overline{Y}_{s}, \overline{Z}_{s}^{0}, u_{s}),
\end{array}\ee
with
$$F_{\bar{v}}(s, x_{0}, x_{1}, y, z_{0}, u):=F(x_{0}, x_{1}, y, z_{0}, u, \bar{v}
(x_{0},x_{1},y,z_{0},\bar{u}(s,x_0,y,z_0)).$$

\noindent From the properties of $f,\ b_{0}$ and those of $\bar{u},\ \bar{v}$
we obtain
\be\label{4.185}\begin{array}{lll}
&&(1)\ |F_{\bar{v}}(s,\xi,u)-F_{\bar{v}}(s,\xi^{'}, u)|\leq C(1+|u|)
|\xi-\xi^{'}|,\ \  \xi=(x_{0},x_{1},y,z),\ \xi^{'}=(x_{0}^{'},x_{1}^{'},
y^{'},z^{'}),\\
&&(2)\ |F_{\bar{v}}(s,\xi,u)|\leq C(1+|u|^{2}).
\end{array}\ee
Putting
\be\label{4.186} R_{s}^{N}(x_{0},y,z_{0},u):=|E[F_{\bar{v}}(s,x_{0},
\overline{X}_{s}^{1},y,z_{0},u)|\mathcal{F}_{T}^{W^{0}}]
-\frac{1}{N}\sum_{\ell=1}^{N}F_{\bar{v}}(s,x_{0},X_{s}^{\ell},y,z_{0},u)|,
\ee
we have from Lemma \ref{lemma3.47}
\be\label{4.187}E[|R_{s}^{N}(x_{0},y,z_{0},u)|^{2}|\mathcal{F}_{T}^{W^{0}}]
\leq C_{K}(\frac{1}{N}+\frac{1}{N}\sum_{\ell=1}^{N}|x_{\ell}-\bar{x}|^{2}),
\ee
for all $N\geq1,\ s\in [t,T],\ x_{0},\ y,\ z_{0}\in {\mathbb R}$ and $u\in K$, where
$K\subset U$\ is an arbitrary compact subset of $U$. Consequently, for all
$u\in \overline{\mathcal {U}}$,
\be\label{4.188} E[|R_{s}^{N}(X_{s}^{0},\overline{Y}_{s},\overline{Z}_{s}^{0},u_{s})
|^{2}|\mathcal {F}_{T}^{W^{0}}]\leq C_{u}(\frac{1}{N}+\frac{1}{N}|x_{l}-
\bar{x}|^{2}),
\ee
and with the notation
\be\label{4.189}R_{s}^{N,u}:=R_{s}^{N}(X_{s}^{0},\overline{Y}_{s},\overline{Z}_{s}^{0},u_{s}),
\ee
it follows from (\ref{4.183}) and (\ref{4.186}) that
\be\label{4.190}\begin{array}{lll}
&&E[F(X_{s}^{0},\overline{X}_{s}^{1},\overline{Y}_{s},\overline{Z}_{s}^{0},u_{s},
\bar{v}_{s})| \mathcal{F}_{T}^{W^{0}}] \leq \frac{1}{N}\sum_{\ell=1}^{N}
F_{\bar{v}}(s,X_{s}^{0},X_{s}^{\ell},\overline{Y}_{s},\overline{Z}_{s}^{0},u_{s})+
R_{s}^{N,u}\\
&=&\frac{1}{N}\sum_{\ell=1}^{N}F(X_{s}^{0},X_{s}^{\ell},\overline{Y}_{s},
\overline{Z}_{s}^{0},u_{s}, \bar{v}(X_{s}^{0},X_{s}^{\ell},\overline{Y}_{s},
\overline{Z}_{s}^{0},\bar{u}(s,X_{s}^{0},\overline{Y}_{s},\overline{Z}_{s}^{0})))+
R_{s}^{N,u}\\
&=&\frac{1}{N}\sum_{\ell=1}^{N}F(X_{s}^{0},X_{s}^{\ell},\overline{Y}_{s},
\overline{Z}_{s}^{0},u_{s}, \bar{v}(X_{s}^{0},X_{s}^{\ell},\overline{Y}_{s},
\overline{Z}_{s}^{0},\bar{u}_{s}))
+R_{s}^{N,u}\\
&\leq&\frac{1}{N}\sum_{\ell=1}^{N}F(X_{s}^{0},X_{s}^{\ell},\overline{Y}_{s},
\overline{Z}_{s}^{0},u_{s}, \bar{v}(X_{s}^{0},X_{s}^{\ell},\overline{Y}_{s},
\overline{Z}_{s}^{0},\bar{u}_{s}^{N,0}))+R_{s}^{N,u}+C_{u}|\bar{u}_{s}-
\bar{u}_{s}^{N,0}|,
\end{array}\ee
where

\centerline{$\bar{u}_{s}^{N,0}:=\bar{u}^N(X^{(N)}_s,\overline{Y}_s^N,
(\overline{Z}_s^{0,N},0,\cdots,0)).$}

\noindent The latter estimate was obtained thanks to the Lipschitz
continuity of $F(x_0,x_1,y, z_0,u,v)$ in $v$ and that of $\bar{v}(x_0,
x_1,y,z_0,u)$ in $u$.

\medskip

Thus, using the uniform Lipschitz continuity of $f$ and $b_0(x_0,x_1,z_0)z_0$ in $(y,z_0)$, we get
\be\label{4.190bis}\begin{array}{lll}
&&E[F(X_{s}^{0},\overline{X}_{s}^{1},\overline{Y}_{s},\overline{Z}_{s}^{0},u_{s},
\bar{v}_{s})| \mathcal{F}_{T}^{W^{0}}]\\
&\leq&\frac{1}{N}\sum_{\ell=1}^{N}F(X_{s}^{0},X_{s}^{\ell},\overline{Y}_{s},
\overline{Z}_{s}^{0},u_{s}, \bar{v}(X_{s}^{0},X_{s}^{\ell},\overline{Y}_{s},
\overline{Z}_{s}^{0},\bar{u}_{s}^{N,0}))+R_{s}^{N,u}+C_{u}|\bar{u}_{s}-
\bar{u}_{s}^{N,0}|\\
&\leq&\frac{1}{N}\sum_{\ell=1}^{N}F(X_{s}^{0},X_{s}^{\ell},\overline{Y}_{s}^{N},
\overline{Z}_{s}^{0,N},u_{s},\bar{v}(X_{s}^{0},X_{s}^{\ell},\overline{Y}_{s}^{N},
\overline{Z}_{s}^{0,N},\bar{u}_{s}^{N,0}))+R_{s}^{N,u}\\
& &\ \ +C(|\overline{Y}_{s}^{N}-\overline{Y}_{s}|+|\overline{Z}_{s}^{0,N}-\overline{Z}_{s}^{0}|)
+C_{u}|\bar{u}_{s}-\bar{u}_{s}^{N,0}|.\\
\end{array}\ee

\noindent Observe that from (\ref{vtilde1}) we have
\be\label{4.191}\begin{array}{lll}
\bar{v}_{\ell,s}^{N,0}&:=&\bar{v}_{\ell}^{N}(X_{s}^{(N)},\overline{Y}_{s}^{N},
(\overline{Z}_{s}^{0,N},0,\cdots,0)))\\
&=&\widetilde{v}_N(X_s^0,X_{s}^{\ell},\overline{Y}_{s}^{N},(\overline{Z}_{s}^{0,N},0,
\cdots,0),\bar{u}_s^{N,0})\\
&=&\bar{v}(X_{s}^{0},X_{s}^{\ell},\overline{Y}_{s}^{N},\overline{Z}^{0,N}_s,
\bar{u}_{s}^{N,0}).\\
\end{array}\ee
The latter equality follows from the fact that, if $z_\ell=0,\, 1\le\ell
\le N$, then the minimizer $\overline{v}(\xi,u)$ of $f(\xi,u,.)$ (see
(\ref{a11})) and the minimizer $\widetilde{v}_N(\xi,u)$ in
(\ref{4.7a}) satisfy the relation
$$\widetilde{v}_N(x_0,x_\ell,y,(z_0,0,\dots,0),u)=\overline{v}
(x_0,x_\ell,y,z_0,u).$$

\smallskip

\noindent Therefore, using definition (\ref{4.169}) of $F$, we can
write
\be\label{4.192}\begin{array}{lll}
&&E[F(X_{s}^{0},\overline{X}_{s}^{1},\overline{Y}_{s},\overline{Z}_{s}^{0},u_{s},
\bar{v}_{s})|\mathcal{F}_{T}^{W^{0}}]\\
&\leq& (\frac{1}{N}\sum_{\ell=1}^{N}f(X_{s}^{0},X_{s}^{\ell},
\overline{Y}_{s}^{N},\overline{Z}_{s}^{0,N}, 0,u_{s},\bar{v}_{\ell,s}^{N,0})
+\frac{1}{N}\sum_{\ell=1}^{N}b_{0}(X_{s}^{0},X_{s}^{\ell},
\overline{Z}_{s}^{0,N})\overline{Z}_{s}^{0,N}.u_{s})\\
& &+R_{s}^{N,u}+C(|\overline{Y}_{s}^{N}-\overline{Y}_{s}|+|\overline{Z}_{s}^{0,N}
-\overline{Z}_{s}^{0}|)+C_{u}|\bar{u}_{s}-\bar{u}_{s}^{N,0}|\\
&=&H_{N}(X_{s}^{(N)},\overline{Y}_{s}^{N},(\overline{Z}_{s}^{0,N},0,\cdots,0),
u_{s},\bar{v}_{s}^{N,0}) +R_{s}^{N,u}+C(|\overline{Y}_{s}^{N}-\overline{Y}_{s}|
+|\overline{Z}_{s}^{0,N}-\overline{Z}_{s}^{0}|) +C_{u}|\bar{u}_{s}-\bar{u}_{s}^{N,0}|\\
&\leq &H_{N}(X_{s}^{(N)},\overline{Y}_{s}^{N},(\overline{Z}_{s}^{0,N},0,\cdots,0),
\bar{{u}}_{s}^{0,N},\bar{v}_{s}^{N,0}) +R_{s}^{N,u}+C(|\overline{Y}_{s}^{N}-
\overline{Y}_{s}|+|\overline{Z}_{s}^{0,N}-\overline{Z}_{s}^{0}|) +C_{u}|\bar{u}_{s}-
\bar{u}_{s}^{N,0}|.\\
\end{array}\ee

Here we have used that $(\bar{u}_{s}^{N},\ \bar{v}_{s}^{N,0}=
(\bar{v}_{1,s}^{N,0},\dots,\bar{v}_{N,s}^{N,0}))$ is the saddle point
control for the Hamiltonian $H_{N}(X_{s}^{(N)},\overline{Y}_{s}^{N},
(\overline{Z}_{s}^{0,N},0,\dots,0),\cdot,\cdot)$. On the other hand,
since the controls $\overline{u}^{N,0}$ and $\overline{u}$ are bounded
by some constant not depending on $N$ (see Lemma \ref{lemma4.8} and take
into account that here $z_\ell=0,\, 1\le \ell\le N$), we obtain
\be\label{4.193}\begin{array}{lll}
&&H_{N}(X_{s}^{(N)},\overline{Y}_{s}^{N},(\overline{Z}_{s}^{0,N},0,\cdots,0),
\bar{u}_{s}^{N,0},\bar{v}_{s}^{N,0})\\
&=&\frac{1}{N}\sum_{\ell=1}^{N}f_{\bar{v}}(X_{s}^{0},X_{s}^{\ell},
\overline{Y}_{s}^{N},\overline{Z}_{s}^{0,N},0, \bar{u}_{s}^{N,0})+\frac{1}{N}
\sum_{\ell=1}^{N}b_{0}(X_{s}^{0},X_{s}^{\ell},\overline{Z}_{s}^{0,N})
\overline{Z}_{s}^{0,N}.\bar{u}_{s}^{N,0}\\
&\leq&\frac{1}{N}\sum_{\ell=1}^{N}f_{\bar{v}}(X_{s}^{0},X_{s}^{\ell},
\overline{Y}_{s}, \overline{Z}_{s}^{0},0,\bar{u}_{s})+\frac{1}{N
}\sum_{\ell=1}^{N}b_{0}(X_{s}^{0},X_{s}^{\ell},
\overline{Z}_{s}^{0})\overline{Z}_{s}^{0}.\bar{u}_{s}\\
&&+C(|\overline{Y}_{s}^{N}-\overline{Y}_{s}|+|\overline{Z}_{s}^{0,N}-\overline{Z}_{s}^{0}|)
+C|\bar{u}_{s}^{N,0}-\bar{u}_{s}|\\
&\leq&E[f_{\bar{v}}(X_{s}^{0},\overline{X}_{s}^{1},\overline{Y}_{s},
\overline{Z}_{s}^{0},0,\bar{u}_{s})+b_{0}(X_{s}^{0},\overline{X}_{s}^{1},
\overline{Z}_{s}^{0})\overline{Z}_{s}^{0}.\bar{u}_{s}|\mathcal{F}_{T}^{W^{0}}]
+C(|\overline{Y}_{s}^{N}-\overline{Y}_{s}|+|\overline{Z}_{s}^{0,N}-\overline{Z}_{s}^{0}|)\\
& &+C|\bar{u}_{s}^{N,0}-\bar{u}_{s}|+R_{s}^{N}\\
&=&E[F(X_{s}^{0},\overline{X}_{s}^{1},\overline{Y}_{s},\overline{Z}_{s}^{0},
\bar{u}_{s},\bar{v}_{s}) |\mathcal{F}_{T}^{W^{0}}]+R_{s}^{N}
+C(|\overline{Y}_{s}^{N}-\overline{Y}_{s}|+|\overline{Z}_{s}^{0,N}-\overline{Z}_{s}^{0}|)+
C|\bar{u}_{s}^{N,0}-\bar{u}_{s}|,
\end{array}\ee
where $E[(R_{s}^{N})^{2}|\mathcal{F}_{T}^{W^{0}}]\leq C(\frac{1}{N}+
\frac{1}{N}\sum_{\ell=1}^{N}|x_{\ell}-\bar{x}|^{2}),$\ and we have used Lemma \ref{lemma3.47} for the latter estimate.
Hence, combining the above estimates of Step 2, we get
\be\label{4.194}\begin{array}{lll}
& & E[F(X_{s}^{0},\overline{X}_{s}^{1},\overline{Y}_{s},\overline{Z}_{s}^{0},
u_{s},\bar{v}_{s})|\mathcal{F}_{T}^{W^{0}}]\\
&\leq & E[F(X_{s}^{0},\overline{X}_{s}^{1},\overline{Y}_{s},\overline{Z}_{s}^{0},
\bar{u}_{s}, \bar{v}_{s})|\mathcal{F}_{T}^{W^{0}}]
+C(|\overline{Y}_{s}^{N}-\overline{Y_{s}}|+|\overline{Z}_{s}^{0,N}-\overline{Z}_{s}^{0}|)
+C|\bar{u}_{s}^{N,0}-\bar{u}_{s}|+R_{s}^{N},\end{array}\ee
where

\centerline{$E[(R_{s}^{N})^{2}|\mathcal{F}_{T}^{W^{0}}]\leq C(
\frac{1}{N}+\frac{1}{N}\sum_{\ell=1}^{N}|x_{\ell}-\bar{x}|^{2}),\
N\geq1,\ s\in [t,T].$}

\smallskip

\noindent But, since due to Lemma \ref{auxlemma11} (recall
that $\overline{u}^{N,0}_s=\overline{u}^N(X_s^{(N)},\overline{Y}^N_s,
(\overline{Z}^{0,N}_s,0\dots,0))$),

\smallskip

\centerline{$|\bar{u}_{s}^{N,0}-\bar{u}_{s}|\le C(|\overline{Y}_{s}^{N}
-\overline{Y_{s}}|+|\overline{Z}_{s}^{0,N}-\overline{Z}_{s}^{0}|)+R_{s}^{N}.$}

\noindent we conclude
\be\label{4.194}\begin{array}{lll}
& & E[F(X_{s}^{0},\overline{X}_{s}^{1},\overline{Y}_{s},\overline{Z}_{s}^{0},
u_{s},\bar{v}_{s})|\mathcal{F}_{T}^{W^{0}}]\\
&\leq & E[F(X_{s}^{0},\overline{X}_{s}^{1},\overline{Y}_{s},\overline{Z}_{s}^{0},
\bar{u}_{s}, \bar{v}_{s})|\mathcal{F}_{T}^{W^{0}}]
+C(|\overline{Y}_{s}^{N}-\overline{Y_{s}}|+|\overline{Z}_{s}^{0,N}-\overline{Z}_{s}^{0}|)
+{R}_{s}^{N},\end{array}\ee
and, thanks to Theorem \ref{Theorem 4.1}, as $N\rightarrow \infty$,
\be\label{4.195}E[F(X_{s}^{0},\overline{X}_{s}^{1},\overline{Y}_{s},
\overline{Z}_{s}^{0},u_{s}, \bar{v}_{s})|\mathcal{F}_{T}^{W^{0}}]\leq
E[F(X_{s}^{0},\overline{X}_{s}^{1},\overline{Y}_{s},\overline{Z}_{s}^{0},
\bar{u}_{s},\bar{v}_{s}) |\mathcal{F}_{T}^{W^{0}}],\  dsdP\mbox{-a.e.}\ee

\noindent Consequently, from the comparison theorem for BSDEs we have
$$\overline{Y}_{s}^{\bar{u},\bar{v}}=\overline{Y}_{s}\geq \overline{Y}_{s}^{u,\bar{v}},
\  s\in [t,T],\ u\in \overline{\mathcal {U}},$$
and the proof is complete now.
\end{proof}

 While we have seen in Theorem \ref{theorem150bb} that
$(\overline{u},\overline{v})$ is a saddle point control of our limit
2-person zero-sum game, it turns out that it is even the unique
saddle point control. Indeed, we have the following uniqueness
result:

\bt Let $(u',v'),\ (u'',v'')\in\overline{\cal U}\times\overline{\cal V}$
be two couples of saddle point controls in the sense that, for all
$(u,v)\in\overline{\cal U}\times\overline{\cal V},$
$$\begin{array}{llllll}
&J(u,v')&\le& J(u',v')&\le& J(u',v),\\
&J(u,v'')&\le& J(u'',v'')&\le& J(u'',v).
\end{array}$$
\noindent Then $(u'_s,\ v'_s)=(u''_s,\ v''_s),$ $dsdP$-a.e.
\et
\begin{proof} Let $(u',v'),\ (u'',v'')\in\overline{\cal U}\times
\overline{\cal V}$ be two couples of saddle point controls in the
above sense.

\noindent\underline{Step 1}. Let $u\in{\cal U}$. Then $J(u',v')
\ge J(u,v')$, $P$-a.s., and it follows that $\overline{Y}^{u',v'}_s\ge\overline{Y}^{u,v'}_s,\ s\in[t,T],$ and
$$E[F(X_s^0,\overline{X}^1_s,\overline{Y}^{u',v'}_s,\overline{Z}^{u',v'}_s,u'_s,v'_s)|{\cal F}_T^{W^0}]\ge E[F(X_s^0,
\overline{X}^1_s,\overline{Y}^{u',v'}_s,\overline{Z}^{u',v'}_s,u_s,v'_s)|
{\cal F}_T^{W^0}], \mbox{dsdP-a.e.}$$

\noindent Indeed, putting $\gamma_s:=I\{E[F(X_s^0,\overline{X}^1_s,
\overline{Y}^{u',v'}_s, \overline{Z}^{u',v'}_s,u'_s,v'_s)|{\cal F}_T^{W^0}]
\ge E[F(X_s^0,\overline{X}^1_s,\overline{Y}^{u',v'}_s, \overline{Z}^{u',
v'}_s,u_s,v'_s)|{\cal F}_T^{W^0}]\}$, the process $\tilde{u}_s:=u'_s\gamma_s
+u_s(1-\gamma_s),\ s\in[t,T],$ defines an admissible control in
$\overline{\cal U},$ and
$$E[F(X_s^0,\overline{X}^1_s,\overline{Y}^{u',v'}_s,\overline{Z}^{u',v'}_s,\tilde{u}_s,v'_s)|{\cal F}_T^{W^0}]\ge
E[F(X_s^0,\overline{X}^1_s,\overline{Y}^{u',v'}_s,\overline{Z}^{u',v'}_s,u'_s,v'_s)|{\cal F}_T^{W^0}],\ s\in[t,T].$$

\noindent Consequently, the comparison theorem yields $\overline{Y}_s^{\tilde{u},v'}\ge \overline{Y}_s^{u',v'},$\ P-a.s., $s\in[t,T]$. But, since on the other hand,
$J(\tilde{u},
v')\le J(u',v'),$ it follows from the converse comparison theorem
that $\overline{Y}^{\tilde{u},v'}_s= \overline{Y}^{u',v'}_s,$\ P-a.s., $s\in[t,T],$ and
$$E[F(X_s^0,\overline{X}^1_s,\overline{Y}^{u',v'}_s,
\overline{Z}^{u',v'}_s,\tilde{u}_s, v'_s)|{\cal F}_T^{W^0}]=
E[F(X_s^0,\overline{X}^1_s,\overline{Y}^{u',v'}_s, \overline{Z}^{u',
v'}_s,u'_s, v'_s)|{\cal F}_T^{W^0}],\ \mbox{dsdP-a.e.}$$

\noindent But this implies
$$E[F(X_s^0,\overline{X}^1_s,\overline{Y}^{u',v'}_s,
 \overline{Z}^{u',v'}_s,u_s, v'_s)|{\cal F}_T^{W^0}]\le E[F(X_s^0,
 \overline{X}^1_s,\overline{Y}^{u',v'}_s, \overline{Z}^{u',v'}_s,u'_s,
 v'_s)|{\cal F}_T^{W^0}],\ \mbox{dsdP-a.e.}, $$

\noindent from where we see that $\overline{Y}^{u',v'}_s\ge
\overline{Y}^{u,v'}_s,\ s\in[t,T],$ and $u'_s$ is a maximum point
of $u\rightarrow E[F(X_s^0, \overline{X}^1_s,\overline{Y}^{u',v'}_s,$ $
\overline{Z}^{u',v'}_s,u, v'_s)|{\cal F}_T^{W^0}]$, dsdP-a.e. On the
other hand, let us also observe that, from the definition of $F$ and the
assumptions on $f$,
$$\begin{array}{rcl}& &\(D_uE[F(X_s^0,\overline{X}^1_s,\overline{Y}^{u',v'}_s,
\overline{Z}^{u',v'}_s,u^1, v'_s)|{\cal F}_T^{W^0}]-D_u
E[F(X_s^0,\overline{X}^1_s,\overline{Y}^{u',v'}_s,
\overline{Z}^{u',v'}_s,u^2, v'_s)|{\cal F}_T^{W^0}], u^1-u^2\)\\
& &\le -\lambda|u^1-u^2|^2,\ u^1, u^2\in U.\end{array}$$

\noindent Consequently, $u'_s$ is the unique maximum point of
$u\rightarrow E[F(X_s^0,\overline{X}^1_s,\overline{Y}^{u',v'}_s,
\overline{Z}^{u',v'}_s,u,v'_s)|{\cal F}_T^{W^0}]$.

\medskip

\noindent\underline{Step 2}. Let $v\in\overline{\cal V}.$\ Following the argument in Step 1, but with putting
$$\gamma_s:=I\{F(X_s^0,\overline{X}^1_s,
\overline{Y}^{u',v'}_s, \overline{Z}^{u',v'}_s,u'_s,v'_s)
\le F(X_s^0,\overline{X}^1_s,\overline{Y}^{u',v'}_s,
\overline{Z}^{u',v'}_s,u'_s,v_s)\}$$
\noindent and $\tilde{v}_s:=v'_s\gamma_s+v_s(1-\gamma_s),
\ s\in[t,T],$ we have $\tilde{v}\in\overline{\cal V},$ and
$$F(X_s^0,\overline{X}^1_s,\overline{Y}^{u',v'}_s,
\overline{Z}^{u',v'}_s,u'_s,\tilde{v}_s)\le F(X_s^0,
\overline{X}^1_s,\overline{Y}^{u',v'}_s, \overline{Z}^{u',v'}_s,
u'_s,v'_s),\ s\in[t,T].$$

\noindent Consequently, from the comparison theorem it follows that
$\overline{Y}^{u',\tilde{v}}_s\le \overline{Y}^{u',v'}_s,$\ P-a.s., $s\in[t,T],$ and since
on the other hand, $J(u',\tilde{v})\ge J(u',v'),$ the converse comparison theorem implies $\overline{Y}^{u',\tilde{v}}_s= \overline{Y}^{u',v'}_s,$\ P-a.s., $s\in[t,T],$ and
$$E[F(X_s^0,\overline{X}^1_s,\overline{Y}^{u',v'}_s,
\overline{Z}^{u',v'}_s,u'_s, \tilde{v}_s)|{\cal F}_T^{W^0}]=
E[F(X_s^0,\overline{X}^1_s,\overline{Y}^{u',v'}_s,
\overline{Z}^{u',v'}_s,u'_s, v'_s)|{\cal F}_T^{W^0}],\ \mbox{dsdP-a.e.}$$

\noindent Thus,

\centerline{$F(X_s^0,\overline{X}^1_s,\overline{Y}^{u',v'}_s,
\overline{Z}^{u',v'}_s,u'_s, \tilde{v}_s)= F(X_s^0,
\overline{X}^1_s,\overline{Y}^{u',v'}_s, \overline{Z}^{u',v'}_s,
u'_s, v'_s),$\ dsdP-a.e.,}

\noindent and

\centerline{$F(X_s^0,\overline{X}^1_s,\overline{Y}^{u',v'}_s,
\overline{Z}^{u',v'}_s,u'_s, v_s)\ge F(X_s^0,\overline{X}^1_s,
\overline{Y}^{u',v'}_s, \overline{Z}^{u',v'}_s,u'_s, v'_s),$\ dsdP-a.e.,}

\noindent which implies that $\overline{Y}^{u',v'}_s\le
\overline{Y}^{u',v}_s,\ s\in[t,T],$\ but also that $v'_s$ is a
minimum point of $v\rightarrow F(X_s^0,\overline{X}^1_s,$ $
\overline{Y}^{u',v'}_s, \overline{Z}^{u',v'}_s, u'_s, v).$
On the other hand, since
$$\begin{array}{lll}& &\(D_vF(X_s^0,\overline{X}^1_s,\overline{Y}^{u',v'}_s,
\overline{Z}^{u',v'}_s,u', v^1)-D_vF(X_s^0,\overline{X}^1_s,
\overline{Y}^{u',v'}_s, \overline{Z}^{u',v'}_s,u', v^2),\ v^1-v^2\)\\
&&\ge \lambda|v^1-v^2|^2,\ v^1,\ v^2\in V,\end{array}$$

\noindent it follows that $v'_s$ is the unique minimum point of $v\rightarrow F(X_s^0,
\overline{X}^1_s,\overline{Y}^{u',v'}_s, \overline{Z}^{u',
v'}_s,u'_s, v).$

\medskip

\noindent\underline{Step 3}.  Let $(u',v'),\ (u'',v'')\in
\overline{\cal U}\times \overline{\cal V}$ be two couples of
saddle point controls. Then, combining our results from the
Steps 1 and 2 we have

\centerline{$\overline{Y}^{u',v'}_s\ge \overline{Y}^{u'',v'}_s\ge \overline{Y}^{u'',v''}_s\ge \overline{Y}^{u',v''}_s
\ge \overline{Y}^{u',v'}_s,\ s\in[t,T],$}  i.e.,

\centerline{$\overline{Y}^{u',v'}_s= \overline{Y}^{u'',v'}_s=\overline{Y}^{u'',v''}_s=\overline{Y}^{u',v''}_s=
\overline{Y}^{u',v'}_s,\ s\in[t,T].$}

\noindent But from $\overline{Y}^{u',v'}_s= \overline{Y}^{u'',v'}_s,\ s\in[t,T],$\
and the uniqueness of the semimartingale decomposition of this
process we deduce that $\overline{Z}^{u',v'}_s= \overline{Z}^{u'',v'}_s$\, dsdP-a.e.,
and also
$$E[F(X_s^0,\overline{X}^1_s,\overline{Y}^{u',v'}_s,
\overline{Z}^{u',v'}_s,u''_s, v'_s)|{\cal F}_T^{W^0}]= E[F(X_s^0,
\overline{X}^1_s,\overline{Y}^{u',v'}_s, \overline{Z}^{u',v'}_s,
u'_s, v'_s)|{\cal F}_T^{W^0}],\ \mbox{dsdP-a.e.}$$

\noindent But this means that also $u''_s$ is a maximum point
of $E[F(X_s^0,\overline{X}^1_s,\overline{Y}^{u',v'}_s,\overline{Z}^{u',v'}_s,., v'_s)|{\cal F}_T^{W^0}]$, dsdP-a.e.,
and, thus, due to Step 1, $u'_s=u''_s,$\ dsdP-a.e.

Let us now show that also $v'_s=v''_s$, $dsdP$-a.e. As we have seen
already above, $\overline{Y}^{u',v'}_s=\overline{Y}^{u'',v''}_s,\
s\in[t,T]$, and from the uniqueness of the semimartingale
decomposition of this process we get $\overline{Z}^{u',v'}_s= \overline{Z}^{u'',v''}_s,$\
dsdP-a.e., and
$$E[F(X_s^0,\overline{X}^1_s,\overline{Y}^{u',v'}_s,
\overline{Z}^{u',v'}_s,u''_s, v''_s)|{\cal F}_T^{W^0}]= E[F(X_s^0,
\overline{X}^1_s,\overline{Y}^{u',v'}_s, \overline{Z}^{u',v'}_s,u',
v'_s)|{\cal F}_T^{W^0}],\ \mbox{dsdP-a.e.}$$

\noindent Since, on the other hand, for $v=v''$ in Step 2,
$$F(X_s^0,\overline{X}^1_s,\overline{Y}^{u',v'}_s,
\overline{Z}^{u',v'}_s,u'_s, v''_s)\ge F(X_s^0,\overline{X}^1_s,
\overline{Y}^{u',v'}_s, \overline{Z}^{u',v'}_s,u'_s, v'_s),\ \mbox{dsdP-a.e.},$$

\noindent it follows that

$$F(X_s^0,\overline{X}^1_s,\overline{Y}^{u',v'}_s,
\overline{Z}^{u',v'}_s,u'_s, v''_s)= F(X_s^0,\overline{X}^1_s,
\overline{Y}^{u',v'}_s, \overline{Z}^{u',v'}_s,u'_s, v'_s),\ \mbox{dsdP-a.e.},$$

\noindent i.e., also $v''_s$ is a minimum point of $F(X_s^0,
\overline{X}^1_s,\overline{Y}^{u',v'}_s, \overline{Z}^{u',v'}_s,
u'_s,.)$, dsdP-a.e. But due to Step 2, the unique minimum point
is $v'_s$. Consequently, $v'_s=v''_s,$\ dsdP-a.e.

Since we have already shown that $u'_s=u''_s,$\ dsdP-a.e., we
conclude that the saddle point controls $(u',v')$ and $(u'',v'')$
coincide.
\end{proof}

\section{Appendix}
\subsection{Appendix 1}
Appendix 1 is devoted to the proof of Proposition $\ref{p3}.$

\begin{proof} (of Proposition \ref{p3}). Without loss of the
generality, we may suppose that $b_0=0$ and $b_1=0$.
We note that with this convention and with the notation
$$\sigma^{1,j,N}_u:=\frac{1}{N}\sum_{\ell=1}^N
\sigma_1(X_u^{0,N},X_u^{j,N},X_u^{\ell,N})
-E[\sigma_1(\overline{X}_u^{0},\overline{X}_u^{j},
\overline{X}_u^{j+1})|{\cal F}_T^{W^0,W^j}],$$

\noindent for $t\le r\le s\le T,$
$$X_r^{j,N}-\overline{X}^j_r=(x_j-\overline{x})+
\int_t^r\sigma^{1,j,N}_udW_u^j,$$

\noindent from where we deduce
\be\label{annexe10}\begin{array}{lll} & &\frac{1}{N}\sum_{j=1}^N
|X_r^{j,N}-\overline{X}^j_r|^2\\
&\le&\frac{2}{N}\sum_{j=1}^N|x_j-\overline{x}|^2+\frac{2}{N}
\sum_{j=1}^N|\int_t^r\sigma^{1,j,N}_u d W_u^j|^2\\
&=&\frac{2}{N}\sum_{j=1}^N|x_j-\overline{x}|^2+2|M^{(N)}_r|^2,
\end{array}\ee
where
$$M^{(N)}_r=\sum_{j=1}^Ne_j\int_t^r\frac{1}{\sqrt{N}}
\sigma^{1,j,N}_udW_u^j,\ r\in[t,T],$$

\noindent is an $\ell_2({\mathbb R}^d)$-valued $\mathbb{F}$-martingale, and
$e_j=(\delta_{j,k})_{k\ge 1}$ is the element of $\ell_2$ with
$\delta_{j,j}=1$ and $\delta_{j,k}=0,\ k\not=j.$ Consequently,
due to the Burkholder-Davis-Gundy Inequality, for some constant
$C_m\in\mathbb{R}$ which can vary from line to line but doesn't
depend on $N$,
\be\label{annexe11}\begin{array}{lll} & &E[\sup_{r\in[t,s]}
(\frac{1}{N}\sum_{j=1}^N|X_r^{j,N}-\overline{X}^j_r|^2
)^m]\\
&\le& C_m(\frac{1}{N}\sum_{j=1}^N|x_j-\overline{x}|^2)^m+C_m E[\sup_{r\in[t,s]}|M^{(N)}_r|^{2m}]\\
&\le& C_m(\frac{1}{N}\sum_{j=1}^N|x_j-\overline{x}|^2)^m+C_m E[\langle M^{(N)}\rangle_s^{m}],
\end{array}\ee
where
$$\langle M^{(N)}\rangle_s=\frac{1}{N}\sum_{j=1}^N
\int_t^s|\sigma^{1,j,N}_u|^2du,\ s\in[t,T],$$

\noindent is the quadratic variation process of the martingale
$M^{(N)}=(M^{(N)}_s)_{s\in[t,T]}.$\ Hence, with the notation
\be\label{3.27}\zeta_{1,j,\ell}(r):=\sigma_1(\overline{X}^0_r,
\overline{X}_r^j, \overline{X}_r^\ell)-E[
\sigma_1(\overline{X}^0_r,\overline{X}_r^j,\overline{X}_r^\ell)|
{\cal F}^{W^0,W^j}_T],\ee
we have
\be\label{annexe12}\begin{array}{lll} & &E[\sup_{r\in[t,s]}
(\frac{1}{N}\sum_{j=1}^N|X_r^{j,N}-\overline{X}^j_r|^2
)^m]\\
&\le&C_m(\frac{1}{N}\sum_{j=1}^N|x_j-\overline{x}|^2
)^m+C_m E[(\frac{1}{N}\sum_{j=1}^N\int_t^s
|\sigma^{1,j,N}_r|^2dr)^{m}]\\
&\le&C_m(\frac{1}{N}\sum_{j=1}^N|x_j-\overline{x}|^2)^m\\
& &+C_m E[(\int_t^s|X_r^{0,N}-\overline{X}_r^{0}|^2dr)^m]\hfill(=: J_s^{1,N})\\
& &+C_m E[(\frac{1}{N}\sum_{j=1}^N\int_t^s|
X_r^{j,N}-\overline{X}_r^j|^2dr+\int_t^s(\frac{1}{N}
\sum_{\ell=1}^N| X_r^{\ell,N}-\overline{X}_r^{\ell}|^2) dr)^m]\hfill (=:J_s^{2,N})\\
& &+C_mE[(\frac{1}{N}\sum_{j=1}^N\int_t^s|
\frac{1}{N}\sum_{\ell=1,\ell\neq j}^N\zeta_{1,j,\ell}(r)|^2dr)^m]\hfill (=: J_s^{3,N})\\
& &+\frac{C_m}{N^{2m}}.\\
\end{array}\ee

\noindent The term $J_s^{2,N}$ will be handled by Gronwall's
Lemma. Let us begin with estimating $J_s^{1,N}.$ Obviously,
putting
\be\label{annexe13}\begin{array}{ll}
& \sigma^{0,N}_r:=\frac{1}{N}\sum_{\ell=1}^N
\sigma_0(X_r^{0,N},X_r^{\ell,N})
-E[\sigma_0(\overline{X}_r^{0},\overline{X}_r^{1})|
{\cal F}_T^{W^0}],\quad\mbox{ and }\\
& \zeta_{0,\ell}(r):=\sigma_0(\overline{X}^0_r,
\overline{X}_r^\ell)-E[\sigma_0 (\overline{X}^0_r,
\overline{X}_r^\ell)|{\cal F}^{W^0}_T],\ \ r\in[t,T],
\end{array}\ee
we obtain,
\be\label{annexe13}\begin{array}{lll} & & E[\sup_{r\in[t
,s]}|X_r^{0,N}-\overline{X}^0_r|^{2m}]\le C_m E[
(\int_t^s|\sigma^{0,N}_r|^2dr )^{m}]\\
&\le & C_m E[(\int_t^s|X_r^{0,N}-\overline{X}_r^{0}|^2dr)^m]+ C_m E[(\frac{1}{N}
\sum_{\ell=1}^N\int_t^s|X_r^{\ell,N}-\overline{X}_r^\ell|^2dr)^m]\\
& &+C_m E[(\int_t^s|\frac{1}{N}\sum_{\ell=1}^N
\zeta_{0,\ell}(r)|^2dr)^m]\hfill
(=: J_s^{4,N}).\\
\end{array}\ee
The following estimates of $J^{3,N}_s$ and $J_s^{4,N}$ allow to
complete the proof by applying Grownwall's Lemma to the system
(\ref{annexe12})-(\ref{annexe13}).\end{proof}

\bl\label{lemma_zeta} For some constant $C_m$ independent of
$N\ge 1$, with the notation $\Lambda_{j,N}:=\{1,\dots,N\}
\setminus\{j\}$, we have
\be\label{3.28} E[\|\sum_{\ell=1}^N\zeta_{0,\ell}(r)
\|^{2m}|{\cal F}_T^{W^0}]\le C_mN^m,\, \, E[\|\sum_{\ell\in\Lambda_{j,N}}\zeta_{1,j,\ell}(r)\|^{2m}|
{\cal F}_T^{W^0,W^{j}}]\le C_mN^m,\ N\ge 1.\ee
\el
\begin{proof} We give the proof only for $E[\|\sum_{\ell\in\Lambda_{j,N}}\zeta_{j,\ell}(r)|^{2m}|{\cal F}_T^{W^0,W^{j}}]$, that for
$E[\|\sum_{\ell=1}^N\zeta_{0,\ell}(r)\|^{2m}|{\cal F}_T^{W^0}]$ is analogous. For simplifying the argument we assume
that $d=1$. For the multi-dimensional case $d>1$ the random
variable $\zeta_{1,j,\ell}(r)$ is matrix-valued and the argument
we develop shall be applied component-wise. For our argument
we will need the set
$$\Gamma_{m,N,j}:=\{(\ell_1,\dots,\ell_{2m})\in
\Lambda_{j,N}^{2m}\, |\, \forall i(1\le i\le 2m),\exists i'\in
\{1,\dots,2m\}\setminus\{i\}: \ell_{i}=\ell_{i'}\}.$$

\noindent Remark that the cardinal number of $\Gamma_{m,N,j}$
satisfies the estimate
$$\mbox{card}(\Gamma_{m,N,j})\le C_mN^m,\, \, N\ge 1,$$

\noindent for some $C_m\in R$ not depending on $N$.

We begin the proof of the lemma by remarking that knowing
${\cal F}^{W^0, W^j}_T$ the processes $\overline{X}^\ell,\,
\ell\ge 1$ with $\ell\not=j$, are conditionally i.i.d.
Consequently, also the random variables $\zeta_{1,j,\ell}(r),\,
\ell\in\{1,\dots,N\} \setminus\{j\},$ are conditionally i.i.d.,
knowing ${\cal F}^{W^0, W^j}_T$. Moreover, $E[\zeta_{1,j,
\ell}(r)|{\cal F}^{W^0, W^j}_T]=0$. Thus,
\be\label{3.29}\begin{array}{lll}
& &E[\|\frac{1}{\sqrt{N}}\sum_{\ell\in\Lambda_{j,N}}
\zeta_{1,j,\ell}(r) \|^{2m}|{\cal F}_T^{W^0,W^{j}}]
=\frac{1}{N^m}\sum_{\ell_1,\dots,\ell_{2m}\in\Lambda_{j,N}}
E[ \prod_{p=1}^{2m}\zeta_{1,j,\ell_p}(r)|{\cal F}_T^{W^0,W^{j}}]\\
& & =\frac{1}{N^m}\sum_{(\ell_1,\dots,\ell_{2m})\in\Gamma_{m,
N,j}}E[ \prod_{p=1}^{2m}\zeta_{1,j,\ell_p}(r)|
{\cal F}_T^{W^0,W^{j}}]\\
& &\le \frac{1}{N^m}\sum_{(\ell_1,\dots,\ell_{2m})\in
\Gamma_{m,N,j}}\prod_{p=1}^{2m}(E[\|\zeta_{1,j,
\ell_p}(r)\|^{2m}|{\cal F}_T^{W^0,W^{j}}])^{\frac{1}{2m}}.\\
\end{array}\ee
\noindent In virtue of the boundedness of $\sigma_1$ and the
estimate of the cardinal number of $\Gamma_{m, N, j}$ this
yields
\be\label{3.30} E[\|\frac{1}{\sqrt{N}}\sum_{\ell\in
\Lambda_{j,N}}\zeta_{1,j,\ell}(r) \|^{2m}|
{\cal F}_T^{W^0,W^{j}}]\le C_m,
\ee

\smallskip

\noindent for some $C_m\in {\mathbb R}$ independent of
$N\ge 1.$ The statement of the lemma follows easily from
this.\end{proof}

\subsection{Appendix 2}
This appendix is devoted to the proof of Proposition \ref{BSDE_a}.

\begin{proof} (of Proposition \ref{BSDE_a}). Taking the
difference between the BSDE for $Y^N$ and that for
$\overline{Y}$ we obtain
\be\label{3.53}\begin{array}{lll}
& Y_s^{N}-\overline{Y}_s=\displaystyle\{\frac{1}{N}
\sum_{\ell=1}^N\Phi(X_T^{0,N}, X_T^{\ell,N})-E[\Phi
(\overline{X}_T^0,\overline{X}_T^1)|{\cal{F}}_T^{W^0}]
\} &(:=\eta^{1,N})\\
&\displaystyle +\int_s^T(\frac{1}{N}\sum_{\ell=1}^N
(f(X_r^{0,N},X_r^{\ell,N},Y_r^N,Z_r^{0,N}, Z_r^{\ell,N})
-f(\overline{X}_r^{0},\overline{X}_r^{\ell}, Y_r^N,Z_r^{0,N},
Z_r^{\ell,N})))dr &(:=\int_s^T\eta_r^{2,N}dr)\\
& +\displaystyle\int_s^T(\frac{1}{N}\sum_{\ell=1}^N
(f(\overline{X}_r^{0},\overline{X}_r^{\ell}, Y_r^N,Z_r^{0,N},
Z_r^{\ell,N})-f(\overline{X}_r^{0},\overline{X}_r^{\ell},
\overline{Y}_r,\overline{Z}_r^{0},0)))dr
&(:=\int_s^T\eta_r^{3,N}dr)\\
&+\displaystyle\int_s^T(\frac{1}{N}\sum_{\ell=1}^N
f(\overline{X}_r^{0},\overline{ X}_r^{\ell}, \overline{Y}_r,
\overline{Z}_r^{0},0)-E[f( \overline{X}_r^0,\overline{X}_r^1,
\overline{Y}_r,\overline{Z}_r^0,0)|{\cal F}^{W^0}_T])
dr &(:=\int_s^T\eta_r^{4,N}dr)\\
&-\displaystyle\int_s^T(Z_r^{0,N}-\overline{Z}_r^0)dW_r^0-
\sum_{\ell=1}^N\int_s^TZ_r^{\ell, N}dW_r^\ell,\ s\in [t, T]. &\\
\end{array}\ee
From Proposition \ref{p3} and Lemma \ref{lemma3.47} we know
already that
\be\label{3.54} E[|\eta^{1,N}|^{2m}]\leq C_m (\frac{1}{N}+\frac{1}{N}\sum_{\ell=1}^N|x_\ell-
\overline{x}|^2)^m,\  N\ge 1,\ee

\smallskip

\noindent for some constant $C_m$ independent of $N$.
Furthermore, again from Proposition \ref{p3} we have
\be\label{3.55}\begin{array}{lll}& & E[\int_t^T|
\eta_r^{2,N}|^{2m}dr]\le C_m E[\sup_{r\in [t,T]}
(|X_r^{0,N}-\overline{X}_r^0|^{2}+ \frac{1}{N}
\sum_{\ell=1}^N|X_r^{\ell,N}- \overline{X}_r^\ell|^{2})^{m}]\\
& & \leq C_m(\frac{1}{N}+\frac{1}{N}
\sum_{\ell=1}^N| x_\ell -\overline{x}|^2)^m.
\end{array}\ee
As concerns the convergence of
\be\label{3.56}\int_s^T\eta_r^{4,N}dr=\int_s^T(
\frac{1}{N}\sum_{\ell=1}^N f(\overline{X}_r^{0},
\overline{X}_r^{\ell},\overline{Y}_r,\overline{Z}_r^{0},
0)-E[f(\overline{X}_r^0, \overline{X}_r^1,\overline{Y}_r,
\overline{Z}_r^0,0)|{\cal F}^{W^0}_T])dr,\ee
the same argument as that for the proof of Lemma
\ref{lemma_zeta} can be used. Indeed, recalling that the
function $f$ is bounded and that the processes $\overline{X}^0,\ \overline{Y}$ and $\overline{Z}^0$ are $\mathbb{F}^{W^0}$-adapted, we can use the fact that the
processes $\overline{X}^\ell,\, \ell\ge 1,$ are conditionally
i.i.d., knowing ${\cal F}^{W^0}_T,$ in order to conclude that
\be\label{3.57} E[\int_t^T|\eta_r^{4,N}|^{2m}dr]
\le\frac{C_m}{N^m},\ N\ge 1.\ee

\medskip

\noindent The above estimates allow now to get by a BSDE
standard argument that
\be\label{3.58}\begin{array}{lll}& & E[\sup_{s\in[t,T]}|
Y_s^N-\overline{Y}_s|^{2m}+(\int_t^T|Z_s^{0,N} -
\overline{Z}^0_s |^2ds+\sum_{\ell=1}^N\int_t^T|Z_s^{\ell,N}
|^2ds)^m]\\
& & \displaystyle\le C_m(\frac{1}{N}+\frac{1}{N}
\sum_{\ell=1}^N| x_\ell -\overline{x}|^2)^m.\\
\end{array}\ee
\noindent The proof is complete.\end{proof}

\subsection{Appendix 3}

\begin{proof} (of Lemma \ref{lemma4.7}). For convenience we omit the dependence on $\xi^{(N)}$, we deduce from
the fact that $(\overline{u}^N,\ \overline{v}^{(N)})$\ is the saddle
point of $H_N$ and the definition of $(\tilde{u}^N, \tilde{v}^{(N)})$ at one hand
$$\begin{array}{lll}
& &H_N(\overline{u}^N,\overline{v}^{(N)}) \ge H_N(\tilde{u}^N,
\overline{v}^{(N)})\ge H_N(\tilde{u}^N,\tilde{v}^{(N)}(\tilde{u}^N))\\
& &\ge H_N(\overline{u}^N,\tilde{v}^{(N)}(\overline{u}^N))\ge
H_N(\overline{u}^N,\overline{v}^{(N)}),\\
\end{array}$$

\noindent i.e., for all $u\in U,$
$$H_N(\tilde{u}^N,\overline{v}^{(N)})=
H_N(\overline{u}^N,\overline{v}^{(N)})\ge H_N(u,\overline{v}^{(N)}).$$

\noindent Consequently, $\tilde{u}^N,\ \overline{u}^N\in U$ are
both maximum points of the strict concave function $H_N(.,
\overline{v}^{(N)})$ (see (\ref{4.2})). This implies $\tilde{u}^N
=\overline{u}^N.$ On the other hand, using this equality we
have
$$H_N(\overline{u}^N,\tilde{v}^{(N)}(\tilde{u}^N))=H_N
(\tilde{u}^N,\tilde{v}^{(N)}(\tilde{u}^N))\le H_N(\tilde{u}^N,v)=
H_N(\overline{u}^N,v),$$

\noindent i.e., $\tilde{v}^{(N)}(\tilde{u}^N)\in V^N$ is a minimum
point of the strict convex function $H_N(\overline{u}^N,.)$
(see (\ref{4.2}) again), and it follows that $\tilde{v}^{(N)}(\tilde{u}^N)=\overline{v}^{(N)}.$
\end{proof}

We continue with the proof of Lemma \ref{lemma4.8}.
\begin{proof} (of Lemma \ref{lemma4.8}). Let $\xi^{(N)}=(x^{(N)}
=(x_0,\dots,x_N),y,z^{(N)}=(z_0,\dots,z_N)),\ \xi^{(N)}_\ell=
(x_0,x_\ell,y,z^{(N)})$ and, given an arbitrary $u\in U$,
$\tilde{v}_{N,\ell}:=\tilde{v}_N(\xi^{(N)}_\ell,u),\, 1\le\ell\le
N.$ Since the function $v^{(N)}\rightarrow H_N(\xi^{(N)},u,
v^{(N)})$ achieves its minimum at $\tilde{v}^{(N)}(\xi^{(N)},u)$,
we have, for all $1\le\ell\le N$,
\be\label{4.8}\begin{array}{lll}& &0=N(D_{v_\ell}H_N)
(\xi^{(N)},u,\tilde{v}^{(N)}(\xi^{(N)},u))\\
& &=(D_{v_\ell}f)(x_0,x_\ell,y,z_0,z_\ell,u,\tilde{v}_{N,\ell})
+\varepsilon_N\sum^N_{i=1}b_1(x_0,x_\ell,z_i)z_i.\end{array}\ee
Thus,
\be\label{4.9}\begin{array}{lll}
0&=&\((D_{v_\ell}f)(x_0,x_\ell,y,z_0,z_\ell,u,\tilde{v}_{N,
\ell})-(D_{v_\ell}f)(x_0,x_\ell,y,z_0,z_\ell,u,0),\ \tilde{v}_{N,
\ell}\)\\
& &+\left((D_{v_\ell}f)(x_0,x_\ell,y,z_0,z_\ell,u,0)+
\varepsilon_N\sum^N_{i=1}b_1(x_0,x_\ell,z_i)z_i,\ \tilde{v}_{N,
\ell}\right),
\end{array}\ee
and using that
\be\label{4.14}|\sum^N_{i=1}b_1(x_0,x_\ell,z_i)z_i|\leq \min
\{C\sum^N_{i=1}|z_i|, CN\},\ \mbox{for all}\ (x_0, x_\ell, z_i)\in\mathbb{R}^3,\ee
we deduce
\be\label{4.16}|\tilde{v}_N(\xi^{(N)}_\ell,u)|\leq C+
\frac{\mu}{\lambda}|u|+C\varepsilon_N\min\{\sum^N_{i=1}|z_i|,N\}.\ee
On the other hand, since $u\rightarrow H_N(\xi^{(N)},u,
\overline{v}^{(N)}(\xi^{(N)}))$ achieves a maximum at $\tilde{u}^N
(\xi^{(N)})(=\overline{u}^N(\xi^{(N)}))$, with the notation
$\xi_\ell=(x_0,x_\ell,y,z_0,z_\ell)$ it holds
\be\label{4.17}\begin{array}{lll}
 0&=&\frac{1}{N}\sum^N_{\ell=1}\((D_u f)(\xi_\ell,\tilde{u}^N
 (\xi^{(N)}),\bar{v}^N_\ell(\xi^{(N)}))-(D_u f)(\xi_\ell,0,
 \bar{v}^N_\ell(\xi^{(N)})), \tilde{u}^N(\xi^{(N)})\)\\
& &\ +\(\frac{1}{N}\sum^N_{\ell=1}(D_u f)(\xi_\ell,0,\bar{v}^N_\ell
(\xi^{(N)}))+\frac{1}{N}\sum^N_{\ell=1}b_0(x_0,x_\ell,z_0)z_0,
\tilde{u}^N(\xi^{(N)})\)\\
&\leq& -\lambda |\tilde{u}^N(\xi^{(N)})|^2+(C+\frac{\mu}{N}
\sum^N_{\ell=1}|\bar{v}^N_\ell(\xi^{(N)})|)|\tilde{u}^N(\xi^{(N)})|.
\end{array}\ee
Hence,
\be\label{4.18}\begin{array}{rcl}|\tilde{u}^N(\xi^{(N)})| &\leq &
C+\frac{\mu}{\lambda}\frac{1}{N}\sum^N_{\ell=1}|\bar{v}^N_\ell
(\xi^{(N)})|\\
&=&C+\frac{\mu}{\lambda}\frac{1}{N}\sum^N_{\ell=1}|\tilde{v}_N
(\xi^{(N)}_\ell,\tilde{u}^N(\xi^{(N)}))|\\
&\leq &C+\frac{\mu}{\lambda}\frac{1}{N}\sum^N_{\ell=1}(C+\frac{\mu}
{\lambda}|\tilde{u}^N(\xi^{(N)})|)
+\varepsilon_N NC \\
&\leq &C+C\varepsilon_N N+(\frac{\mu}{\lambda})^2|\tilde{u}^N
(\xi^{(N)})|.
\end{array}\ee
From $\mu<\lambda$, it follows that
\be\label{4.19}|\tilde{u}^N(\xi^{(N)})|\leq C(1+\varepsilon_N N),
\ee
and since $\bar{v}^{(N)}(\xi^{(N)})=(\tilde{v}_N(\xi^{(N)}_l,
\tilde{u}^N(\xi^{(N)})))_{1\leq l\leq N},$ we conclude that
\be\begin{array}{rcl}
& &|\bar{v}^N_\ell(\xi^{(N)})|\leq C(1+\varepsilon_N N),\ \ 1
\leq \ell\leq N,\ \mbox{i.e.,}\\
& &|\bar{v}^{(N)}(\xi^{(N)})|\leq \sqrt{N}C(1+\varepsilon_N N).
\end{array}\ee
On the other hand, from ($\ref{4.17}$),
$$\displaystyle 0\leq -\lambda|\tilde{u}^N(\xi^{(N)})
|^2+(C+\frac{\mu}{N}\sum^N_{\ell=1}|\tilde{v}_N(\xi^{(N)}_\ell,
\tilde{u}^N(\xi^{(N)}))|)|\tilde{u}^N(\xi^{(N)})|,$$

\noindent we see that
$$\displaystyle |\tilde{u}^N(\xi^{(N)})|\leq
C+\frac{\mu}{\lambda}(C+\frac{\mu}{\lambda}|\tilde{u}^N(\xi^{(N)})|+C\varepsilon_N\sum^N_{i=1}|z_i|),$$

\smallskip

\noindent from where, since $0<\mu<\lambda$, we get that
\be\label{4.22}|\tilde{u}^N(\xi^{(N)})|\leq C(1+
\varepsilon_N\sum^N_{i=1}|z_i|).\ee
Finally, putting (\ref{4.16}), (\ref{4.19}) and (\ref{4.22})
together, we complete the proof of the growth conditions (i) and (ii)
for the saddle point control.

It remains to prove the Lipschitz properties (iii) and (iv).
For this end, we observe that, since
$$ 0=(D_{v_\ell}f)(\xi_\ell,u,
\tilde{v}_N(\xi^{(N)}_\ell,u))+\varepsilon_N\sum^N_{i=1}
b_1(x_0,x_\ell,z_i)z_i, \quad 1\leq \ell\leq N,$$

\noindent we get from Assumption Ai), for all $u,\ \tilde{u}\in U,$
$$\displaystyle\begin{array}{rcl}
0&=&\((D_{v_\ell}f)(\xi_\ell,u,\tilde{v}_N(\xi^{(N)}_\ell,u))
-(D_{v_\ell}f)(\tilde{\xi}_\ell,\tilde{u},\tilde{v}_N
(\tilde{\xi}^{(N)}_\ell,\tilde{u})),\tilde{v}_N(
\xi^{(N)}_\ell,u)-\tilde{v}_N(\tilde{\xi}^{(N)}_\ell,\tilde{u})\)\\
& &+\varepsilon_N \sum^N_{i=1}\(b_1(x_0,x_\ell,z_i)z_i-b_1
(\tilde{x}_0,\tilde{x}_\ell,\tilde{z}_i)\tilde{z}_i,\ \tilde{v}_N(\xi^{(N)}_\ell,u)-\tilde{v}_N(\tilde{\xi}^{(N)}_\ell,\tilde{u})\)\\
&\geq&  \lambda
|\tilde{v}_N(\xi^{(N)}_\ell,u)-\tilde{v}_N(
\tilde{\xi}^{(N)}_\ell,\tilde{u})|^2\\
& &-|(D_{v_l}f)(\xi_\ell,u,\tilde{v}_N(\tilde{\xi}^{(N)}_\ell,
\tilde{u}))-(D_{v_\ell}f)(\tilde{\xi}_\ell,\tilde{u},
\tilde{v}_N(\tilde{\xi}^{(N)}_\ell,\tilde{u}))|
|\tilde{v}_N(\xi^{(N)}_\ell,u)-\tilde{v}_N(
\tilde{\xi}^{(N)}_\ell,\tilde{u})|\\
& &-\varepsilon_N\sum^N_{i=1}|b_1(x_0,x_\ell,z_i)z_i-b_1
(\tilde{x}_0,\tilde{x}_\ell,\tilde{z}_i)\tilde{z}_i|
|\tilde{v}_N(\xi^{(N)}_\ell,u)-\tilde{v}_N
(\tilde{\xi}^{(N)}_\ell,\tilde{u})|\\
&\geq & \lambda
|\tilde{v}_N(\xi^{(N)}_\ell,u)-\tilde{v}_N
(\tilde{\xi}^{(N)}_\ell,\tilde{u})|^2-(C|\eta_\ell
-\tilde{\eta}_\ell|+\mu
|u-\tilde{u}|)|\tilde{v}_N(\xi^{(N)}_\ell,u)-\tilde{v}_N
(\tilde{\xi}^{(N)}_\ell,\tilde{u})|\\
& &-\varepsilon_N \big\{N(|x_0-\tilde{x}_0|+|x_\ell-
\tilde{x}_\ell|)+\sum^N_{i=1}|z_i-\tilde{z}_i|\big\}|
\tilde{v}_N(\xi^{(N)}_\ell,u)-\tilde{v}_N
(\tilde{\xi}^{(N)}_\ell,\tilde{u})|,
\end{array}$$

\noindent and, consequently,
\be\label{4.26}\begin{array}{rcl}
|\tilde{v}_N(\xi^{(N)}_\ell,u)-\tilde{v}_N
(\tilde{\xi}^{(N)}_\ell,\tilde{u})|&\leq &
C|\eta_\ell-\tilde{\eta}_\ell|+\frac{\mu}{\lambda}|u-\tilde{u}|+
C\varepsilon_N N(|x_0-\tilde{x}_0|+|x_\ell-\tilde{x}_\ell|)\\
& &\displaystyle +C\varepsilon_N\sum^N_{i=1}|z_i-\tilde{z}_i|,
\end{array}\ee
where $\eta_\ell:=(x_0, x_\ell, y, z_0),\ \tilde{\eta}_\ell:=(\tilde{x}_0, \tilde{x}_\ell, \tilde{y}, \tilde{z}_0)$. With a similar argument, using again Assumption Ai) but exploiting now the strict concavity of $f$ in $u$, we deduce
from
$$\displaystyle0=\frac{1}{N}\sum^N_{\ell=1}(D_uf)
(\xi_\ell,\bar{u}^N(\xi^{(N)}),\bar{v}^N_\ell(\xi^{(N)}))+
\frac{1}{N}\sum^N_{\ell=1}b_0(x_0,x_\ell,z_0)z_0$$

\noindent that
\be\label{4.29}\begin{array}{rcl}|\bar{u}^N(\xi^{(N)})-
\bar{u}^N(\tilde{\xi}^{(N)})|&\leq&
\frac{\mu}{N\lambda}\sum^N_{\ell=1}|
\bar{v}^N_\ell(\xi^{(N)})-\bar{v}^N_\ell(\tilde{\xi}^{(N)})|\\
& &+C\big\{|x_0-\tilde{x}_0|+|y-\tilde{y}|+|z_0-\tilde{z}_0|+
\frac{1}{N}\sum^N_{\ell=1}(|x_\ell-\tilde{x}_\ell|
+|z_\ell-\tilde{z}_\ell|)\big\}.\end{array}\ee Hence, recalling that
$\bar{v}_\ell^N(\xi^{(N)})=\tilde{v}_N (\xi^{(N)}_\ell,\bar{u}^N(\xi^{(N)}))$
and using (\ref{4.26}) as well as $\mu/\lambda<1$, we obtain
$$\begin{array}{rcl}& &|\bar{u}^N(\xi^{(N)})-
\bar{u}^N(\tilde{\xi}^{(N)})|\\
&\leq& C(|x_0-\tilde{x}_0|+|y-\tilde{y}|+|z_0-\tilde{z}_0|)+
C(1+\varepsilon_N N)\{|x_0-\tilde{x}_0|+\frac{1}{N}
\sum^N_{\ell=1}(|x_\ell-\tilde{x}_\ell|+|z_\ell-
\tilde{z}_\ell|)\},
\end{array}$$

\smallskip

\noindent and combining the latter result with (\ref{4.26})
we have
\be\label{4.32}\begin{array}{rcl}& &|\bar{v}^N_\ell
(\xi^{(N)})-\bar{v}^N_\ell(\tilde{\xi}^{(N)})|=
|\tilde{v}_N(\xi^{(N)}_\ell,\bar{u}^N(\xi^{(N)}))-
\tilde{v}_N(\tilde{\xi}^{(N)}_\ell,\bar{u}^N
(\tilde{\xi}^{(N)}))|\\
&\leq& C|\eta_\ell-\tilde{\eta}_\ell|+\frac{\mu}{\lambda}
|\bar{u}^N(\xi^{(N)})-\bar{u}^N(\tilde{\xi}^{(N)})|+
C\varepsilon_N N(|x_0-\tilde{x}_0|+|x_\ell-
\tilde{x}_\ell|)+C\varepsilon_N\sum^N_{i=1}|z_i-
\tilde{z}_i|\\
&\leq& C|\eta_\ell-\tilde{\eta}_\ell|+C(1+\varepsilon_N N)
(|x_0-\tilde{x}_0|+|x_\ell-\tilde{x}_\ell|+\frac{1}{N}
\sum^N_{i=1}(|x_i-\tilde{x}_i|+|z_i-\tilde{z}_i|)).
\end{array}\ee
With the relations (\ref{4.29}) and (\ref{4.32}) we have
gotten the Lipschitz property stated for the saddle point
control in the lemma. The proof is complete.\end{proof}

\subsection{Appendix 4}

\begin{proof} (of Lemma \ref{lemma4.4}). The statement 1) for $\overline{H}$\ is a direct
consequence of the estimates for $f$ in (\ref{4.66}) and for $\overline{v}$ (see Lemma \ref{lemma4.3}). In order to prove statement 2), it suffices to consider
the function
$$\overline{H}_{0}(s,\xi,u):=E[f_{\bar{v}}
(\xi,\overline{X}_{s}^{1},u)|{\cal F}_T^{W^0}],\ \ (s,\xi,u)\in [0,T]\times\mathbb{R}^3\times U,$$

\smallskip

\noindent  and to determine its derivative with
respect to $u$. For this note that, for all $h\in U$,
$$\begin{array}{rcl}
&&\overline{H}_{0}(s,\xi,u+h)-\overline{H}_{0}(s,
\xi,u)=E[\int_0^1((D_{u}f_{\bar{v}})(\xi,
\overline{X}_{s}^{1},u+\delta h),h){\rm d}\delta
|{\cal F}_T^{W^0}]\\
&=&\(E[(D_{u}f)(\xi,\overline{X}_{s}^{1},u,\bar{v}
(\xi,\overline{X}_{s}^{1},u))|{\cal F}_T^{W^0}], h\)\\
& &+E[\int_0^1{\((D_{u}f)(\xi,\overline{X}_{s}^{1},
u+\delta h,{\bar{v}}(\xi,\overline{X}_{s}^{1},u
+\delta h))-(D_{u}f)(\xi,\overline{X}_{s}^{1},u,
{\bar{v}}(\xi,\overline{X}_{s}^{1},u)), h\)}{\rm d}
\delta |{\cal F}_T^{W^0}]\\
&=&\(E[(D_{u}f)(\xi,\overline{X}_{s}^{1},u,\bar{v}
(\xi,\overline{X}_{s}^{1},u))|{\cal F}_T^{W^0}], h\)
+\rho(h),
\end{array}$$

\noindent where, due to Lemma \ref{lemma4.3},
$$ |\rho(h)|\leq |h|E[\int_0^1{(C|h|
+\mu|\bar{v}(\xi,\overline{X}_{s}^{1},u+\delta
h)-\bar{v}(\xi,\overline{X}_{s}^{1},u)|)}{\rm d}
\delta|{\cal F}_T^{W^0}]\leq C|h|^2.$$

\noindent Hence,
$$ (D_{u}\overline{H}_{0})(s,\xi,u)
=E[(D_{u}f)(\xi,\overline{X}_{s}^{1},u,\bar{v}
(\xi,\overline{X}_{s}^{1},u))|{\cal F}_T^{W^0}],$$

\smallskip

\noindent and statement 2) of the lemma follows.

Now by using this relation and the assumptions on $f$ and the
Lipschitz property of $\overline{v}$ (Lemma
\ref{lemma4.3}), we get
$$\begin{array}{rcl}
 &&\((D_{u}\overline{H})(s,\xi,u)-(D_{u}
 \overline{H})(s,\xi,u'), u-u'\)=
 \((D_{u}\overline{H}_{0})(s,\xi,u)-(D_{u}
 \overline{H}_{0})(s,\xi,u'), u-u'\)\\
&=&E[\((D_{u}f)(\xi,\overline{X}_{s}^{1},u,
\bar{v}(\xi,\overline{X}_{s}^{1},u))-
(D_{u}f)(\xi,\overline{X}_{s}^{1},u',\bar{v}
(\xi,\overline{X}_{s}^{1},u)), u-u'\)|{\cal F}_T^{W^0}]\\
& &+E[\((D_{u}f)(\xi,\overline{X}_{s}^{1},u',
\bar{v}(\xi,\overline{X}_{s}^{1},u))-
(D_{u}f)(\xi,\overline{X}_{s}^{1},u',
\bar{v}(\xi,\overline{X}_{s}^{1},u')), u-u'\)|{\cal F}_T^{W^0}]\\
&\leq & -\lambda|u-u'|^2+\mu E[|\bar{v}(\xi,
\overline{X}_{s}^{1},u)-\bar{v}(\xi,
\overline{X}_{s}^{1},u')||{\cal F}_T^{W^0}]
|u-u'|\\
&\leq & -\lambda|u-u'|^2+\frac{\mu^2}{\lambda} |u-u'|^2=
-(\lambda-\frac{\mu^2}{\lambda})|u-u'|^2.
\end{array}$$

\noindent The proof is complete now.
\end{proof}

We continue with the proof of Lemma \ref{lemma4.61}.

\begin{proof} (of Lemma \ref{lemma4.61}). From the
 definition of $\overline{u}(s,\xi)$ as the maximum
 point of the $C^1$-function $\overline{H}(s,\xi,.):
 U\rightarrow\mathbb{R}$ it follows that, for all
$(s,\xi)\in[t,T]\times \mathbb{R}^3,$
$$\begin{array}{rcl}
0&=&\((D_{u}\overline{H})(s,\xi,\bar{u}(s,\xi)),
\bar{u}(s,\xi)\)\\
&=&\((D_{u}\overline{H})(s,\xi,\bar{u}(s,\xi))
-(D_{u}\overline{H})(s,\xi,0), \bar{u}(s,\xi)\)\\
& &+\(E[(D_{u}f)(\xi,\overline{X}_{s}^{1},0,\bar{v}
(\xi,\overline{X}_{s}^{1},0))|{\cal{F}}_T^{W^0}]+E[b_{0}(x_{0},
\overline{X}_{s}^{1},z_{0})z_{0}|{\cal{F}}_T^{W^0}], \bar{u}(s,\xi)\)\\
&\leq& -(\lambda-\frac{\mu^2}{\lambda})|\bar{u}(s,\xi)|^2+C(1+E[\mu
|\bar{v}(\xi,\overline{X}_{s}^{1},0)|])|\bar{u}(s,
\xi)|\\
 &\leq& -(\lambda-\frac{\mu^2}{\lambda})|\bar{u}(s,\xi)|^2+C
|\bar{u}(s,\xi)|,
\end{array}$$

\noindent (see Lemmas \ref{lemma4.4} and
\ref{lemma4.3}), i.e.,
$$|\bar{u}(s,\xi)|\leq C,\ (s,\ \xi)\in[t,T]\times \mathbb{R}^3.$$

\noindent For the same reason, using Lemmas \ref{lemma4.4} and \ref{lemma4.3} again,
we also have
$$\begin{array}{rcl}
0&=&\((D_{u}\overline{H})(s,\xi,\bar{u}
(s,\xi))-(D_{u}\overline{H})(s,\xi',\bar{u}
(s,\xi')),\ \bar{u}(s,\xi)-\bar{u}(s,\xi')\)\\
&\leq&-(\lambda-\frac{\mu^2}{\lambda})|\bar{u}(s,\xi)-\bar{u}(s,
\xi')|^2+\((D_{u}\overline{H})(s,\xi,\bar{u}(s,
\xi'))-(D_{u}\overline{H})(s,\xi',\bar{u}(s,\xi')),\ \bar{u}(s,\xi)-\bar{u}(s,\xi')\)\\
&\leq&-(\lambda-\frac{\mu^2}{\lambda})|\bar{u}(s,\xi)-\bar{u}(s,\xi')|^2+\{ E[|(D_{u}f)(\xi,\overline{X}_{s}^{1},
\bar{u}(s,\xi'),\bar{v}(\xi,\overline{X}_{s}^{1},
\bar{u}(s,\xi')))\\
& &-(D_{u}f)(\xi',\overline{X}_{s}^{1},\bar{u}
(s,\xi'),\bar{v} (\xi',\overline{X}_{s}^{1},\bar{u}
(s,\xi')))||{\cal F}_T^{W^0}]\\
& &+E[|b_{0}(x_{0},\overline{X}_{s}^{1},z_{0})
z_{0}-b_{0}(x_{0}', \overline{X}_{s}^{1},z_{0}')
z_{0}'||{\cal F}_T^{W^0}]\}|\bar{u}(s,\xi)-
\bar{u}(s,\xi')|\\
&\leq&-(\lambda-\frac{\mu^2}{\lambda})|\bar{u}(s,\xi)-\bar{u}(s,
\xi')|^2\\
& &+C\big(|\xi-\xi'|+E[\mu|\bar{v}(\xi,
\overline{X}_{s}^{1},\bar{u}(s,\xi'))-\bar{v}(\xi',
\overline{X}_{s}^{1},\bar{u}(s,\xi'))|
|{\cal F}_T^{W^0}]\big)|\bar{u}(s,\xi)-\bar{u}
(s,\xi')|\\
&\leq&-\frac12(\lambda-\frac{\mu^2}{\lambda})|\bar{u}(s,\xi)-
\bar{u}(s,\xi')|^{2}+C|\xi-\xi'|^2.
\end{array}$$

\smallskip

\noindent This proves the Lipschitz continuity of
$\overline{u}(s,.),$ uniformly with respect to
$s\in[0,T].$
\end{proof}

\subsection{Appendix 5}

Let us begin with the proof of Lemma \ref{auxlemma11}.
\begin{proof} (of Lemma \ref{auxlemma11}).
\underline{Step 1.} Recalling the definition of
$\tilde{v}_N$ in (\ref{4.7a}) and that of $\overline{v}$
in (\ref{a11}), we see that
\be\label{4.88}
\tilde{v}_{N}(x_{0},x_{1},y,(z_{0},0,\dots,0),u)=
\overline{v}(x_{0},x_{1},y,z_{0},u).
\ee
On the other hand, from (\ref{4.26}),
$$|\tilde{v}_{N}(x_{0},x_{1},y,(z_{0},
z_{1},\dots,z_{N}),u)-\tilde{v}_{N}(x_{0},x_{1},y,
(z_{0},0,\dots,0),u)|\leq C\varepsilon_{N}\sum_{i=1}^N|z_{i}|,$$

\noindent for all $(x_{0},x_{1},y,z_{0},(z_{1},\dots,z_{N}),
u)\in\mathbb{R}^{N+4}\times U.$\ Consequently,
\be\label{4.91}
|\tilde{v}_{N}(x_{0},x_{1},y,z^{(N)},u)-\overline{v}
(x_{0},x_{1},y,z_{0},u)|\leq C\varepsilon_{N}\sum_{i=1}^N|z_{i}|.
\ee

\noindent\underline{Step 2}. The objective of this step
is to estimate the difference between the controls
$\overline{u}^N$ and $\overline{u}$.

Let us use here in our computations the notations
$\xi:=(x_{0},y,z_{0})$ and  $(\xi,x_{1}):=(x_{0},x_{1},
y,z_{0})$. Recall also the notations $\xi^{(N)}=(x^{(N)}, y, z^{(N)})$ and $\xi^{(N,0)}=(x^{(N)},y,(z^{0,N},0,\dots,
0)).$

From the definition of $\overline{u}$ (see (\ref{4.58}))
we get (see also (\ref{4.57-1}):
\be\label{4.93}\begin{array}{rcl}
0&=&(D_{u}\overline{H})(s,\xi,\bar{u}(s,\xi))\\
&=&E[(D_{u}f)(\xi,\overline{X}_{s}^1,\bar{u}(s,\xi),
\bar{v}(\xi,\overline{X}_{s}^1,\bar{u}(s,\xi)))
+b_{0}(x_{0},\overline{X}_{s}^1,z_{0})z_{0}|
{\cal F}_T^{W^0}].
\end{array}\ee
On the other hand, from the definition of $\overline{u}^N$
as one of the both saddle point feedback controls for
$H_N(\xi^{(N)},.,.)$, from Lemma \ref{lemma4.7} and from
(\ref{4.88}) we obtain
\be\label{4.94}\begin{array}{rcl}
0&=&(D_{u}H_{N})(\xi^{(N,0)},\overline{u}^{N}(\xi^{(N,0)}),
\overline{v}^{N}(\xi^{(N,0)}))\\
&=&\frac{1}{N}\sum_{\ell=1}^{N}(D_{u}f)(\xi,x_\ell,0,
\overline{u}^{N}(\xi^{(N,0)}),\tilde{v}_{N}(\xi,x_\ell,0,
\overline{u}^{N}(\xi^{(N,0)})))\\
& &+\frac{1}{N}\sum_{\ell=1}^{N}b_{0}(x_{0},x_{\ell},z_{0})
z_{0}\\
&=&\frac{1}{N}\sum_{\ell=1}^{N}(D_{u}f_{\overline{v}})
(\xi,x_{\ell},0,\overline{u}^{N}(\xi^{(N,0)}))+
\frac{1}{N}\sum_{\ell=1}^{N}b_{0}(x_{0},x_{\ell},z_{0})z_{0}.
\end{array}\ee
\noindent Let us use the notations
$\overline{\Theta}_s=(X^0_s,\overline{Y}_s,\overline{Z}^0_s,
0)$, $\overline{\Theta}^N_s=(X^0_s,\overline{Y}^N_s,
\overline{Z}^{0,N}_s,0)$, and, with abusing notation we also
write $(\overline{\Theta}^N_s,X_s^\ell)=(X^0_s,X_s^\ell,
\overline{Y}^N_s, \overline{Z}^{0,N}_s,0).$ Moreover, let
$\Xi^{(N)}_s=(X^{(N)}_s,\overline{Y}_s^N,\overline{Z}_s^{N})$
and $\Xi^{(N,0)}_s=(X^{(N)}_s,\overline{Y}_s^N,
(\overline{Z}_s^{0,N}, 0,\dots,0)).$ Then subtracting
(\ref{4.93}) from (\ref{4.94}) yields
$$\begin{array}{lll}\displaystyle
0&=&\((\frac{1}{N}\sum_{\ell=1}^{N}(D_{u}f_{\overline{v}})
(\overline{\Theta}^{N},X_{s}^{\ell},
\overline{u}^{N}(\Xi^{(N,0)}_s))
+\frac{1}{N}\sum_{\ell=1}^{N}b_{0}(X_{s}^{0},X_{s}^{\ell},
\overline{Z}_{s}^{0,N})\overline{Z}_{s}^{0,N})\\
& &-E[(D_{u}f_{\overline{v}})(\overline{\Theta}_{s},
\overline{X}_{s}^1,\overline{u}(s,\overline{\Theta}_{s}))
+b_{0}(X_{s}^{0},\overline{X}_{s}^1,\overline{Z}_{s}^{0})
\overline{Z}_{s}^{0}|{\mathcal{F}}_{T}^{W^{0}}],\
\overline{u}^{N}
(\Xi^{(N,0)}_s)-\bar{u}(s,\overline{\Theta}_{s})\)\\
& =&\frac{1}{N}\sum_{\ell=1}^{N}\((D_{u}f_{\overline{v}})
(\overline{\Theta}_{s}^{N},X_{s}^{\ell},\overline{u}^{N}
(\Xi^{(N,0)}_s))-(D_{u}f_{\overline{v}})
(\overline{\Theta}_{s}^{N},X_{s}^{\ell},\overline{u}(s,
\overline{\Theta}_{s})),\overline{u}^{N}(\Xi^{(N,0)}_s)
-\overline{u}(s,\overline{\Theta}_{s})\) \hfill (=:I_1^N)\\
& &+\frac{1}{N}\sum_{\ell=1}^{N}\((D_{u}f_{\overline{v}})
(\overline{\Theta}_{s}^{N},X_{s}^{\ell},\overline{u}(s,
\overline{\Theta}_{s}))-(D_{u}f_{\overline{v}})
(\overline{\Theta}_{s},X_{s}^{\ell},\overline{u}(s,
\overline{\Theta}_{s})), \overline{u}^{N}(\Xi^{(N,0)}_s)
-\overline{u}(s,\overline{\Theta}_{s})\)\hfill (=: I_2^N)\\
& &+ \(\frac{1}{N}\sum_{\ell=1}^{N}b_{0}(X_{s}^{0},
X_{s}^{\ell},\overline{Z}_{s}^{0,N})\overline{Z}_{s}^{0,N}
-E[b_{0}(X_{s}^{0},\overline{X}_{s}^1,\overline{Z}_{s}^{0})
\overline{Z}_{s}^{0}|{\mathcal{F}}_{T}^{W^{0}}],
\overline{u}^{N}(\Xi^{(N,0)}_s)-\overline{u}(s,
\overline{\Theta}_{s})\)\hfill (=: I_3^N)\\
& &+\(\frac{1}{N}\sum_{\ell=1}^N(D_u f_{\overline{v}})
(\overline{\Theta}_s,X_s^\ell,\overline{u}(s,
\overline{\Theta}_s))-E[(D_u f_{\overline{v}})
(\overline{\Theta}_s,\overline{X}_s^1,\overline{u}(s,
\overline{\Theta}_s))|\mathcal{F}_T^{W^0}],
\overline{u}^N(\Xi^{(N,0)}_s)-\overline{u}(s,
\overline{\Theta}_s)\).\hfill (=: I_4^N)
\end{array}$$

\smallskip

\noindent Let us estimate the expressions $I_k^N,\
 1\le k\le 4.$ We begin with that of $I_1^N$.

\noindent $\bullet$ Estimate for $I_1^N$: By
standard estimates using our assumptions on $f$
we have
$$\begin{array}{rcl}
& &\((D_u f_{\overline{v}})(\xi,x_1,0,u)-
(D_u f_{\bar{v}})(\xi,x_1,0,u'), u-u'\)\\
&=&\((D_u f)(\xi,x_1,0,u,\overline{v}(\xi,x_1,0,u))
-(D_u f)(\xi,x_1,0,u',\overline{v}(\xi,x_1,0,u)),
u-u'\)\\
& &+\((D_u f)(\xi,x_1,0,u',\overline{v}(\xi,x_1,0,
u))-(D_u f)(\xi,x_1,0,u',\overline{v}(\xi,x_1,0,
u')),u-u'\)\\
&\leq&-\lambda|u-u'|^2+\mu|\overline{v}(\xi,x_1,
0,u)-\overline{v}(\xi,x_1,0,u')||u-u'|\\
&\leq&-(\lambda-\frac{\mu^2}{\lambda})|u-u'|^2
\end{array}$$
(recall that we have supposed that $\mu<\lambda)$.
Consequently,
\be\label{4.101}
I_1^N\leq -(\lambda-\frac{\mu^2}{\lambda})|\bar{u}^N(\Xi_s^{(N,0)})
-\overline{u}(s,\overline{\Theta}_s)|^2.\ee

\noindent $\bullet$ Estimate for $I_2^N:$\ From
(\ref{4.44}) and Assumption Aii) on $D_uf$ we have,
for arbitrarily small given $\delta>0$ and a constant
$C_\delta$ only depending on $\delta$,
\be\label{4.102}\begin{array}{rcl}
I_2^N &\leq& C|\overline{\Theta}_s^N-
\overline{\Theta}_s||\overline{u}^N(\Xi^{(N,0)}_s)
-\overline{u}(s,\overline{\Theta}_s)| \\
&\leq& C_\delta(|\overline{Y}_s^N-\overline{Y}_s|^2
+|\overline{Z}_s^{0,N}-\overline{Z}_s^0|^2)+\delta|
\overline{u}^N(\Xi^{(N,0)}_s)-\overline{u}(s,
\overline{\Theta}_s)|^2.
\end{array}\ee

\noindent $\bullet$ Estimate for $I_3^N:$ Using the
Lipschitz continuity of $z_0\mapsto b_0(x_0,x_1,z_0)
z_0,$ uniformly with respect to $(x_0,x_1)$, we obtain
that for any small $\delta>0$ there is a constant
$C_\delta$ such that
\be\label{4.103}\begin{array}{rcl}
I_3^N &\leq& C_\delta|\overline{Z}_s^{0,N}-
\overline{Z}^0_s|^2+\delta|\overline{u}^N(\Xi^{(N,
0)}_s)-\overline{u}(s,\overline{\Theta}_s)|^2 \\
& &+ C_\delta|\frac{1}{N}\sum_{\ell=1}^Nb_0(X_s^0,
X_s^\ell,\overline{Z}_s^0)\overline{Z}_s^0
-E[b_0(X_s^0,\overline{X}_s^1,\overline{Z}_s^0)
\overline{Z}_s^0|\mathcal{F}_T^{W^0}]|^2.
\end{array}\ee
Recalling that $|b_0(x_0,x_1,z_0)z_0-b_0(x_0,
x'_1,z_0)z_0|\leq C|x_1-x'_1|$, we can use
(\ref{3.47x}), in order to deduce that for all
$m\ge 1$, there is some constant $C_m$ such that
$$E[|\frac{1}{N}\sum_{\ell=1}^Nb_0(x_0,
X_s^\ell,z_0)z_0-E[b_0(x_0,\overline{X}_s^1,z_0)
z_0|{\cal F}_T^{W^0}]|^{2m}|{\cal F}_T^{W^0}]
\leq C_m(\frac{1}{N}+\frac{1}{N}\sum_{\ell=1}^N
|x_\ell-\overline{x}|^2)^m,$$

\smallskip

\noindent for all $N\ge 1,\ s\in[t,T],\ (x_0,z_0)
\in {\mathbb R}^2.$ Hence, as $X_s^0$ and
$\overline{Z}^0_s$ are ${\cal F}_s^{W^0}$-measurable,
\be\label{4.105}\begin{array}{l}
E[|\frac{1}{N}\sum_{\ell=1}^Nb_0
(X_s^0,X_s^\ell,\overline{Z}_s^0)\overline{Z}_s^0
-E[b_0(X_s^0,\overline{X}_s^1,\overline{Z}_s^0
)\overline{Z}_s^0|{\cal F}_T^{W^0}]|^{2m}
|{\cal F}_T^{W^0}]\\
\leq C_m(\frac{1}{N}+\frac{1}{N}\sum_{\ell=1}^N
|x_\ell-\overline{x}|^2)^m.
\end{array}\ee
\noindent $\bullet$ Estimate for $I_4^N:$ Obviously,
for all $\delta>0$ there is $C_\delta>0$ such that
\be\label{4.106}\begin{array}{rcl}
I_4^N&\leq&\delta|\overline{u}^N(\Xi^{(N,0)}_s)-
\overline{u}(s,\overline{\Theta}_s)|^2 \\
& & +C_\delta|\frac{1}{N}\sum_{\ell=1}^N(D_u
f_{\overline{v}})(\overline{\Theta}_s,X_s^\ell,
\overline{u}(s,\overline{\Theta}_s))
-E[(D_u f_{\overline{v}})(\overline{\Theta}_s,
\overline{X}_s^1,\overline{u}(s,
\overline{\Theta}_s))|\mathcal{F}_T^{W^0}]|^2.
\end{array}\ee
\noindent Noting that  $|(D_uf_{\overline{v}})
(\xi,x_1,0,u)-(D_uf_{\overline{v}})(\xi,x'_1,0,u)
|\leq C|x_1-x'_1|$ and observing that
$\overline{\Theta}_s$ is
$\mathcal{F}_s^{W^0}$-measurable, we obtain
similarly to the estimate for $I_3^N$ from
(\ref{3.47x}) that, for all $m\ge 1$ there is
some $C_m\in \mathbb{R}$ with
\be\label{4.110}\begin{array}{rcl}
&&E[|\frac{1}{N}\sum_{\ell=1}^N(D_u
f_{\overline{v}})(\overline{\Theta}_s,X_s^\ell,
\overline{u}(s,\overline{\Theta}_s))
-E[(D_u f_{\overline{v}})(\overline{\Theta}_s,
\overline{X}_s^1,\overline{u}(s,
\overline{\Theta}_s))|\mathcal{F}_T^{W^0}]
|^{2m}| \mathcal{F}_T^{W^0}]\\
&\leq&C_m(\frac{1}{N}+\frac{1}{N}
\sum_{\ell=1}^N|x_\ell-\overline{x}|^2)^m.
\end{array}\ee

\medskip

\noindent Now, choosing $\delta=\frac16
(\lambda-\frac{\mu^2}{\lambda})>0$, and combing the above
estimates for $I_k^N,\ 1\le k\le 4,$ we obtain
\be\label{4.111}\begin{array}{rcl}
0 &=& I_1^N+I_2^N+I_3^N+I_4^N \\
 &\leq&-3\delta|\overline{u}^N(\Xi^{(N,0)}_s)
 -\overline{u}(s,\overline{\Theta}_s)|^2\\
& &+C(|\overline{Y}_s^N-\overline{Y}_s|^2+|
\overline{Z}_s^{0,N}-\overline{Z}_s^0|^2)+
|R^N(s,\overline{\Theta}_s, \overline{u}
(s,\overline{\Theta}_s))|^2,
\end{array}\ee
\noindent where
\be\begin{array}{lcl}
R^N(s,\xi,u):&=&C\frac{1}{N}\sum_{\ell=1}^N|(D_u
f_{\overline{v}})(\xi,X_s^\ell,u)-E[(D_u
f_{\overline{v}})(\xi,\overline{X}^\ell_s,u)
|\mathcal{F}_T^{W^0}]|\\
& &+C\frac{1}{N}\sum_{\ell=1}^N|b_0(x_0,X_s^\ell,
z_0)z_0-E[b_0(x_0,\overline{X}^1_s,z_0)z_0
|\mathcal{F}_T^{W^0}]|,
\end{array}\ee
\noindent and
$$E[|R^N(s,\overline{\Theta}_s,
\overline{u}(s,\overline{\Theta}_s))|^{2m}|
{\cal F}_T^{W^0}]\le C_m(\frac{1}{N}+\frac{1}{N}
\sum_{\ell=1}^N|x_\ell-\overline{x}|^2)^m.$$

\smallskip

\noindent We recall that $\overline{u}^N_s
=\overline{u}^N(\Xi^{(N)}_s)\ (=\overline{u}^N
(X^{(N)}_s,\overline{Y}^N_s, \overline{Z}^{(N)})),$
and we put $\overline{u}^{N,0}_s=\overline{u}^N
(\Xi^{(N,0)}_s)\ (=\overline{u}^N(X^{(N)}_s,
\overline{Y}^N_s, \allowbreak (\overline{Z}^{0,N},
0,\dots,0)))$ and $\overline{u}_s=\overline{u}
(s,\overline{\Theta}_s)\ (=\overline{u}(s,X_s^0,
\overline{Y}_s,\overline{Z}^{0}_s))$.

Then, taking into account that by Lemma
\ref{lemma4.8} (iii) we have
$$|\bar{u}_s^N-\bar{u}_s^{N,0}|\leq
C(1+\varepsilon_N\cdot N)\frac{1}{N}\sum
\limits_{l=\ell}^N| \overline{Z}_s^{\ell,N}|,$$

\noindent we obtain from (\ref{4.111})
\be\label{4.115}\begin{array}{rcl}
& &|\overline{u}_s^N-\overline{u}_s|\leq|\overline{u}_s^N-\overline{u}_s^{N,0}|+|
\overline{u}_s^{N,0}-\overline{u}_s|\\
&\leq& C(|\overline{Y}_s^N-\overline{Y}_s|
+|\overline{Z}_s^{0,N}-\overline{Z}_s^0|+
(1+\varepsilon_N\cdot N)\frac{1}{N}
\sum_{\ell=1}^N|\overline{Z}_s^{\ell,N}|)+C|R^N(s,\overline{\Theta}_s,
\overline{u}_s)|,
\end{array}\ee
\noindent where
$$E[R^N(s,\overline{\Theta}_s,
\overline{u}_s)|^{2m}|{\cal F}_T^{W^0}]
\le C_m(\frac{1}{N}+\frac{1}{N}
\sum_{\ell=1}^N|x_\ell-\overline{x}|^2)^m.$$

\noindent\underline{Step 3}. Basing on the
results of above steps we prove now the
limit behavior of the controls process
$\overline{v}^{(N)}_s=(\overline{v}^{1,N}_s,
\dots,\overline{v}^{N,N}_s)$, as $N$ tends to
$+\infty.$ For this end, we recall from Lemma
\ref{lemma4.7} that
$$\overline{v}^{\ell,N}_s=
\overline{v}^N_\ell(X^{(N)}_s,\overline{Y}^N_s,
\overline{Z}^{(N)}_s)=\tilde{v}_N(X^0_s,
X^\ell_s,\overline{Y}_s^N,\overline{Z}^{(N)}_s,
\overline{u}_s),\, s\in[t,T],\ 1\le \ell\le N.$$

\noindent Thus, due to the estimates (\ref{4.91}) in Step 1
and (\ref{4.115}) in Step 2 as well as Lemma \ref{lemma4.3} we have
$$\begin{array}{rcl}
&&|\overline{v}^{\ell,N}_s-\overline{v}^\ell_s|\\
&\leq&|\widetilde{v}_N(X_s^0,X_s^\ell,
\overline{Y}_s^N,\overline{Z}_s^{(N)},
\overline{u}_s^N)- \overline{v}(X_s^0,X_s^\ell,
\overline{Y}_s^N,\overline{Z}_s^{0,N},
\overline{u}_s^N)|\\
&&+|\overline{v}(X_s^0,X_s^\ell,
\overline{Y}_s^N,\overline{Z}_s^{0,N},
\overline{u}_s^N)-\overline{v}(X_s^0,X_s^\ell,
\overline{Y}_s,\overline{Z}_s^0,
\overline{u}_s)|\\
&\leq&C\varepsilon_N\sum_{\ell=1}^N|\overline{Z}_s^{\ell,N}|
+C(|\overline{Y}_s^N-\overline{Y}_s|+
|\overline{Z}_s^{0,N}-\overline{Z}_s^0|+
|\overline{u}_s^N-\overline{u}_s|)\\
&\leq &C(|\overline{Y}_s^N-\overline{Y}_s|+
|\overline{Z}_s^{0,N}-\overline{Z}_s^0| +
(1+\varepsilon_N\cdot
N)\frac{1}{N}\sum_{\ell=1}^N|
\overline{Z}_s^{\ell,N}|)+C|R^N(s,
\overline{\Theta}_s,\overline{u}_s)|
\end{array}$$
\noindent We recall that an estimate for $R^N(s,
\overline{\Theta}_s,\overline{u}_s)$ is given
in (\ref{4.115}) in Step 2. The proof is complete
now.
\end{proof}
Let us come, finally, to the proof of Lemma
\ref{lemma4.138}.

\begin{proof} (of Lemma \ref{lemma4.138}).
Let us keep notations introduced in the
preceding proof. So we recall that, in particular,
$\overline{u}^{N,0}_s=\overline{u}^N(\Xi^{(N,0)}_s)\
(=\overline{u}^N(X^{(N)}_s,\overline{Y}^N_s,
\allowbreak (\overline{Z}^{0,N},0,\dots,0)))$, and
we introduce in the same sense the notation
$\overline{v}^{\ell,N,0}_s=\overline{v}^N_\ell
(X^{(N)}_s,\overline{Y}_s^N,(\overline{Z}^{0,N}_s,0,
\dots,0))$. Then, using our assumptions on $f$ (see
(\ref{4.66})) we get

\be\begin{array}{lll} &&|f(X_s^0,X_s^\ell,
\overline{Y}_s^N,\overline{Z}_s^{0,N},0,
\overline{u}_s^{N,0},\overline{v}^{\ell,N,0}_s)
-f(X_s^0,X_s^\ell,\overline{Y}_s,\overline{Z}_s^0,0,
\overline{u}_s,\overline{v}^\ell_s)|\\
&\leq&C(|\overline{Y}_s^N-\overline{Y}_s|+
|\overline{Z}_s^{0,N}-\overline{Z}_s^0|
+|\overline{u}_s^{N,0}-\overline{u}_s|
+|\overline{v}^{\ell,N,0}_s-\overline{v}^\ell_s|).\\
\end{array}\ee

\noindent Indeed, from (\ref{4.23}) of Lemma
\ref{lemma4.8} we know that the processes
$\overline{u}^{N,0}$ and $\overline{v}^{\ell,N,0}$
are bounded by a constant not depending on $N$.
On the other hand, from Lemma \ref{lemma4.61} we
have the boundedness of process $\overline{u}$,
and from the Lemmas \ref{lemma4.61} and \ref{lemma4.3}
we obtain also that of the process
$\overline{v}^\ell.$

Moreover, from Lemma \ref{auxlemma11} it follows
that
\be\label{4.125}
|\overline{u}_s^{N,0}-\overline{u}_s|+
|\overline{v}_{s}^{\ell,N,0}-\overline{v}^\ell_s|
\leq C(|\overline{Y}_s^N-\overline{Y}_s|+
|\overline{Z}_s^{0,N}-\overline{Z}_s^0|)+R^N_s.
\ee
Hence,
\be\label{4.126}\begin{array}{rcl}
&&|f(X_s^0,X_s^\ell,\overline{Y}_s^N,
\overline{Z}_s^{0,N},0,\overline{u}_s^{N,0},
\overline{v}^{\ell,N,0}_s)
-f(X_s^0,X_s^\ell,\overline{Y}_s,
\overline{Z}_s^0,0,\overline{u}_s,
\overline{v}_s^\ell)|\\
&\leq&C(|\overline{Y}_s^N-\overline{Y}_s|+
|\overline{Z}_s^{0,N}-\overline{Z}_s^0|)
+R^N_s,
\end{array}\ee
with $R^N_s$ satisfying the estimate given
in Lemma \ref{auxlemma11}. We also note that,
by using the assumptions on $b_0$, we can show
with similar arguments that
\be\label{4.127}\begin{array}{rcl}
&&|b_0(X_s^0,X_s^\ell,\overline{Z}_s^{0,N})
\overline{Z}_s^{0,N}\overline{u}_s^{N,0}
-b_0(X_s^0,X_s^\ell,\overline{Z}_s^0)
\overline{Z}_s^0\overline{u}_s| \\
&\leq&C|\overline{Z}_s^{0,N}-\overline{Z}_s^0|
+C|\overline{u}_s^{N,0}-\overline{u}_s|\\
&\leq& C(|\overline{Y}_s^N-\overline{Y}_s^0|
+|\overline{Z}_s^{0,N}-\overline{Z}_s^0|)
+R^N_s,\, s\in[t,T],\ N\ge 1.
\end{array}\ee
\noindent Consequently, recalling the notations
introduced for this proof and for the preceding
one, and by using the fact that $(\overline{u}^N
(\xi^{(N)}),\overline{v}^N(\xi^{(N)}))$ is a
saddle point of the Hamiltonian $H_N(\xi^{(N)},
\cdot,\cdot)$ and $(\overline{u}^N(\xi^{(N,0)}),
\overline{v}^N(\xi^{(N,0)}))$ is one of
$H_N(\xi^{(N,0)},\cdot,\cdot)(=H_N(x^{(N)},y,
(z^{0},0,\dots,0,.,.)))$, we observe
\be\label{4.129}\begin{array}{lll}
&&\overline{H}_N(\Xi_s^{(N)}):=H_N(\Xi_s^{(N)},
\overline{u}_s^N,\overline{v}_s^{N})
\leq H_N(\Xi_s^{(N)},\overline{u}_s^N,
\overline{v}_s^{N,0}) \\
&=&H_N(\Xi_s^{(N,0)},\overline{u}_s^N,
\overline{v}_s^{N,0})+\frac{1}{N}
\sum_{\ell=1}^N\big(f(X_s^0,X_s^\ell,
\overline{Y}_s^N,\overline{Z}_s^{0,N},
\overline{Z}_s^{\ell,N},\overline{u}_s^N,
\overline{v}^{\ell,N,0}_s) \\
& &-f(X_s^0,X_s^\ell,\overline{Y}_s^N,
\overline{Z}_s^{0,N},0,\overline{u}_s^N,
\overline{v}^{\ell,N,0}_s)\big)
+\varepsilon_N\frac{1}{N}\sum_{\ell=1}^N
(\sum_{i=1}^Nb_1(X_s^0,X_s^\ell,
\overline{Z}_s^{i,N})\overline{Z}_s^{i,N})
\overline{v}^{\ell,N,0}_s\\
&\leq&H_N(\Xi_s^{(N,0)},\overline{u}_s^N,
\overline{v}_s^{N,0})
+C\left(\frac{1}{N}+\varepsilon_N\right)
\sum_{\ell=1}^N|\overline{Z}_s^{\ell,N}|\\
&\leq&H_N(\Xi_s^{(N,0)},\overline{u}_s^{N,0},
\overline{v}_s^{N,0})
+C\left(\frac{1}{N}+\varepsilon_N\right)
\sum_{\ell=1}^N|\overline{Z}_s^{\ell,N}|.
\end{array}\ee

\noindent On the other hand, estimating
$\overline{H}_N(\Xi_s^{(N)})$ in the opposite
direction by using similar arguments as above
and the fact that $|\overline{v}^{\ell,N}_s|
\leq C(1+N\varepsilon_N)$ (see Lemma
\ref{lemma4.8}), we obtain
\be\label{4.133}\begin{array}{lll}
&&\overline{H}_N(\Xi_s^{(N)})=H_N(\Xi_s^{(N)},
\overline{u}_s^N,\overline{v}_s^{N})\geq H_N
(\Xi_s^{N},\overline{u}_s^{N,0},
\overline{v}_s^{N})\\
&=&H_N(\Xi_s^{(N,0)},\overline{u}_s^{N,0},
\overline{v}_s^N)+\frac{1}{N}\sum_{\ell=1}^N
(f(X_s^0,X_s^\ell, \overline{Y}_s^N,
\overline{Z}_s^{0,N},\overline{Z}_s^{\ell,N},
\overline{u}_s^{N,0},\overline{v}^{\ell,N}_s)\\
& &-f(X_s^0,X_s^\ell,\overline{Y}_s^N,
\overline{Z}_s^{0,N},0,\overline{u}_s^{N,0},
\overline{v}^{\ell,N}_s)
+\varepsilon_N\frac{1}{N}\sum_{\ell=1}^N
(\sum_{i=1}^Nb_1(X_s^0,X_s^\ell,
\overline{Z}_s^{i,N})\overline{Z}_s^{i,N})
\overline{v}^{\ell,N}_s \\
&\geq&H_N(\Xi_s^{(N,0)},\overline{u}_s^{N,0},
\overline{v}_s^{N})-C\frac{1}{N}\sum_{\ell=1}^N|
\overline{Z}_s^{\ell,N}|-C\varepsilon_N
(1+N\varepsilon_N)\sum_{i=1}^N|
\overline{Z}_s^{i,N}|\\
&\geq&H_N(\Xi_s^{(N,0)},\overline{u}_s^{N,0},
\overline{v}_s^{N})-C(\frac{1}{N}+\varepsilon_N
+N\varepsilon_N^2)\sum_{\ell=1}^N|
\overline{Z}_s^{\ell,N}|\\
&\geq&H_N(\Xi_s^{(N,0)},\overline{u}_s^{N,0},
\overline{v}_s^{N,0})-C(\frac{1}{N}+
\varepsilon_N+N\varepsilon_N^2)\sum_{\ell=1}^N
|\overline{Z}_s^{\ell,N}|.
\end{array}\ee
\noindent By combining (\ref{4.129}) and
(\ref{4.133}) we obtain
\be\label{4.134}
|\overline{H}_N(\Xi_s^{(N)})-\overline{H}_N
(\Xi_s^{(N,0)})|\leq C(\frac{1}{N}+\varepsilon_N+
N\varepsilon_N^2) \sum_{\ell=1}^N|
\overline{Z}_s^{\ell,N}|.
\ee
\noindent Finally, from (\ref{4.126}) and
(\ref{4.127}) we see
\be\label{4.135}\begin{array}{lll}
&&|\overline{H}_N(\Xi_s^{(N,0)})-(\frac{1}{N}
\sum_{\ell=1}^Nf(X_s^0,X_s^\ell,\overline{Y}_s,
\overline{Z}_s^0,0,\overline{u}_s,
\overline{v}^\ell_s)+\frac{1}{N}\sum_{\ell=1}^Nb_0
(X_s^0,X_s^\ell,\overline{Z}_s^0)\overline{Z}_s^0
\overline{u}_s)|\\
&\leq&C(|\overline{Y}_s^N-\overline{Y}_s|+
|\overline{Z}_s^{0,N}-\overline{Z}_s^0|)+R^N_s,
\ \ s\in[t,T],\ N\geq 1.
\end{array}\ee
\noindent Recalling that the functions $x_\ell
\rightarrow f(x_0,x_l,y,z_0,u),$ and $x_l\rightarrow
b_0(x_0,x_l,z_0)z_0$ are
Lipschitz, uniformly with respect to $(x_0,y,z_0),$
and that $|\overline{u}_s|\leq C,$ we can apply Lemma
\ref{lemma3.47}, and we get for
$$\begin{array}{lll}
& &\overline{R}_1^N(s,x_0,y,z_0,u):=|\frac{1}{N}
\sum_{\ell=1}^Nf(x_0,X_s^\ell,y,z_0,0,u,\overline{v}^\ell_s)-
E[f(x_0,\overline{X}_s^1,y,z_0,0,u,\overline{v}^\ell_s)|
{\cal F}_T^{W^0}]|;\\
& &\overline{R}_2^{N}(s,x_0,z_0):=|\frac{1}{N}
\sum_{\ell=1}^Nb_0(x_0,X_s^\ell,z_0)z_0-E[b_0(x_0,
\overline{X}_s^1,z_0)z_0|{\cal F}_T^{W^0}]|,\\
& &\overline{R}^N(s,x_0,y,z_0,u):=\overline{R}_1^N
(s,x_0,y,z_0,u)+\overline{R}_2^N(s,x_0,z_0)u,\\
\end{array}$$

\noindent  the estimate
$$E[|\overline{R}^N(s,X_s^0,
\overline{Y}_s,\overline{Z}_s^0,\overline{u}_s)
|^{2m}|\mathcal{F}_T^{W^0}]\leq C_m(\frac{1}{N}
+\frac{1}{N}\sum_{\ell=1}^N|x_\ell-\overline{x}
|^2)^m,\mbox{ where }s\in[t,T],\ N\geq1,\ m\geq1.$$

\smallskip

\noindent The statement of the lemma follows now easily from the
latter estimates and (\ref{4.135}).
\end{proof}

\end{document}